\documentclass{amsart}
\usepackage[latin9]{inputenc}
\usepackage{ifsym}
\usepackage{amstext}
\usepackage{amsthm}
\usepackage{amssymb}

\makeatletter
\numberwithin{equation}{section}
\numberwithin{figure}{section}
\theoremstyle{plain}
\newtheorem{thm}{\protect\theoremname}
\theoremstyle{plain}
\newtheorem{lem}[thm]{\protect\lemmaname}
\theoremstyle{plain}
\newtheorem{prop}[thm]{\protect\propositionname}
\theoremstyle{definition}
\newtheorem{defn}[thm]{\protect\definitionname}
\theoremstyle{plain}
\newtheorem{cor}[thm]{\protect\corollaryname}

\@ifundefined{date}{}{\date{}}

\usepackage{amscd}
\usepackage{thmdefs}
\usepackage[dvips]{epsfig}

\input{tcilatex}
\theoremstyle{definition}
\theoremstyle{remark}
\numberwithin{equation}{section}

\input tcilatex

\makeatother

\providecommand{\corollaryname}{Corollary}
\providecommand{\definitionname}{Definition}
\providecommand{\lemmaname}{Lemma}
\providecommand{\propositionname}{Proposition}
\providecommand{\theoremname}{Theorem}

\begin{document}
\title{Operators Induced by Certain Hypercomplex Systems}
\author{Daniel Alpay and Ilwoo Cho}
\address{Chapman Univ., Dept. of Math., 1 University Dr., Orange, CA, 92866,
U. S. A. / St. Ambrose Univ., Dept. of Math. and Stat., 421 Ambrose
Hall, 518 W. Locust St., Davenport, Iowa, 52803, U. S. A.}
\email{alpay@chapman.edu / choilwoo@sau.edu}
\begin{abstract}
In this paper, we consider natural Hilbert-space representations $\left\{ \left(\mathbb{C}^{2},\pi_{t}\right)\right\} _{t\in\mathbb{R}}$
of the hypercomplex system $\left\{ \mathbb{H}_{t}\right\} _{t\in\mathbb{R}}$,
and study the realizations $\pi_{t}\left(h\right)$ of hypercomplex
numbers $h\in\mathbb{H}_{t}$, as $\left(2\times2\right)$-matrices
acting on $\mathbb{C}^{2}$, for an arbitrarily fixed scale $t\in\mathbb{R}$.
Algebraic, operator-theoretic, spectral-analytic, and free-probabilistic
properties of them are considered.
\end{abstract}

\maketitle

\section{Introduction}

In this paper, we study representations of the hypercomplex numbers
$\left(a,b\right)$ of complex numbers $a$ and $b$, constructing
a ring, 
\[
\mathbb{H}_{t}=\left(\mathbb{C}^{2},\;+,\;\cdot_{t}\right),
\]
scaled by a real number $t\in\mathbb{R}$, where ($+$) is the usual
vector addition on the 2-dimensional vector space $\mathbb{C}^{2}$,
and ($\cdot_{t}$) is the $t$-scaled vector-multiplication on $\mathbb{C}^{2}$,
defined by
\[
\left(a_{1},b_{1}\right)\cdot_{t}\left(a_{2},b_{2}\right)=\left(a_{1}a_{2}+tb_{1}\overline{b_{2}},\:a_{1}b_{2}+b_{1}\overline{a_{2}}\right),
\]
where $\overline{z}$ are the conjugates of $z$ in $\mathbb{C}$.

Motivated by the canonical Hilbert-space representation $\left(\mathbb{C}^{2},\pi\right)$
of the quaternions $\mathbb{H}$, introduced in {[}2{]}, {[}3{]} and
{[}19{]}, we consider the canonical representation,
\[
\Pi_{t}=\left(\mathbb{C}^{2},\;\pi_{t}\right),
\]
of the ring $\mathbb{H}_{t}$, and understand each element $h=\left(a,b\right)$
of $\mathbb{H}_{t}$ as its realization,
\[
\pi_{t}\left(h\right)\overset{\textrm{denote}}{=}[h]_{t}\overset{\textrm{def}}{=}\left(\begin{array}{cc}
a & tb\\
\overline{b} & \overline{a}
\end{array}\right)\;\mathrm{in\;}M_{2}\left(\mathbb{C}\right),
\]
where $M_{2}\left(\mathbb{C}\right)=B\left(\mathbb{C}^{2}\right)$
is the matricial algebra (or, the operator algebra acting on $\mathbb{C}^{2}$)
of all $\left(2\times2\right)$-matrices over $\mathbb{C}$ (respectively,
all bounded linear transformations, or simply operators on $\mathbb{C}^{2}$),
for each $t\in\mathbb{R}$. Under our setting, one can check that
the ring $\mathbb{H}_{-1}$ is nothing but the noncommutative field
$\mathbb{H}$ of all quaternions (e.g., {[}2{]}, {[}3{]} and {[}19{]}),
and the ring $\mathbb{H}_{1}$ is the ring of all bicomplex numbers
(e.g., {[}1{]}). 

The spectral-analytic, operator-theoretic (or, matrix-theoretic),
and free-probabilistic properties of $\mathbb{H}_{t}$ are considered
and characterized under the canonical representation $\Pi_{t}$. In
particular, certain decompositional properties on $\mathbb{H}_{t}$
are studied algebraically, and spectral-theoretically. And then, it
is considered how those properties affect the spectral-analytic, operator-theoretic,
and free-probabilistic properties of hypercomplex numbers of $\mathbb{H}_{t}$,
for $t\in\mathbb{R}$. 

\subsection{Motivation}

The quaternions $\Bbb{H}$ is an interesting object not only in pure
mathematics (e.g., {[}5{]}, {[}10{]}, {[}11{]}, {[}12{]}, {[}13{]}
{[}14{]}, {[}17{]}, {[}19{]}, {[}23{]}), but also in applied mathematics
(e.g., {[}4{]}, {[}7{]}, {[}15{]}, {[}16{]}, {[}20{]} and {[}21{]}).
Independently, spectral analysis on $\mathbb{H}$ is considered in
{[}2{]} and {[}3{]}, under representation, ``over $\mathbb{C}$,''
different from the usual quaternion-eigenvalue problems of quaternion-matrices
studied in {[}13{]}, {[}15{]} and 16{[}{]}.

Motivated by the generalized setting of the quaternions so-called
the split-quaternions of {[}1{]}, and by the main results of {[}2{]}
and {[}3{]}, we study a new type of hypercomplex numbers induced by
the pairs of $\mathbb{C}^{2}$. Especially, we construct a system
of the scaled hypercomplex rings $\left\{ \mathbb{H}_{t}\right\} _{t\in\mathbb{R}}$,
and study how the hypercomplex numbers act as $\left(2\times2\right)$-matrices
over $\mathbb{C}$ for given scales $t\in\mathbb{R}$, under our canonical
Hilbert-space representations $\left\{ \Pi_{t}=\left(\mathbb{C}^{2},\pi_{t}\right)\right\} _{t\in\mathbb{R}}$.
We are interested in algebraic, operator-theoretic, spectral-theoretic,
free-probabilistic properties of $\mathbb{H}_{t}$ under $\Pi_{t}$,
for $t\in\mathbb{R}$. Are they similar to those of the quaternions
$\mathbb{H}_{-1}=\mathbb{H}$, shown in {[}2{]} and {[}3{]}? The answers
are determined differently case-by-case, up to scales (See below).

\subsection{Overview}

In Section 2, we define our main objects, the scaled hypercomplex
rings $\left\{ \mathbb{H}_{t}\right\} _{t\in\mathbb{R}}$, and their
canonical Hilbert-space representations $\left\{ \Pi_{t}\right\} _{t\in\mathbb{R}}$.
We understand each hypercomplex number of $\mathbb{H}_{t}$ as an
operator, a $\left(2\times2\right)$-matrix over $\mathbb{C}$. We
concentrate on studying the invertibility on $\mathbb{H}_{t}$, for
an arbitrarily fixed scale $t$. It is shown that if $t<0$, then
$\mathbb{H}_{t}$ forms a noncommutative field like the quaternions
$\mathbb{H}=\mathbb{H}_{-1}$, however, if $t\geq0$, then it becomes
a ring with unity, which is not a noncommutative field.

In Section 3, the spectral theory on (the realizations of) $\mathbb{H}_{t}$
is studied over $\mathbb{C}$. After finding the spectra of hypercomplex
numbers, we define so-called the $t$-spectral forms whose main diagonal
entries are from the spectra, and off-diagonal entries are 0's. As
we have seen in {[}2{]} and {[}3{]}, such spectral forms are similar
to the realizations of quaternions of $\mathbb{H}_{-1}$. However,
if a scale $t\in\mathbb{R}\setminus\left\{ -1\right\} $ is arbitrary,
then such a similarity does not hold in general. We focus on studying
such a similarity in detail.

In Section 4, we briefly discuss about how the usual adjoint on $M_{2}\left(\mathbb{C}\right)$
acts on the sub-structure $\mathcal{H}_{2}^{t}$ of $\mathit{M_{2}\left(\mathbb{C}\right)}$,
consisting of all realizations of $\mathbb{H}_{t}$, for a scale $t\in\mathbb{R}$.
Different from the quaternionic case of {[}2{]} and {[}3{]}, in general,
the adjoints (conjugate-transposes) of many matrices of $\mathcal{H}_{2}^{t}$
are not contained in $\mathcal{H}_{2}^{t}$, especially, if $t\neq-1$.
It shows that a bigger, operator-algebraically-better $*$-algebraic
structure generated by $\mathcal{H}_{2}^{t}$ is needed in $M_{2}\left(\mathbb{C}\right)$,
to consider operator-theoretic, and free-probabilistic properties
on $\mathcal{H}_{2}^{t}$.

In the final Section 5, on the $C^{*}$-algebraic structure of Section
4, we study operator-theoretic, and free-probabilistic properties
up to the usual trace, and the normalized trace.

\section{The Scaled Hypercomplex Systems $\left\{ \mathbb{H}_{t}\right\} _{t\in\mathbb{R}}$}

In this section, we define a ring $\mathbb{H}_{t}$ of hypercomplex
numbers, and establish the corresponding canonical Hilbert-space representations
$\Pi_{t}$, for an arbitrary fixed scale $t\in\mathbb{R}$. Throughout
this section, we let 
\[
\mathbb{C}^{2}=\left\{ \left(a,b\right):a,b\in\mathbb{C}\right\} 
\]
be the Cartesian product of two copies of the complex field $\mathbb{C}$.
One may understand $\mathbb{C}^{2}$ as the usual 2-dimensional Hilbert
space equipped with its canonical orthonormal basis, $\left\{ \left(1,0\right),\;\left(0,1\right)\right\} .$

\subsection{A $t$-Scaled Hypercomplex Ring $\mathbb{H}_{t}$}

In this section, we fix an arbitrary real number $t$ in the real
field $\mathbb{R}$. On the vector space $\mathbb{C}^{2}$ (over $\mathbb{C}$),
define the $t$-scaled vector-multiplication ($\cdot_{t}$) by

\medskip{}

\hfill{}$\left(a_{1},b_{1}\right)\cdot_{t}\left(a_{2},b_{2}\right)\overset{\textrm{def}}{=}\left(a_{1}a_{2}+tb_{1}\overline{b_{2}},\:a_{1}b_{2}+b_{1}\overline{a_{2}}\right),$\hfill{}(2.1.1)

\medskip{}

\noindent for $\left(a_{l},b_{l}\right)\in\mathbb{C}^{2}$, for all
$l=1,2$, where $\overline{z}$ are the conjugates of $z$ in $\mathbb{C}$.
It is not difficult to check that such an operation ($\cdot_{t}$)
is closed on $\mathbb{C}^{2}$. Moreover, it satisfies that

\medskip{}

$\;\;\;\;\;\;$$\left(\left(a_{1},b_{1}\right)\cdot_{t}\left(a_{2},b_{2}\right)\right)\cdot_{t}\left(a_{3},b_{3}\right)$

\medskip{}

$\;\;\;\;\;\;\;\;\;\;\;\;$$=\left(a_{1}a_{2}+tb_{1}\overline{b_{2}},\:a_{1}b_{2}+b_{1}\overline{a_{2}}\right)\cdot_{t}\left(a_{3},b_{3}\right)$

\medskip{}

$\;\;\;\;\;\;\;\;\;\;\;\;$$=\left(a_{1}a_{2}a_{3}+t\left(b_{1}\overline{b_{2}}a_{3}+a_{1}b_{2}\overline{b_{3}}+b_{1}\overline{a_{2}}\overline{b_{3}}\right),\right.$

$\;\;\;\;\;\;\;\;\;\;\;\;\;\;\;\;\;\;\;\;$$\left.\;a_{1}a_{2}b_{3}+a_{1}b_{2}\overline{a_{3}}+b_{1}\overline{a_{2}}\overline{a_{3}}+tb_{1}\overline{b_{2}}b_{3}\right),$

\medskip{}

\noindent and

\medskip{}

$\;\;\;\;\;\;$$\left(a_{1},b_{1}\right)\cdot_{t}\left(\left(a_{2},b_{2}\right)\cdot_{t}\left(a_{3},b_{3}\right)\right)$

\medskip{}

$\;\;\;\;\;\;\;\;\;\;\;\;$$=\left(a_{1},b_{1}\right)\cdot_{t}\left(a_{2}a_{3}+tb_{2}\overline{b_{3}},\:a_{2}b_{3}+b_{2}\overline{a_{3}}\right)$

\medskip{}

$\;\;\;\;\;\;\;\;\;\;\;\;$$=\left(a_{1}\left(a_{2}a_{3}+tb_{2}\overline{b_{3}}\right)+tb_{1}\left(\overline{a_{2}}\overline{b_{3}}+\overline{b_{2}}a_{3}\right),\right.$

$\;\;\;\;\;\;\;\;\;\;\;\;\;\;\;\;\;\;\;\;$$\left.\;a_{1}\left(a_{2}b_{3}+b_{2}\overline{a_{3}}\right)+b_{1}\left(\overline{a_{2}}\overline{a_{3}}+t\overline{b_{2}}b_{3}\right)\right)$,

\medskip{}

\noindent implying the equality,\hfill{}(2.1.2)
\[
\left(\left(a_{1},b_{1}\right)\cdot_{t}\left(a_{2},b_{2}\right)\right)\cdot_{t}\left(a_{3},b_{3}\right)=\left(a_{1},b_{1}\right)\cdot_{t}\left(\left(a_{2},b_{2}\right)\cdot_{t}\left(a_{2},b_{3}\right)\right),
\]
in $\mathbb{C}^{2}$, for $\left(a_{l},b_{l}\right)\in\mathbb{C}^{2}$,
for all $l=1,2,3$.

Furthermore, if $\vartheta=\left(1,0\right)\in\mathbb{C}^{2}$, then

\medskip{}

\hfill{}$\vartheta\cdot_{t}\left(a,b\right)=\left(a,b\right)=\left(a,b\right)\cdot_{t}\vartheta$\hfill{}(2.1.3)

\medskip{}

\noindent by (2.1.1), for all $\left(a,b\right)\in\mathbb{C}^{2}$.

By (2.1.2) and (2.1.3), if
\[
\mathbb{C}^{2\times}=\mathbb{C}^{2}\setminus\left\{ \left(0,0\right)\right\} ,
\]
then the pair $\left(\mathbb{C}^{2\times},\:\cdot_{t}\right)$ forms
a monoid (i.e., semigroup with its identity $\left(1,0\right)$).
\begin{lem}
Let $\mathbb{C}^{2\times}=\mathbb{C}^{2}\setminus\left\{ \left(0,0\right)\right\} $,
and ($\cdot_{t}$) be the closed operation (2.1.1) on $\mathbb{C}^{2}$.
Then the algebraic structure $\left(\mathbb{C}^{2\times},\:\cdot_{t}\right)$
forms a monoid with its identity $\left(1,0\right)$.
\end{lem}

\begin{proof}
The proof is done by (2.1.2) and (2.1.3).
\end{proof}
Therefore, one can obtain the following ring structure.
\begin{prop}
The algebraic triple $\left(\mathbb{C}^{2},+,\cdot_{t}\right)$ forms
a unital ring with its unity (or, the multiplication-identity) $\left(1,0\right)$,
where ($+$) is the usual vector addition on $\mathbb{C}^{2}$, and
($\cdot_{t}$) is the vector multiplication (2.1.1).
\end{prop}

\begin{proof}
Clearly, the algebraic pair $\left(\mathbb{C}^{2},+\right)$ is an
abelian group under the usual addition ($+$) with its ($+$)-identity
$\left(0,0\right)$. While, by Lemma 1, the pair $\left(\mathbb{C}^{2\times},\cdot_{t}\right)$
forms a monoid (and hence, a semigroup). Observe now that

\medskip{}

$\;\;\;\;$$\left(a_{1},b_{1}\right)\cdot_{t}\left(\left(a_{2},b_{2}\right)+\left(a_{3},b_{3}\right)\right)=\left(a_{1},b_{1}\right)\cdot_{t}\left(a_{2}+a_{3},b_{2}+b_{3}\right)$

\medskip{}

$\;\;\;\;\;\;$$=\left(a_{1}\left(a_{2}+a_{3}\right)+tb_{1}\left(\overline{b_{2}}+\overline{b_{3}}\right),\:a_{1}\left(b_{2}+b_{3}\right)+b_{1}\left(\overline{a_{2}}+\overline{a_{3}}\right)\right)$

\medskip{}

$\;\;\;\;\;\;$$=\left(a_{1}a_{2}+a_{1}a_{3}+tb_{1}\overline{b_{2}}+tb_{1}\overline{b_{3}},\:a_{1}b_{2}+a_{1}b_{3}+b_{1}\overline{a_{2}}+b_{1}\overline{a_{3}}\right)$

\medskip{}

$\;\;\;\;\;\;$$=\left(a_{1}a_{2}+tb_{1}\overline{b_{2}},\:a_{1}b_{2}+b_{1}\overline{a_{2}}\right)+\left(a_{1}a_{3}+tb_{1}\overline{b_{3}},\:a_{1}b_{3}+b_{1}\overline{a_{3}}\right)$

\medskip{}

$\;\;\;\;\;\;$$=\left(a_{1},b_{1}\right)\cdot_{t}\left(a_{2},b_{2}\right)+\left(a_{1},b_{1}\right)\cdot_{t}\left(a_{3},b_{3}\right)$,

\medskip{}

\noindent and, similarly,\hfill{}(2.1.4)
\[
\left(\left(a_{1},b_{1}\right)+\left(a_{2},b_{2}\right)\right)\cdot_{t}\left(a_{3},b_{3}\right)=\left(a_{1},b_{1}\right)\cdot_{t}\left(a_{3},b_{3}\right)+\left(a_{2},b_{2}\right)\cdot_{t}\left(a_{3},b_{3}\right),
\]
in $\mathbb{C}^{2}$. So, the operations ($+$) and ($\cdot_{t}$)
are left-and-right distributive by (2.1.4).

Therefore, the algebraic triple $\left(\mathbb{C}^{2},+,\cdot_{t}\right)$
forms a unital ring with its unity $\left(1,0\right)$.
\end{proof}
The above proposition characterizes the algebraic structure of $\left(\mathbb{C}^{2},+,\cdot_{t}\right)$
as a well-defined unital ring for a fixed $t\in\mathbb{R}$. Remark
here that, since a scale $t$ is arbitrary in $\mathbb{R}$, in fact,
we obtain the unital rings $\left\{ \mathbb{H}_{t}\right\} _{t\in\mathbb{R}}$.
\begin{defn}
For a fixed $t\in\mathbb{R}$, the ring $\left(\mathbb{C}^{2},+,\cdot_{t}\right)$
is called the hypercomplex ring with its scale $t$ (in short, the
$t$-scaled hypercomplex ring). By $\mathbb{H}_{t}$, we denote the
$t$-scaled hypercomplex ring.
\end{defn}

\subsection{The Canonical Representation $\Pi_{t}=\left(\mathbb{C}^{2},\;\pi_{t}\right)$
of $\mathbb{H}_{t}$}

In this section, we fix $t\in\mathbb{R}$, and the corresponding $t$-scaled
hypercomplex ring,
\[
\mathbb{H}_{t}=\left(\mathbb{C}^{2},\:+,\:\cdot_{t}\right),
\]
where ($\cdot_{t}$) is the vector-multiplication (2.1.1). We consider
a natural finite-dimensional-Hilbert-space representation $\Pi_{t}$
of $\mathbb{H}_{t}$, and understand each hypercomplex number $h\in\mathbb{H}_{t}$
as an operator acting on a Hilbert space determined by $\Pi_{t}$.
In particular, as in the quaternionic case of {[}2{]}, {[}3{]} and
{[}19{]}, a 2-dimensional-Hilbert-space representation of the hypercomplex
ring $\mathbb{H}_{t}$ is established naturally.

Define now a morphism,
\[
\pi_{t}:\mathbb{H}_{t}\rightarrow B\left(\mathbb{C}^{2}\right)=M_{2}\left(\mathbb{C}\right),
\]
by\hfill{}(2.2.1)
\[
\pi_{t}\left(\left(a,b\right)\right)=\left(\begin{array}{cc}
a & tb\\
\overline{b} & \overline{a}
\end{array}\right),\;\forall\left(a,b\right)\in\mathbb{H}_{t},
\]
where $B\left(H\right)$ is the operator algebra consisting of all
bounded (or, continuous linear) operators on a Hilbert space $H$,
and $M_{k}\left(\mathbb{C}\right)$ is the matricial algebra of all
$\left(k\times k\right)$-matrices over $\mathbb{C}$, isomorphic
to $B\left(\mathbb{C}^{k}\right)$, for all $k\in\mathbb{N}$ (e.g.,
{[}8{]} and {[}9{]}).

By definition, the function $\pi_{t}$ of (2.2.1) is an injective
map from $\mathbb{H}_{t}$ into $M_{2}\left(\mathbb{C}\right)$. Indeed,
if
\[
\left(a_{1},b_{1}\right)\neq\left(a_{2},b_{2}\right)\;\mathrm{in\;}\mathbb{H}_{t},
\]
then\hfill{}(2.2.2)
\[
\pi_{t}\left(\left(a_{1},b_{1}\right)\right)=\left(\begin{array}{cc}
a_{1} & tb_{1}\\
\overline{b_{1}} & \overline{a_{1}}
\end{array}\right)\neq\left(\begin{array}{cc}
a_{2} & tb_{2}\\
\overline{b_{2}} & \overline{a_{2}}
\end{array}\right)=\pi_{t}\left(\left(a_{2},b_{2}\right)\right),
\]
in $M_{2}\left(\mathbb{C}\right)$. Furthermore, it satisfies that

\medskip{}

$\;\;\;\;$$\pi_{t}\left(\left(a_{1},b_{1}\right)+\left(a_{2},b_{2}\right)\right)=\left(\begin{array}{ccc}
a_{1}+a_{2} &  & t\left(b_{1}+b_{2}\right)\\
\\
\overline{b_{1}+b_{2}} &  & \overline{a_{1}+a_{2}}
\end{array}\right)$

\medskip{}

$\;\;\;\;\;\;\;\;$$=\left(\begin{array}{cc}
a_{1} & tb_{1}\\
\overline{b_{1}} & \overline{b_{2}}
\end{array}\right)+\left(\begin{array}{cc}
a_{2} & tb_{2}\\
\overline{b_{2}} & \overline{a_{2}}
\end{array}\right)=\pi_{t}\left(\left(a_{1},b_{1}\right)\right)+\pi_{t}\left(\left(a_{2},b_{2}\right)\right).$\hfill{}(2.2.3)

\medskip{}

\noindent Also, one has

\medskip{}

$\;\;\;\;$$\pi_{t}\left(\left(a_{1},b_{1}\right)\cdot_{t}\left(a_{2},b_{2}\right)\right)=\pi_{t}\left(\left(a_{1}a_{2}+tb_{1}\overline{b_{2}},\:a_{1}b_{2}+b_{1}\overline{a_{2}}\right)\right)$

\medskip{}

\noindent by (2.1.1)

\medskip{}

$\;\;\;\;\;\;\;\;\;\;$$=\left(\begin{array}{ccc}
a_{1}a_{2}+tb_{1}\overline{b_{2}} &  & t\left(a_{1}b_{2}+b_{1}\overline{a_{2}}\right)\\
\\
\overline{a_{1}b_{2}+b_{1}\overline{a_{2}}} &  & \overline{a_{1}a_{2}+tb_{1}\overline{b_{2}}}
\end{array}\right)$

\medskip{}

$\;\;\;\;\;\;\;\;\;\;$$=\left(\begin{array}{cc}
a_{1} & tb_{1}\\
\overline{b_{1}} & \overline{a_{1}}
\end{array}\right)\left(\begin{array}{cc}
a_{2} & tb_{2}\\
\overline{b_{2}} & \overline{a_{2}}
\end{array}\right)=\pi_{t}\left(\left(a_{1},b_{1}\right)\right)\pi_{t}\left(\left(a_{2},b_{2}\right)\right),$\hfill{}(2.2.4)

\medskip{}

\noindent where the multiplication ($\cdot$) in the far-right-hand
side of (2.2.4) is the usual matricial multiplication on $M_{2}\left(\mathbb{C}\right)$.

Since our $t$-scaled hypercomplex ring $\mathbb{H}_{t}=\left(\mathbb{C}^{2},+,\cdot_{t}\right)$
is identified with the 2-dimensional space $\mathbb{C}^{2}$ (set-theoretically),
one may / can understand this ring $\mathbb{H}_{t}$ as a topological
ring equipped with the usual topology for $\mathbb{C}^{2}$, for any
$t\in\mathbb{R}$. From below, we regard the ring $\mathbb{H}_{t}$
as a topological unital ring under the usual topology for $\mathbb{C}^{2}$.
\begin{lem}
The pair $\left(\mathbb{C}^{2},\:\pi_{t}\right)$ is an injective
Hilbert-space representation of the $t$-scaled hypercomplex ring
$\mathbb{H}_{t}$, where $\pi_{t}$ is an action (2.2.1).
\end{lem}

\begin{proof}
The morphism $\pi_{t}:\mathbb{H}_{t}\rightarrow M_{2}\left(\mathbb{C}\right)$
of (2.2.1) is a well-defined injective function by (2.2.2). Moreover,
this map $\pi_{t}$ satisfies the relations (2.2.3) and (2.2.4), and
hence, it is a(n algebraic) ring-action of $\mathbb{H}_{t}$, acting
on the 2-dimensional vector space $\mathbb{C}^{2}$. So, the pair
$\left(\mathbb{C}^{2},\pi_{t}\right)$ forms an algebraic representation
of $\mathbb{H}_{t}$. By regarding $\mathbb{H}_{t}$ and $M_{2}\left(\mathbb{C}\right)$
as topological spaces equipped with their usual topologies, then it
is not difficult to check that the ring-action $\pi_{t}$ is continuous
from $\mathbb{H}_{t}$ (which is homeomorphic to $\mathbb{C}^{2}$
as a topological space) into $M_{2}\left(\mathbb{C}\right)$ (which
is $*$-isomorphic to the $C^{*}$-algebra $B\left(\mathbb{C}^{2}\right)$).
Thus, the algebraic representation $\left(\mathbb{C}^{2},\pi_{t}\right)$
forms a Hilbert-space representation of $\mathbb{H}_{t}$ acting on
$\mathbb{C}^{2}$ via $\pi_{t}$.
\end{proof}
The above lemma shows that the $t$-scaled hypercomplex ring $\mathbb{H}_{t}$
is realized in the matricial algebra $M_{2}\left(\mathbb{C}\right)$
as 
\[
\pi_{t}\left(\mathbb{H}_{t}\right)=\left\{ \left(\begin{array}{cc}
a & tb\\
\overline{b} & \overline{a}
\end{array}\right)\in M_{2}\left(\mathbb{C}\right):\left(a,b\right)\in\mathbb{H}_{t}\right\} ,
\]
as an embedded topological ring in $M_{2}\left(\mathbb{C}\right)$.
\begin{defn}
The realization $\pi_{t}\left(\mathbb{H}_{t}\right)$ of the $t$-scaled
hypercomplex ring $\mathbb{H}_{t}$ is called the $t$-scaled (hypercomplex-)realization
of $\mathbb{H}_{t}$ (in $M_{2}\left(\mathbb{C}\right)$), for a scale
$t\in\mathbb{R}$. And we denote $\pi_{t}\left(\mathbb{H}_{t}\right)$
by $\mathcal{H}_{2}^{t}$. i.e.,
\[
\mathcal{H}_{2}^{t}\overset{\textrm{denote}}{=}\pi_{t}\left(\mathbb{H}_{t}\right)=\left\{ \left(\begin{array}{cc}
a & tb\\
\overline{b} & \overline{a}
\end{array}\right):\left(a,b\right)\in\mathbb{H}_{t}\right\} .
\]
Also, by $[\xi]_{t}$, we denote $\pi_{t}\left(\xi\right)\in\mathcal{H}_{2}^{t}$,
for all $\xi\in\mathbb{H}_{t}$.
\end{defn}

By the above lemma and definition, we obtain the following result.
\begin{thm}
For $t\in\mathbb{R}$, the corresponding $t$-scaled hypercomplex
ring $\mathbb{H}_{t}$ is topological-ring-isomorphic to the $t$-scaled
realization $\mathcal{H}_{2}^{t}$ in $M_{2}\left(\mathbb{C}\right)$.
i.e.,

\medskip{}

\hfill{}$\mathbb{H}_{t}\overset{\textrm{T.R}}{=}\mathit{\mathcal{H}_{2}^{t}\;\;\;\mathrm{in\;\;\;}M_{2}\left(\mathbb{C}\right)},$\hfill{}(2.2.5)

\medskip{}

\noindent where ``$\overset{\textrm{T.R}}{=}$'' means ``being
topological-ring-isomorphic to.''
\end{thm}

\begin{proof}
The relation (2.2.5) is proven by Lemma 4 and the injectivity (2.2.2)
of $\pi_{t}$.
\end{proof}
By the above theorem, one can realize that $\mathbb{H}_{t}$ and $\mathcal{H}_{2}^{t}$
as an identical topological ring, for a fixed $t\in\mathbb{R}$. Recall
that the relation (2.2.5) is independently shown in {[}2{]} and {[}3{]},
only for the quaternionic case where $t=-1$. 

\subsection{Scaled Hypercomplex Monoids}

Throughout this section, we fix a scale $t\in\mathbb{R}$, and the
corresponding $t$-scaled hypercomplex ring,
\[
\mathbb{H}_{t}=\left(\mathbb{C}^{2},\:+,\:\cdot_{t}\right),
\]
which is isomorphic to the $t$-scaled realization,
\[
\mathcal{H}_{2}^{t}=\left\{ \left(\begin{array}{cc}
a & tb\\
\overline{b} & \overline{a}
\end{array}\right)\in M_{2}\left(\mathbb{C}\right):\left(a,b\right)\in\mathbb{H}_{t}\right\} ,
\]
in $M_{2}\left(\mathbb{C}\right)$. Let
\[
\mathbb{H}_{t}^{\times}\overset{\textrm{denote}}{=}\mathbb{H}_{t}\setminus\left\{ \left(0,0\right)\right\} ,
\]
set-theoretically, where $\left(0,0\right)\in\mathbb{H}_{t}$ is the
($+$)-identity of the abelian group $\left(\mathbb{C}^{2},+\right)$.
Thus, by Proposition 2, this set forms a well-defined semigroup,
\[
\mathbb{H}_{t}^{\times}\overset{\textrm{denote}}{=}\left(\mathbb{H}_{t}^{\times},\;\cdot_{t}\right),
\]
equipped with its ($\cdot_{t}$)-identity $\left(1,0\right)$, and
hence, the pair $\mathbb{H}_{t}^{\times}$ is the maximal monoid embedded
in $\mathbb{H}_{2}^{t}$ up to the operation ($\cdot_{t}$).
\begin{defn}
The maximal monoid $\mathbb{H}_{t}^{\times}=\left(\mathbb{H}_{t}^{\times},\:\cdot_{t}\right)$,
embedded in the $t$-scaled hypercomplex ring $\mathbb{H}_{t}$, is
called the $t$-scaled hypercomplex monoid.
\end{defn}

By (2.2.5), it is trivial that:
\begin{cor}
The $t$-scaled hypercomplex monoid $\mathbb{H}_{t}^{\times}$ is
monoid-isomorphic to the monoid $\mathcal{H}_{2}^{t\times}\overset{\textrm{denote}}{=}\left(\mathcal{H}_{2}^{t\times},\:\cdot\right)$,
equipped with its identity,
\[
I_{2}=\left(\begin{array}{cc}
1 & 0\\
0 & 1
\end{array}\right)=\left(\begin{array}{cc}
1 & t\cdot0\\
0 & 1
\end{array}\right)=\left[\left(1,0\right)\right]_{t},
\]
the $\left(2\times2\right)$-identity matrix of $M_{2}\left(\mathbb{C}\right)$,
where ($\cdot$) is the usual matricial multiplication inherited from
that on $M_{2}\left(\mathbb{C}\right)$. i.e.,

\medskip{}

\hfill{}$\mathbb{H}_{t}^{\times}=\left(\mathbb{H}_{t}^{\times},\:\cdot_{t}\right)\overset{\textrm{Monoid}}{=}\left(\mathcal{H}_{2}^{t\times},\:\cdot\right)=\mathcal{H}_{2}^{t\times},$\hfill{}(2.3.1)

\medskip{}

\noindent where ``$\overset{\textrm{Monoid}}{=}$'' means ``being
monoid-isomorphic.''
\end{cor}

\begin{proof}
The isomorphic relation (2.3.1) is proven by the proof of Proposition
2, and that of Theorem 6. 
\end{proof}

\subsection{Invertibility on $\mathbb{H}_{t}$}

In this section, by identifying our $t$-scaled hypercomplex ring
$\mathbb{H}_{t}$ as its isomorphic realization $\mathcal{H}_{2}^{t}$,
we consider invertibility of elements of $\mathbb{H}_{t}$, for an
arbitrarily fixed $t\in\mathbb{R}$.

Observe first that, for any $\left(a,b\right)\in\mathbb{H}_{t}$ realized
to be $\left[\left(a,b\right)\right]_{t}\in\mathcal{H}_{2}^{t}$,
one can get that
\[
det\left(\left[\left(a,b\right)\right]_{t}\right)=det\left(\begin{array}{cc}
a & tb\\
\overline{b} & \overline{a}
\end{array}\right)=\left|a\right|^{2}-t\left|b\right|^{2},
\]
i.e.,\hfill{}(2.4.1)
\[
det\left(\left[\left(a,b\right)\right]_{t}\right)=\left|a\right|^{2}-t\left|b\right|^{2},
\]
where $det:M_{2}\left(\mathbb{C}\right)\rightarrow\mathbb{C}$ is
the determinant, and $\left|.\right|$ is the modulus on $\mathbb{C}$.
\begin{thm}
Let $\left(a,b\right)\in\mathbb{H}_{t}$, realized to be $\left[\left(a,b\right)\right]_{t}\in\mathcal{H}_{2}^{t}$.

\noindent (2.4.2) $det\left(\left[\left(a,b\right)\right]_{t}\right)=\left|a\right|^{2}-t\left|b\right|^{2}$.

\noindent (2.4.3) If either $\left|a\right|^{2}>t\left|b\right|^{2}$,
or $\left|a\right|^{2}<t\left|b\right|^{2}$, then $\left[\left(a,b\right)\right]_{t}$
is invertible ``in $M_{2}\left(\mathbb{C}\right)$,'' with its inverse
matrix,
\[
\left[\left(a,b\right)\right]_{t}^{-1}=\frac{1}{\left|a\right|^{2}-t\left|b\right|^{2}}\left(\begin{array}{cc}
\overline{a} & t\left(-b\right)\\
\overline{\left(-b\right)} & a
\end{array}\right).
\]

\noindent (2.4.4) If $\left|a\right|^{2}-t\left|b\right|^{2}\neq0$,
then $\left(a,b\right)\in\mathbb{H}_{t}$ is invertible in the sense
that there exists a unique $\left(c,d\right)\in\mathbb{H}_{t}$, such
that
\[
\left(a,b\right)\cdot_{t}\left(c,d\right)=\left(1,0\right)=\left(c,d\right)\cdot_{t}\left(a,b\right).
\]
In particular, one has that
\[
\left(c,d\right)=\left(\frac{\overline{a}}{\left|a\right|^{2}-t\left|b\right|^{2}},\:\frac{-b}{\left|a\right|^{2}-t\left|b\right|^{2}}\right)\in\mathbb{C}^{2}
\]

\noindent (2.4.5) Assume that $\left(a,b\right)$ is invertible in
$\mathbb{H}_{t}$ in the sense of (2.4.4). Then the inverse is also
contained ``in $\mathbb{H}_{t}$.''
\end{thm}

\begin{proof}
The statement (2.4.2) is shown by (2.4.1).

Note-and-recall that a matrix $A\in M_{n}\left(\mathbb{C}\right)$
is invertible in $M_{n}\left(\mathbb{C}\right)$, if and only if $det\left(A\right)\neq0$,
for all $n\in\mathbb{N}$. Therefore, 
\[
det\left(\left[\left(a,b\right)\right]_{t}\right)\neq0\Longleftrightarrow\left[\left(a,b\right)\right]_{t}\textrm{ is invertible in }M_{2}\left(\mathbb{C}\right).
\]
So, by (2.4.2),
\[
\left|a\right|^{2}-t\left|b\right|^{2}\neq0,\Longleftrightarrow\left[\left(a,b\right)\right]_{t}\textrm{ is invertible in }M_{2}\left(\mathbb{C}\right).
\]
Moreover, $\left|a\right|^{2}-t\left|b\right|^{2}\neq0$, if and only
if
\[
\left[\left(a,b\right)\right]_{t}^{-1}=\left(\begin{array}{cc}
a & tb\\
\overline{b} & \overline{a}
\end{array}\right)^{-1}=\frac{1}{\left|a\right|^{2}-t\left|b\right|^{2}}\left(\begin{array}{cc}
\overline{a} & -tb\\
-\overline{b} & a
\end{array}\right),
\]
in $M_{2}\left(\mathbb{C}\right)$. Therefore, the statement (2.4.3)
holds true in $M_{2}\left(\mathbb{C}\right)$.

By (2.4.3), one has $det\left(\left[\left(a,b\right)\right]_{t}\right)\neq0$,
if and only if
\[
\left[\left(a,b\right)\right]_{t}^{-1}=\left(\begin{array}{ccc}
\frac{\overline{a}}{\left|a\right|^{2}-t\left|b\right|^{2}} &  & t\left(\frac{-b}{\left|a\right|^{2}-t\left|b\right|^{2}}\right)\\
\\
\overline{\left(\frac{-b}{\left|a\right|^{2}-t\left|b\right|^{2}}\right)} &  & \frac{a}{\left|a\right|^{2}-t\left|b\right|^{2}}
\end{array}\right)\in M_{2}\left(\mathbb{C}\right),
\]
and it is actually contained ''in $\mathcal{H}_{2}^{t}$,'' satisfying
\[
\pi_{t}^{-1}\left(\begin{array}{ccc}
\frac{\overline{a}}{\left|a\right|^{2}-t\left|b\right|^{2}} &  & t\left(\frac{-b}{\left|a\right|^{2}-t\left|b\right|^{2}}\right)\\
\\
\overline{\left(\frac{-b}{\left|a\right|^{2}-t\left|b\right|^{2}}\right)} &  & \frac{a}{\left|a\right|^{2}-t\left|b\right|^{2}}
\end{array}\right)=\left(\frac{\overline{a}}{\left|a\right|^{2}-t\left|b\right|^{2}},\:\frac{-b}{\left|a\right|^{2}-t\left|b\right|^{2}}\right),
\]
in $\mathbb{H}_{t}$, by the injectivity of $\pi_{t}$. It shows that
$\left[\left(a,b\right)\right]_{t}^{-1}$ exists in $M_{2}\left(\mathbb{C}\right)$,
if and only if it is contained ``in $\mathcal{H}_{2}^{t}$.'' i.e.,
if $\left[\left(a,b\right)\right]_{t}$ is invertible, then its inverse
is also contained in $\mathcal{H}_{2}^{t}$, too, and vice versa.
So, the statements (2.2.4) and (2.2.5) hold.
\end{proof}
The above theorem not only characterizes the invertibility of the
monoidal elements of the $t$-scaled hypercomplex monoid $\mathbb{H}_{t}^{\times}$,
but also confirms that the inverses (if exist) are contained in the
monoid $\mathbb{H}_{t}^{\times}$. i.e.,
\[
\left(a,b\right)^{-1}\textrm{ exists,}\Longleftrightarrow\left(a,b\right)^{-1}=\left(\frac{\overline{a}}{\left|a\right|^{2}-t\left|b\right|^{2}},\:\frac{-b}{\left|a\right|^{2}-t\left|b\right|^{2}}\right),
\]
''in $\mathbb{H}_{t}^{\times}$,'' equivalently,
\[
\left[\left(a,b\right)^{-1}\right]_{t}=\left[\left(a,b\right)\right]_{t}^{-1}\;\mathrm{in\;}\mathcal{H}_{2}^{\times}.
\]

\begin{cor}
Let $\left(a,b\right)\in\mathbb{H}_{t}^{\times}$. Then it is invertible,
if and only if

\medskip{}

\hfill{}$\left[\left(a,b\right)^{-1}\right]_{t}=\left[\left(\frac{\overline{a}}{\left|a\right|^{2}-t\left|b\right|^{2}},\:\frac{-b}{\left|a\right|^{2}-t\left|b\right|^{2}}\right)\right]_{t}=\left[\left(a,b\right)\right]_{t}^{-1},$\hfill{}(2.4.6)

\medskip{}

\noindent in $\mathcal{H}_{2}^{\times}$, where $\left[\left(a,b\right)\right]_{t}^{-1}$
means the matricial inverse in $M_{2}\left(\mathbb{C}\right)$.
\end{cor}

\begin{proof}
The proof of (2.4.6) is immediately done by (2.4.3), (2.4.4) and (2.4.5).
\end{proof}
The above corollary can be re-stated by that: if $\xi\in\mathbb{H}_{t}^{\times}$
is invertible, then
\[
\pi_{t}\left(\xi^{-1}\right)=\left(\pi_{t}\left(\xi\right)\right)^{-1}\;\mathrm{in\;}\mathcal{H}_{2}^{t\times}.
\]

Now consider the cases where

\medskip{}

\hfill{}$\left|a\right|^{2}-t\left|b\right|^{2}=0\Longleftrightarrow\left|a\right|^{2}=t\left|b\right|^{2},$\hfill{}(2.4.7)

\medskip{}

\noindent in $\mathbb{R}$. As we have seen above, the condition (2.4.7)
holds for $\left(a,b\right)\in\mathbb{H}_{t}$, if and only if $\left(a,b\right)$
is not invertible in $\mathbb{H}_{t}$ (and hence, its realization
$\left[\left(a,b\right)\right]_{t}$ is not invertible in $M_{2}\left(\mathbb{C}\right)$,
and hence, in $\mathcal{H}_{2}^{t}$). Clearly, we are not interested
in the ($+$)-identity $\left(0,0\right)$ of $\mathbb{H}_{t}$ automatically
satisfying the condition (2.4.7). So, without loss of generality,
we focus on elements $\left(a,b\right)$ of the $t$-scaled hypercomplex
monoid $\mathbb{H}_{t}^{\times}$ (or, its realizations $\left[\left(a,b\right)\right]_{t}$
of $\mathcal{H}_{2}^{t\times}$), satisfying the condition (2.4.7).

Recall that an algebraic triple, $\left(X,+,\cdot\right)$, is a noncommutative
field, if (i) $\left(X,+\right)$ is an abelian group, (ii) $\left(X^{\times},\cdot\right)$
forms a non-abelian group, and (iii) the operations ($+$) and ($\cdot$)
are left-and-right distributive. For instance, the quaternions $\mathbb{H}=\mathbb{H}_{-1}$
is a noncommutative field (e.g., {[}2{]} and {[}3{]}).
\begin{thm}
Suppose the fixed scale $t\in\mathbb{R}$ is negative, i.e., $t<0$
in $\mathbb{R}$. Then ``all'' elements $\left(a,b\right)$ of the
$t$-scaled hypercomplex monoid $\mathbb{H}_{t}^{\times}$ are invertible
in $\mathbb{H}_{t}$, with their inverses,
\[
\left(\frac{\overline{a}}{\left|a\right|^{2}-t\left|b\right|^{2}},\;\frac{-b}{\left|a\right|^{2}-t\left|b\right|^{2}}\right)\in\mathbb{H}_{t}^{\times}.
\]
i.e.,\hfill{}(2.4.8)
\[
t<0\;\mathrm{in\;}\mathbb{R}\Longrightarrow\mathbb{H}_{t}\textrm{ is a noncommutative field.}
\]
\end{thm}

\begin{proof}
Suppose the scale $t\in\mathbb{R}$ is negative. Then, for any $\left(a,b\right)\in\mathbb{H}_{t}^{\times}$,
\[
\left|a\right|^{2}\neq t\left|b\right|^{2}\Longleftrightarrow\left|a\right|^{2}-t\left|b\right|^{2}>0,
\]
since $\left(a,b\right)\neq\left(0,0\right)$. i.e., if $t<0$, then
every element $\left(a,b\right)\in\mathbb{H}_{t}^{\times}$ does ``not''
satisfy the condition (2.4.7). It implies that if $t<0$, then every
element $\left(a,b\right)\in\mathbb{H}_{t}^{\times}$ is invertible
in $\mathbb{H}_{t}^{\times}$, by (2.4.4) and (2.4.5); and the inverse
is determined to be (2.4.6) in $\mathbb{H}_{t}^{\times}$. Thus, the
pair $\mathbb{H}_{t}^{\times}=\left(\mathbb{H}_{t}^{\times},\:\cdot_{t}\right)$
forms a group which is not abelian by (2.1.1) and (2.2.4).

Therefore, if $t<0$ in $\mathbb{R}$, then the $t$-scaled hypercomplex
ring $\mathbb{H}_{t}$ becomes a noncommutative field, proving the
statement (2.4.8).
\end{proof}
The above theorem characterizes that the algebraic structure of scaled
hypercomplex rings $\left\{ \mathbb{H}_{t}\right\} _{t<0}$ as noncommutative
fields.
\begin{thm}
Suppose $t=0$ in $\mathbb{R}$. Then an element $\left(a,b\right)$
of the $0$-scaled hypercomplex monoid $\mathbb{H}_{0}^{\times}$
is invertible in $\mathbb{H}_{0}$, with their inverses,
\[
\left(\frac{\overline{a}}{\left|a\right|^{2}},\;\frac{-b}{\left|a\right|^{2}}\right)\in\mathbb{H}_{0}^{\times},
\]
if and only if $a\neq0$ in $\mathbb{C}$, if and only if only the
elements of the subset,

\medskip{}

\hfill{}$\left\{ \left(a,b\right)\in\mathbb{H}_{0}^{\times}:a\neq0\right\} \textrm{ of }\mathbb{H}_{0}^{\times}$\hfill{}(2.4.9)

\medskip{}

\noindent are invertible in $\mathbb{H}_{0}^{\times}$, if and only
if $\left(0,b\right)\in\mathbb{H}_{0}^{\times}$ are not invertible
in $\mathbb{H}_{0}^{\times}$, for all $b\in\mathbb{C}$.
\end{thm}

\begin{proof}
Assume that we have the zero scale, i.e., $t=0$ in $\mathbb{R}$.
Then, by (2.4.7),
\[
\left|a\right|^{2}=0\cdot\left|b\right|^{2}\Longleftrightarrow\left|a\right|^{2}=0\Longleftrightarrow a=0\textrm{ in }\mathbb{C},
\]
if and only if $\left(0,b\right)\in\mathbb{H}_{0}^{\times}$ are not
invertible in $\mathbb{H}_{0}^{\times}$, for all $b\in\mathbb{C}$,
if and only if all elements $\left(a,b\right)$, contained in the
subset (2.4.9), are invertible in $\mathbb{H}_{0}^{\times}$.

Observe that $\left(a,b\right)$ is contained in the subset (2.4.9)
of $\mathbb{H}_{0}^{\times}$, if and only if

\medskip{}

$\;\;\;\;$$\left[\left(a,b\right)\right]_{0}\left[\left(\frac{\overline{a}}{\left|a\right|^{2}},\:\frac{-b}{\left|a\right|^{2}}\right)\right]_{0}=\left(\begin{array}{cc}
a & 0\\
\overline{b} & \overline{a}
\end{array}\right)\left(\begin{array}{ccc}
\frac{\overline{a}}{\left|a\right|^{2}} &  & 0\\
\\
\frac{\overline{-b}}{\left|a\right|^{2}} &  & \frac{a}{\left|a\right|^{2}}
\end{array}\right)$

\medskip{}

$\;\;\;\;\;\;\;\;\;\;\;\;\;\;\;\;$$=\left(\begin{array}{cc}
1 & 0\\
0 & 1
\end{array}\right)=\left(\begin{array}{ccc}
\frac{\overline{a}}{\left|a\right|^{2}} &  & 0\\
\\
\frac{\overline{-b}}{\left|a\right|^{2}} &  & \frac{a}{\left|a\right|^{2}}
\end{array}\right)\left(\begin{array}{cc}
a & 0\\
\overline{b} & \overline{a}
\end{array}\right)$

\medskip{}

$\;\;\;\;\;\;\;\;\;\;\;\;\;\;\;\;$$=\left[\left(\frac{\overline{a}}{\left|a\right|^{2}},\:\frac{-b}{\left|a\right|^{2}}\right)\right]_{0}\left[\left(a,b\right)\right]_{0},$

\medskip{}

\noindent in $\mathbb{H}_{0}^{\times}$. Therefore, if exists, $\left(a,b\right)^{-1}=\left(\frac{\overline{a}}{\left|a\right|^{2}},\:\frac{-b}{\left|a\right|^{2}}\right)$
in $\mathbb{H}_{0}^{\times}$.
\end{proof}
The above theorem shows that if we have the zero-scale in $\mathbb{R}$,
then our $0$-scaled hypercomplex ring $\mathbb{H}_{0}$ cannot be
a noncommutative field. It directly illustrates that the algebra on
the quaternions $\mathbb{H}=\mathbb{H}_{-1}$, and the algebra on
the scaled-hypercomplex rings $\left\{ \mathbb{H}_{t}\right\} _{t\in\mathbb{R}\setminus\left\{ -1\right\} }$
can be different in general, especially, when $t\geq0$.
\begin{thm}
Suppose the scale $t\in\mathbb{R}$ is positive, i.e., $t>0$ in $\mathbb{R}$.
Then an element $\left(a,b\right)\in\mathbb{H}_{t}^{\times}$ is invertible
in $\mathbb{H}_{t}^{\times}$ with its inverse,
\[
\left(\frac{\overline{a}}{\left|a\right|^{2}-t\left|b\right|^{2}},\;\frac{-b}{\left|a\right|^{2}-t\left|b\right|^{2}}\right)\in\mathbb{H}_{t}^{\times},
\]
if and only if $\left|a\right|^{2}\neq t\left|b\right|^{2}$ in $\mathbb{R}_{0}^{+}=\left\{ r\in\mathbb{R}:r\geq0\right\} $,
if and only if $\left(a,b\right)$ is contained in the subset,

\medskip{}

\hfill{}$\left\{ \left(a,b\right):\left|a\right|^{2}\neq t\left|b\right|^{2}\;\mathrm{in\;}\mathbb{R}_{0}^{+}\right\} ,$\hfill{}(2.4.10)

\medskip{}

\noindent of $\mathbb{H}_{t}^{\times}$. As application, if $t>0$
in $\mathbb{R}$, then the all elements of

\medskip{}

\noindent \hfill{}$\left\{ \left(a,0\right)\in\mathbb{H}_{t}:a\in\mathbb{C}^{\times}\right\} \cup\left\{ \left(0,b\right)\in\mathbb{H}_{t}:b\in\mathbb{C}^{\times}\right\} ,$\hfill{}(2.4.11)

\medskip{}

\noindent are invertible in $\mathbb{H}_{t}$, where $\mathbb{C}^{\times}=\mathbb{C}\setminus\left\{ 0\right\} $.
\end{thm}

\begin{proof}
Assume that $t>0$ in $\mathbb{R}$, and $\mathbb{H}_{t}^{\times}$,
the corresponding $t$-scaled hypercomplex monoid. Then $\left(a,b\right)\in\mathbb{H}_{t}^{\times}$
is invertible in $\mathbb{H}_{t}^{\times}$, if and only if the condition
(2.4.7) does not hold, if and only if
\[
\left|a\right|^{2}\neq t\left|b\right|^{2}\Longleftrightarrow\mathrm{either\;}\left|a\right|^{2}>t\left|b\right|^{2},\;\mathrm{or\;}\left|a\right|^{2}<t\left|b\right|^{2},
\]
in $\mathbb{R}_{0}^{+}$, since $t>0$. Therefore, if $t>0$ in $\mathbb{R}$,
then an element $\left(a,b\right)$ is invertible in $\mathbb{H}_{t}^{\times}$,
if and only if
\[
\mathrm{either\;}\left|a\right|^{2}>t\left|b\right|^{2},\;\mathrm{or\;}\left|a\right|^{2}<t\left|b\right|^{2}\;\mathrm{in\;}\mathbb{R}_{0}^{+},
\]
if and only if $\left(a,b\right)$ is contained in the subset (2.4.10)
in $\mathbb{H}_{t}^{\times}$.

In particular, for $t>0$ in $\mathbb{R}$, (i) if $\left(a,0\right)\in\mathbb{H}_{t}^{\times}$
with $a\in\mathbb{C}^{\times}$, then $\left|a\right|^{2}>0$; and
(ii) if $\left(0,b\right)\in\mathbb{H}_{t}^{\times}$ with $b\in\mathbb{C}^{\times}$,
then $0<t\left|b\right|^{2}$. Therefore, the subset (2.4.11) is properly
contained in the subset (2.4.10) in $\mathbb{H}_{t}^{\times}$, whenever
$t>0$. So, all elements, formed by $\left(a,0\right),$or by $\left(0,b\right)$
with $a,b\in\mathbb{C}^{\times}$, are invertible in $\mathbb{H}_{t}^{\times}$.
\end{proof}
The above theorem characterizes the invertibility on the $t$-scaled
hypercomplex monoid $\mathbb{H}_{t}^{\times}$, where the scale $t$
is positive in $\mathbb{R}$. Theorems 11, 12 and 13 refine Theorem
8, case-by-case. We again summarize the main results.
\begin{cor}
Let $\mathbb{H}_{t}^{\times}$ be the $t$-scaled hypercomplex monoid.
If $t<0$, then all nonzero elements of $\mathbb{H}_{t}^{\times}$
are invertible; and if $t=0$, then
\[
\left\{ \left(a,b\right)\in\mathbb{H}_{0}^{\times}:a\neq0\right\} 
\]
is the invertible proper subset of $\mathbb{H}_{0}^{\times}$; and
if $t>0$, then
\[
\left\{ \left(a,b\right):\left|a\right|^{2}\neq t\left|b\right|^{2}\;\mathrm{in\;}\mathbb{R}_{0}^{+}\right\} 
\]
is the invertible proper subset of $\mathbb{H}_{t}^{\times}$, where
``invertible subset of $\mathbb{H}_{t}^{\times}$'' means ``a subset
of $\mathbb{H}_{t}^{\times}$ containing of all invertible elements.'' 
\end{cor}

\begin{proof}
This corollary is nothing but a summary of Theorems 11, 12 and 13.
\end{proof}

\subsection{Decompositions of the Nonnegatively-Scaled Hypercomplex Rings}

In this section, we consider a certain decomposition of the $t$-scaled
hypercomplex ring $\mathbb{H}_{t}$, for an arbitrary fixed ``positive''
scale $t>0$ in $\mathbb{R}$. Recall that, as we have seen in Section
2.4, the negatively-scaled hypercomplex rings $\left\{ \mathbb{H}_{s}\right\} _{s<0}$
are noncommutative fields by (2.4.8), equivalently, the negatively-scaled
hypercomplex monoids $\left\{ \mathbb{H}_{s}^{\times}\right\} _{s<0}$
are non-abelian groups. However, if $t\geq0$, then $\mathbb{H}_{t}$
cannot be a noncommutative field in general, by (2.4.9) and (2.4.10).
We here concentrate on such cases.

Let $t\geq0$ and $\mathbb{H}_{t}$, the corresponding $t$-scaled
hypercomplex ring. Partition $\mathbb{H}_{t}$ by
\[
\mathbb{H}_{t}=\mathbb{H}_{t}^{inv}\sqcup\mathbb{H}_{t}^{sing}
\]
with\hfill{}(2.5.1)
\[
\mathbb{H}_{t}^{inv}=\left\{ \left(a,b\right):\left|a\right|^{2}\neq t\left|b\right|^{2}\right\} ,
\]
and
\[
\mathbb{H}_{t}^{sing}=\left\{ \left(a,b\right):\left|a\right|^{2}=t\left|b\right|^{2}\right\} ,
\]
where $\sqcup$ is the disjoint union. By (2.4.9) and (2.4.10), $\left(a,b\right)\in\mathbb{H}_{t}^{inv}$,
if and only if it is invertible, equivalently, $\left(a,b\right)\in\mathbb{H}_{t}^{sing}$,
if and only if it is not invertible, in $\mathbb{H}_{t}$.

Recall-and-note that the determinant is a multiplicative map on $M_{n}\left(\mathbb{C}\right)$,
for all $n\in\mathbb{N}$, in the sense that:

\medskip{}

\hfill{}$det\left(AB\right)=det\left(A\right)det\left(B\right),\;\forall A,B\in M_{n}\left(\mathbb{C}\right).$\hfill{}(2.5.2)

\medskip{}

\noindent Thus, by (2.5.2), one has

\medskip{}

\hfill{}$\xi,\eta\in\mathbb{H}_{t}^{inv}\Rightarrow det\left(\left[\xi\cdot_{t}\eta\right]_{t}\right)=det\left(\left[\xi\right]_{t}\left[\eta\right]_{t}\right)\neq0$.\hfill{}(2.5.3)
\begin{lem}
Let $t\geq0$ in $\mathbb{R}$. Then the subset $\mathbb{H}_{t}^{inv}\overset{\textrm{denote}}{=}\left(\mathbb{H}_{t}^{inv},\cdot_{t}\right)$
of the $t$-scaled hypercomplex monoid $\mathbb{H}_{t}^{\times}$
forms a non-abelian group. i.e., $\mathbb{H}_{t}^{inv}$ is not only
a sub-monoid, but also an embedded group in $\mathbb{H}_{t}^{\times}$.
\end{lem}

\begin{proof}
By (2.5.2), if $\xi,\eta\in\mathbb{H}_{t}^{inv}$, then $\xi\cdot_{t}\eta\in\mathbb{H}_{t}^{inv}$,
too. i.e., the operation ($\cdot_{t}$) is closed, and associative
on $\mathbb{H}_{t}^{inv}$. Also, the ($\cdot_{t}$)-identity $\left(1,0\right)$
is contained in $\mathbb{H}_{t}^{inv}$ by (2.5.1). Therefore, the
sub-structure $\left(\mathbb{H}_{t}^{inv},\cdot_{t}\right)$ forms
a sub-monoid of $\mathbb{H}_{t}^{\times}$. But, by (2.4.8) and (2.5.3),
each element $\xi\in\mathbb{H}_{t}^{inv}$ has its ($\cdot_{t}$)-inverse
$\xi^{-1}$ contained in $\mathbb{H}_{t}^{inv}$. It shows that $\mathbb{H}_{t}^{inv}$
forms a non-abelian group in the monoid $\mathbb{H}_{t}^{\times}$.
\end{proof}
By the partition (2.5.1) and the multiplicativity (2.5.3), one can
obtain the following equivalent result of the above theorem.
\begin{lem}
Let $t\geq0$ in $\mathbb{R}$. Then the pair 
\[
\mathbb{H}_{t}^{\times sing}\overset{\textrm{denote}}{=}\left(\mathbb{H}_{t}^{sing}\cap\mathbb{H}_{t}^{\times},\:\cdot_{t}\right)=\left(\mathbb{H}_{t}^{sing}\setminus\left\{ \left(0,0\right)\right\} ,\:\cdot_{t}\right)
\]
 forms a semigroup without identity in the $t$-scaled hypercomplex
monoid $\mathbb{H}_{t}^{\times}$.
\end{lem}

\begin{proof}
By (2.5.2) and (2.5.3), the operation ($\cdot_{t}$) is closed and
associative on the set,
\[
\mathbb{H}_{t}^{\times sing}\overset{\textrm{def}}{=}\mathbb{H}_{t}^{\times}\cap\mathbb{H}_{t}^{sing}=\mathbb{H}_{t}^{sing}\setminus\left\{ \left(0,0\right)\right\} .
\]
However, the ($\cdot_{t}$)-identity $\left(1,0\right)$ is not contained
in $\mathbb{H}_{t}^{\times sing}$, since $I_{2}=\left[\left(1,0\right)\right]_{t}$
is in $\mathbb{H}_{t}^{inv}$. So, in the monoid $\mathbb{H}_{t}^{\times}$,
the sub-structure $\left(\mathbb{H}_{t}^{\times sing},\:\cdot_{t}\right)$
forms a semigroup (without identity).
\end{proof}
The above lemma definitely includes the fact that: $\left(\mathbb{H}_{t}^{sing},\cdot_{t}\right)$
is just a semigroup (without identity), which is not a sub-monoid
of $\mathbb{H}_{t}^{\times}$ (and hence, not a group).

The above two algebraic characterizations show that the set-theoretical
decomposition (2.5.1) induces an algebraic decomposition of the $t$-scaled
hypercomplex monoid $\mathbb{H}_{t}^{\times}$,
\[
\mathbb{H}_{t}^{\times}=\left(\mathbb{H}_{t}^{inv},\:\cdot_{t}\right)\sqcup\left(\mathbb{H}_{t}^{\times sing},\:\cdot_{t}\right),
\]
where\hfill{}(2.5.4)
\[
\mathbb{H}_{t}^{inv}=\left\{ \left(a,b\right)\in\mathbb{H}_{t}^{\times}:\left|a\right|^{2}\neq t\left|b\right|^{2}\right\} ,
\]
and
\[
\mathbb{H}_{t}^{\times sing}=\left\{ \left(a,b\right)\in\mathbb{H}_{t}^{\times}:\left|a\right|^{2}=t\left|b\right|^{2}\right\} ,
\]
whenever $t\geq0$ in $\mathbb{R}$.
\begin{thm}
For $t\geq0$ in $\mathbb{R}$, the $t$-scaled hypercomplex monoid
$\mathbb{H}_{t}^{\times}$ is algebraically decomposed to be
\[
\mathbb{H}_{t}^{\times}=\mathbb{H}_{t}^{inv}\sqcup\mathbb{H}_{t}^{\times sing},
\]
where $\mathbb{H}_{t}^{inv}$ is the group, and $\mathbb{H}_{t}^{\times sing}$
is the semigroup without identity in (2.5.4).
\end{thm}

\begin{proof}
The algebraic decomposition,
\[
\mathbb{H}_{t}^{\times}=\mathbb{H}_{t}^{inv}\sqcup\mathbb{H}_{t}^{\times sing},
\]
of the $t$-scaled hypercomplex monoid $\mathbb{H}_{t}^{\times}$
is obtained by the set-theoretic decomposition (2.5.1) of $\mathbb{H}_{t}^{\times}$,
the above two lemmas, and (2.5.4).
\end{proof}
By the above theorem, one can have the following concepts whenever
a given scale $t$ is nonnegative in $\mathbb{R}$.
\begin{defn}
Let $t\geq0$ in $\mathbb{R}$, and $\mathbb{H}_{t}^{\times}$, the
$t$-scaled hypercomplex monoid. The algebraic block,
\[
\mathbb{H}_{t}^{inv}=\left(\left\{ \left(a,b\right)\in\mathbb{H}_{t}^{\times}:\left|a\right|^{2}\neq t\left|b\right|^{2}\right\} ,\:\cdot_{t}\right),
\]
is called the group-part of $\mathbb{H}_{t}^{\times}$ (or, of $\mathbb{H}_{t}$),
and the other algebraic block,
\[
\mathbb{H}_{t}^{\times sing}=\left(\left\{ \left(a,b\right)\in\mathbb{H}_{t}^{\times}:\left|a\right|^{2}=t\left|b\right|^{2}\right\} ,\:\cdot_{t}\right),
\]
is called the semigroup-part of $\mathbb{H}_{t}^{\times}$ (or, of
$\mathbb{H}_{t}$).
\end{defn}

By the above definition, Theorem 17 can be re-stated that: if a scale
$t$ is nonnegative in $\mathbb{R}$, then the $t$-scaled hypercomplex
monoid $\mathbb{H}_{t}^{\times}$ is decomposed to be the group-part
$\mathbb{H}_{t}^{inv}$ and the semigroup-part $\mathbb{H}_{t}^{\times sing}$. 

One may / can say that if $t<0$ in $\mathbb{R}$, then the semigroup-part
$\mathbb{H}_{t}^{\times sing}$ is empty in $\mathbb{H}_{t}^{\times}$.
Indeed, for any scale $t\in\mathbb{R}$, the $t$-scaled hypercomplex
monoid $\mathbb{H}_{t}$ is decomposed to be (2.5.4). As we have seen
in this section, if $t\geq0$, then the semigroup-part $\mathbb{H}_{t}^{\times sing}$
is nonempty, meanwhile, as we considered in Section 2.4, if $t<0$,
then the semigroup-part $\mathbb{H}_{t}^{\times sing}$ is empty,
equivalently, the $t$-scaled hypercomplex monoid $\mathbb{H}_{t}^{\times}$
is identified with its group-part $\mathbb{H}_{t}^{inv}$, i.e., $\mathbb{H}_{t}^{\times}=\mathbb{H}_{t}^{inv}$
in $\mathbb{H}_{t}$, whenever $t<0$.
\begin{cor}
For every $t\in\mathbb{R}$, the $t$-scaled hypercomplex monoid $\mathbb{H}_{t}^{\times}$
is partitioned by
\[
\mathbb{H}_{t}^{\times}=\mathbb{H}_{t}^{inv}\sqcup\mathbb{H}_{t}^{\times sing},
\]
where the group-part $\mathbb{H}_{t}^{inv}$ and the semigroup-part
$\mathbb{H}_{t}^{\times sing}$ are in the sense of (2.5.4). In particular,
if $t<0$, then
\[
\mathbb{H}_{t}^{\times sing}=\textrm{Ø}\Longleftrightarrow\mathbb{H}_{t}^{\times}=\mathbb{H}_{t}^{inv};
\]
meanwhile, if $t\geq0$, then $\mathbb{H}_{t}^{\times sing}$ is a
non-empty proper subset of $\mathbb{H}_{t}^{\times}$.
\end{cor}

\begin{proof}
It is shown conceptually by the discussion of the very above paragraph.
Also, see Theorems 11 and 17.
\end{proof}

\section{Spectral Analysis on $\left\{ \mathbb{H}_{t}\right\} _{t\in\mathbb{R}}$
Under $\left\{ \left(\mathbb{C}^{2},\pi_{t}\right)\right\} _{t\in\mathbb{R}}$}

Throughout this section, we fix an arbitrary scale $t\in\mathbb{R}$,
and the corresponding $t$-scaled hypercomplex ring,
\[
\mathbb{H}_{t}=\left(\mathbb{C}^{2},\:+,\:\cdot_{t}\right),
\]
containing its hypercomplex monoid $\mathbb{H}_{t}^{\times}=\left(\mathbb{H}_{t}^{\times},\:\cdot_{t}\right)$.
In Section 2, we showed that for a scale $t\in\mathbb{R}$, the monoid
$\mathbb{H}_{t}^{\times}$ is partitioned by
\[
\mathbb{H}_{t}^{\times}=\mathbb{H}_{t}^{inv}\sqcup\mathbb{H}_{t}^{\times sing},
\]
where $\mathbb{H}_{t}^{inv}$ is the group-part, and $\mathbb{H}_{t}^{\times sing}$
is the semigroup-part of $\mathbb{H}_{t}$. In particular, if $t<0$,
then the semigroup-part $\mathbb{H}_{t}^{\times sing}$ is empty in
$\mathbb{H}_{t}^{\times}$, equivalently, $\mathbb{H}_{t}^{\times}=\mathbb{H}_{t}^{inv}$
in $\mathbb{H}_{t}$, meanwhile, if $t\geq0$, then $\mathbb{H}_{t}^{\times sing}$
is a non-empty proper subset of $\mathbb{H}_{t}^{\times}$.

Motivated by such an analysis of invertibility on $\mathbb{H}_{t}$,
we here consider spectral analysis on $\mathbb{H}_{t}$.

\subsection{Hypercomplex-Spectral Forms on $\mathbb{H}_{t}$}

For $t\in\mathbb{R}$, let $\mathbb{H}_{t}$ be the $t$-scaled hypercomplex
ring realized to be
\[
\mathcal{H}_{2}^{t}=\pi_{t}\left(\mathbb{H}_{t}\right)=\left\{ \left(\begin{array}{cc}
a & tb\\
\overline{b} & \overline{a}
\end{array}\right)\in M_{2}\left(\mathbb{C}\right):\left(a,b\right)\in\mathbb{H}_{t}\right\} ,
\]
in $M_{2}\left(\mathbb{C}\right)$ under the Hilbert-space representation
$\Pi_{t}=\left(\mathbb{C}^{2},\pi_{t}\right)$ of $\mathbb{H}_{t}$.

Let $\left(a,b\right)\in\mathbb{H}_{t}$ be an arbitrary element with
\[
\pi_{t}\left(a,b\right)=\left[\left(a,b\right)\right]_{t}=\left(\begin{array}{cc}
a & tb\\
\overline{b} & \overline{a}
\end{array}\right)\in\mathcal{H}_{2}^{t}.
\]
Then, in a variable $z$ on $\mathbb{C}$,

\medskip{}

$\;\;\;\;$$det\left(\left[\left(a,b\right)\right]_{t}-z\left[\left(1,0\right)\right]_{t}\right)=det\left(\begin{array}{ccc}
a-z &  & tb\\
\\
\overline{b} &  & \overline{a}-z
\end{array}\right)$

\medskip{}

$\;\;\;\;\;\;\;\;\;\;\;\;\;\;\;\;$$=\left(a-z\right)\left(\overline{a}-z\right)-t\left|b\right|^{2}$

\medskip{}

$\;\;\;\;\;\;\;\;\;\;\;\;\;\;\;\;$$=\left|a\right|^{2}-az-\overline{a}z+z^{2}-t\left|b\right|^{2}$

\medskip{}

$\;\;\;\;\;\;\;\;\;\;\;\;\;\;\;\;$$=z^{2}-\left(a+\overline{a}\right)z+\left(\left|a\right|^{2}-t\left|b\right|^{2}\right)$

\medskip{}

$\;\;\;\;\;\;\;\;\;\;\;\;\;\;\;\;$$=z^{2}-2Re\left(a\right)z+det\left(\left[\left(a,b\right)\right]_{t}\right)$,\hfill{}(3.1.1)

\medskip{}

\noindent where $Re\left(a\right)$ is the real part of $\mathit{a}$
in $\mathbb{C}$, and
\[
det\left(\left[\left(a,b\right)\right]_{t}\right)=\left|a\right|^{2}-t\left|b\right|^{2},
\]
by (2.4.2). Thus, the equation,
\[
det\left(\left[\left(a,b\right)\right]_{t}-z\left[\left(1,0\right)\right]_{t}\right)=0,
\]
in a variable $z$ on $\mathbb{C}$, has its solutions,
\[
z=\frac{2Re\left(a\right)\pm\sqrt{4Re\left(a\right)^{2}-4det\left(\left[\left(a,b\right)\right]_{t}\right)}}{2},
\]
$\Longleftrightarrow$\hfill{}(3.1.2)
\[
z=Re\left(a\right)\pm\sqrt{Re\left(a\right)^{2}-det\left(\left[\left(a,b\right)\right]_{t}\right)}.
\]

Recall that a matrix $A\in M_{n}\left(\mathbb{C}\right)$, for any
$n\in\mathbb{N}$, has its spectrum,
\[
spec\left(A\right)=\left\{ \lambda\in\mathbb{C}:det\left(A-\lambda I_{n}\right)=0\right\} ,
\]
equivalently,\hfill{}(3.1.3)
\[
spec\left(A\right)=\left\{ \lambda\in\mathbb{C}:\exists\eta\in\mathbb{C}^{n},\;s.t.,\;A\eta=\lambda\eta\right\} ,
\]
if and only if
\[
spec\left(A\right)=\left\{ \lambda\in\mathbb{C}:A-\lambda I_{n}\textrm{ is not invertible in }M_{n}\left(\mathbb{C}\right)\right\} ,
\]
as a nonempty discrete (compact) subset of $\mathbb{C}$, where $I_{n}$
is the identity matrix of $M_{n}\left(\mathbb{C}\right)$ (e.g., {[}8{]}).
More generally, if $T\in B\left(H\right)$ is an operator on a Hilbert
space $H$, then the spectrum $\sigma\left(T\right)$ of $T$ is defined
to be a nonempty compact subset,
\[
\sigma\left(T\right)=\left\{ z\in\mathbb{C}:T-zI_{H}\textrm{ is not invertible on }H\right\} ,
\]
where $I_{H}$ is the identity operator of $B\left(H\right)$. Remark
that if $H$ is infinite-dimensional, then $\sigma\left(T\right)$
is not a discrete subset of $\mathbb{C}$ as in (3.1.3), in general
(e.g., {[}9{]}).
\begin{thm}
Let $\left(a,b\right)\in\mathbb{H}_{t}$ realized to be $\left[\left(a,b\right)\right]_{t}\in\mathcal{H}_{2}^{t}$.
Then
\[
spec\left(\left[\left(a,b\right)\right]_{t}\right)=\left\{ Re\left(a\right)\pm\sqrt{Re\left(a\right)^{2}-det\left(\left[\left(a,b\right)\right]_{t}\right)}\right\} ,
\]
in $\mathbb{C}$. More precisely, if
\[
a=x+yi,\;b=u+vi\in\mathbb{C},
\]
with $x,y,u,v\in\mathbb{R}$ and $i=\sqrt{-1}$ in $\mathbb{C}$,
then

\medskip{}

\hfill{}$spec\left(\left[\left(a,b\right)\right]_{t}\right)=\left\{ x\pm i\sqrt{y^{2}-tu^{2}-tv^{2}}\right\} \;\mathrm{in\;}\mathbb{C}.$\hfill{}(3.1.4)
\end{thm}

\begin{proof}
The realization $\left[\left(a,b\right)\right]_{t}=\left(\begin{array}{cc}
a & tb\\
\overline{b} & \overline{a}
\end{array}\right)\in\mathcal{H}_{2}^{t}$ of a hypercomplex number $\left(a,b\right)\in\mathbb{H}_{t}$ has
its spectrum,
\[
spec\left(\left[\left(a,b\right)\right]_{t}\right)=\left\{ Re\left(a\right)\pm\sqrt{Re\left(a\right)^{2}-\left(\left|a\right|^{2}-t\left|b\right|^{2}\right)}\right\} ,
\]
in $\mathbb{C}$, by (3.1.2) and (3.1.3). If
\[
a=x+yi,\;\mathrm{and\;}b=u+vi\;\mathrm{in\;}\mathbb{C},
\]
with $x,y,u,v\in\mathbb{R}$ and $i=\sqrt{-1}$ in $\mathbb{C}$,
then
\[
Re\left(a\right)=x,
\]
and
\[
\left|a\right|^{2}-t\left|b\right|^{2}=\left(x^{2}+y^{2}\right)-t\left(u^{2}+v^{2}\right),
\]
in $\mathbb{R}$, and hence,
\[
spec\left(\left[\left(a,b\right)\right]_{t}\right)=\left\{ x\pm\sqrt{-y^{2}+tu^{2}+tv^{2}}\right\} ,
\]
if and only if
\[
spec\left(\left[\left(a,b\right)\right]_{t}\right)=\left\{ x\pm i\sqrt{y^{2}-tu^{2}-tv^{2}}\right\} ,
\]
in $\mathbb{C}$. Therefore, the set-equality (3.1.4) holds.
\end{proof}
From below, for our purposes, we let
\[
a=x+yi\;\mathrm{and\;}b=u+vi\;\mathrm{in\;}\mathbb{C},
\]
with\hfill{}(3.1.5)
\[
x,y,u,v\in\mathbb{R},\;\mathrm{and\;}i=\sqrt{-1}.
\]
The above theorem can be refined by the following result.
\begin{cor}
Let $\left(a,b\right)\in\mathbb{H}_{t}$, realized to be $\left[\left(a,b\right)\right]_{t}\in\mathcal{H}_{2}^{t}$,
satisfy (3.1.5).

\noindent (3.1.6) If $Im\left(a\right)^{2}=t\left|b\right|^{2}$ in
$\mathbb{R}$, where $Im\left(a\right)$ is the imaginary part of
$a$ in $\mathbb{C}$, then
\[
spec\left(\left[\left(a,b\right)\right]_{t}\right)=\left\{ x\right\} =\left\{ Re\left(a\right)\right\} \;\mathrm{in\;}\mathbb{R}.
\]
(3.1.7) If $Im\left(a\right)^{2}<t\left|b\right|^{2}$ in $\mathbb{R}$,
then
\[
spec\left(\left[\left(a,b\right)\right]_{t}\right)=\left\{ x\pm\sqrt{tu^{2}+tv^{2}-y^{2}}\right\} \;\mathrm{in\;}\mathbb{R}.
\]
(3.1.8) If $Im\left(a\right)^{2}>t\left|b\right|^{2}$ in $\mathbb{R}$,
then
\[
spec\left(\left[\left(a,b\right)\right]_{t}\right)=\left\{ x\pm i\sqrt{y^{2}-tu^{2}-tv^{2}}\right\} \;\mathrm{in\;}\mathbb{C}\setminus\mathbb{R}.
\]
\end{cor}

\begin{proof}
For $\left(a,b\right)\in\mathbb{H}_{t}$, satisfying (3.1.5), one
has
\[
spec\left(\left[\left(a,b\right)\right]_{t}\right)=\left\{ x\pm i\sqrt{y^{2}-tu^{2}-tv^{2}}\right\} ,
\]
by (3.1.4). So, one can verify that: (i) if $y^{2}-tu^{2}-tv^{2}=0$,
equivalently, if
\[
Im\left(a\right)^{2}=t\left|b\right|^{2}\;\mathrm{in\;}\mathbb{R},
\]
then $spec\left(\left[\left(a,b\right)\right]_{t}\right)=\left\{ x\pm i\sqrt{0}\right\} =\left\{ x\right\} $
in $\mathbb{R}$; (ii) if $y^{2}-tu^{2}-tv^{2}<0,$ equivalently,
if
\[
Im\left(a\right)^{2}<t\left|b\right|^{2}\;\mathrm{in\;}\mathbb{R},
\]
then
\[
x\pm i\sqrt{y^{2}-tu^{2}-tv^{2}}=x\pm i\sqrt{-\left|y^{2}-tu^{2}-tv^{2}\right|},
\]
implying that
\[
x\pm i\sqrt{y^{2}-tu^{2}-tv^{2}}=x\pm i^{2}\sqrt{tu^{2}+tv^{2}-y^{2}},
\]
and hence, 
\[
spec\left(\left[\left(a,b\right)\right]_{t}\right)=\left\{ x\mp\sqrt{tu^{2}+tv^{2}-y^{2}}\right\} \;\mathrm{in\;}\mathbb{R};
\]
and, finally, (iii) if $y^{2}-tu^{2}-tv^{2}>0$, equivalently, if
\[
Im\left(a\right)^{2}>t\left|b\right|^{2}\;\mathrm{in\;}\mathbb{R},
\]
then
\[
spec\left(\left[\left(a,b\right)\right]_{t}\right)=\left\{ x\pm i\sqrt{y^{2}-tu^{2}-tv^{2}}\right\} ,
\]
contained in $\mathbb{C}\setminus\mathbb{R}$.

Therefore, the refined statements (3.1.6), (3.1.7) and (3.1.8) of
the spectrum (3.1.4) of $\left[\left(a,b\right)\right]_{t}$ hold
true.
\end{proof}
By the above corollary, one immediately obtains the following result.
\begin{cor}
Suppose $\left(a,b\right)\in\mathbb{H}_{t}$. If $Im\left(a\right)^{2}\leq t\left|b\right|^{2}$,
then 
\[
spec\left(\left[\left(a,b\right)\right]_{t}\right)\subset\mathbb{R};
\]
meanwhile, if $Im\left(b\right)^{2}>t\left|b\right|^{2}$, then
\[
spec\left(\left[\left(a,b\right)\right]_{t}\right)\subset\left(\mathbb{C}\setminus\mathbb{R}\right),\;\mathrm{in\;}\mathbb{C}.
\]
\end{cor}

\begin{proof}
It is shown by (3.1.6), (3.1.7) and (3.1.8).
\end{proof}
Also, we have the following result.
\begin{thm}
Assume that the fixed scale $t\in\mathbb{R}$ is negative, i.e., $t<0$
in $\mathbb{R}$. If
\[
\left(a,b\right)\in\mathbb{H}_{t},\;\mathrm{with\;}b\neq0\;\mathrm{in\;}\mathbb{C},
\]
then

\hfill{}$spec\left(\left[\left(a,b\right)\right]_{t}\right)\subset\left(\mathbb{C}\setminus\mathbb{R}\right)\;\mathrm{in\;}\mathbb{C}.$\hfill{}(3.1.9)

\medskip{}

\noindent Meanwhile, if $b=0$ in $\mathbb{C}$ for $\left(a,b\right)\in\mathbb{H}_{t}$,
then
\[
a\in\mathbb{R}\Longrightarrow spec\left(\left[\left(a,0\right)\right]_{t}\right)=\left\{ a\right\} \;in\;\mathbb{R},
\]
and\hfill{}(3.1.10)
\[
a\in\mathbb{C}\setminus\mathbb{R}\Longrightarrow spec\left(\left[\left(a,0\right)\right]_{t}\right)=\left\{ a,\overline{a}\right\} \;\mathrm{in\;}\mathbb{C}\setminus\mathbb{R}.
\]
\end{thm}

\begin{proof}
Assume that the scale $t$ is given to be negative in $\mathbb{R}$.
Then, for any $\left(a,b\right)\in\mathbb{H}_{t}$, one immediately
obtains that
\[
Im\left(a\right)^{2}\geq t\left|b\right|^{2},
\]
because the left-hand side, $Im\left(a\right)^{2}$, is nonnegative,
but the right-hand side, $t\left|b\right|^{2}$ is either negative
or zero in $\mathbb{R}$ by the negativity of $t$.

Suppose $b\neq0$ in $\mathbb{C}$, equivalently, $\left|b\right|^{2}>0$,
implying $t\left|b\right|^{2}<0$ in $\mathbb{R}$. Then
\[
Im\left(a\right)^{2}>t\left|b\right|^{2}\;\mathrm{in\;}\mathbb{R}.
\]
Thus, by (3.1.8), the spectra, $spec\left(\left[\left(a,b\right)\right]_{t}\right)$,
of the realizations $\left[\left(a,b\right)\right]_{t}$ of $\left(a,b\right)\in\mathbb{H}_{t}$,
with $b\neq0$, is contained in $\mathbb{C}\setminus\mathbb{R}$.
It proves the relation (3.1.9).

Meanwhile, if $a=Re\left(a\right)$, and $b=0$ in $\mathbb{C}$,
then
\[
0=Im\left(a\right)^{2}\leq0=t\cdot0\;\mathrm{in\;}\mathbb{R},
\]
implying that
\[
spec\left(\left[\left(a,0\right)\right]_{t}\right)\subset\mathbb{R}\;\mathrm{in\;}\mathbb{C},
\]
by (3.1.6). However, if $Im\left(a\right)\neq0$, and $b=0$, then
\[
Im\left(a\right)^{2}>0=t\cdot0\;\mathrm{in\;}\mathbb{R},
\]
and hence,
\[
spec\left(\left[\left(a,0\right)\right]_{t}\right)\subset\left(\mathbb{C}\setminus\mathbb{R}\right)\;\mathrm{in\;}\mathbb{C}.
\]
So, the relation (3.1.10) is proven.
\end{proof}
The above theorem specifies Theorem 19 for the case where $t<0$ in
$\mathbb{R}$, by (3.1.9) and (3.1.10).
\begin{thm}
Assume that $t=0$ in $\mathbb{R}$. If $\left(a,b\right)\in\mathbb{H}_{0}$
with $Im\left(a\right)\neq0$ in $\mathbb{C}$, then

\medskip{}

\hfill{}$spec\left(\left[\left(a,b\right)\right]_{t}\right)\subset\left(\mathbb{C}\setminus\mathbb{R}\right)\;\mathrm{in\;}\mathbb{C}.$\hfill{}(3.1.11)

\medskip{}

\noindent Meanwhile, if $Im\left(a\right)=0$, then

\medskip{}

\hfill{}$spec\left(\left[\left(a,b\right)\right]_{t}\right)\subset\mathbb{R}\;\mathrm{in\;}\mathbb{C}.$\hfill{}(3.1.12)
\end{thm}

\begin{proof}
Suppose the fixed scale $t$ is zero in $\mathbb{R}$. Then, for any
hypercomplex number $\left(a,b\right)\in\mathbb{H}_{0}$, one has
\[
\left[\left(a,b\right)\right]_{0}=\left(\begin{array}{cc}
a & 0\\
\overline{b} & \overline{a}
\end{array}\right)\in\mathcal{H}_{2}^{0},
\]
and hence,
\[
Im\left(a\right)^{2}\geq0=0\cdot\left|b\right|^{2}\;\mathrm{in\;}\mathbb{R}.
\]
In particular, if $Im\left(a\right)\neq0$ in $\mathbb{C}$, then
the above inequality becomes
\[
Im\left(a\right)^{2}>0\;\mathrm{in\;}\mathbb{R},
\]
implying that 
\[
spec\left(\left[\left(a,b\right)\right]_{t}\right)\subset\left(\mathbb{C}\setminus\mathbb{R}\right)\;\mathrm{in\;}\mathbb{C},
\]
by (3.1.8). i.e., for all $\left(a,b\right)\in\mathbb{H}_{0}$, with
$a\in\mathbb{C}$ with $Im\left(a\right)\neq0$, and $b\in\mathbb{C}$
arbitrary, the spectra of the realizations of such $\left(a,b\right)$
are contained in $\mathbb{C}\setminus\mathbb{R}$. It shows the relation
(3.1.11) holds.

Meanwhile, if $Im\left(a\right)=0$ in $\mathbb{C}$, then one has
\[
Im\left(a\right)^{2}=0\geq0=0\cdot\left|b\right|^{2}\;\mathrm{in\;}\mathbb{R}.
\]
So, by (3.1.6), we have
\[
spec\left(\left[\left(a,b\right)\right]_{t}\right)\subset\mathbb{R}\;\mathrm{in}\;\mathbb{C}.
\]
Therefore, the relation (3.1.12) holds true, too.
\end{proof}
The above theorem specifies Theorem 19 for the case where a scale
$t$ is zero in $\mathbb{R}$, by (3.1.11) and (3.1.12).
\begin{thm}
Assume that the fixed scale $t$ is positive in $\mathbb{R}$. Then
the $t$-scaled hypercomplex ring $\mathbb{H}_{t}$ is decomposed
to be
\[
\mathbb{H}_{t}=\mathbb{H}_{t}^{+}\sqcup\mathbb{H}_{t}^{-0},
\]
with\hfill{}(3.1.13)
\[
\mathbb{H}_{t}^{+}=\left\{ \left(a,b\right)\in\mathbb{H}_{t}:Im\left(a\right)^{2}>t\left|b\right|^{2}\right\} ,
\]
and
\[
\mathbb{H}_{t}^{-0}=\left\{ \left(a,b\right)\in\mathbb{H}_{t}:Im\left(a\right)^{2}\leq t\left|b\right|^{2}\right\} ,
\]
where $\sqcup$ is the disjoint union. Moreover, if $\left(a,b\right)\in\mathbb{H}_{t}^{+}$,
then

\medskip{}

\hfill{}$spec\left(\left[\left(a,b\right)\right]_{t}\right)\subset\left(\mathbb{C}\setminus\mathbb{R}\right);$\hfill{}(3.1.14)

\medskip{}

\noindent meanwhile, if $\left(a,b\right)\in\mathbb{H}_{t}^{-0}$,
then

\medskip{}

\hfill{}$spec\left(\left[\left(a,b\right)\right]_{t}\right)\subset\mathbb{R}\;\;\mathrm{in\;\;}\mathbb{C}.$\hfill{}(3.1.15)
\end{thm}

\begin{proof}
Suppose that $t>0$ in $\mathbb{R}$. Then one can decompose the $t$-scaled
hypercomplex ring $\mathbb{H}_{t}$ by
\[
\mathbb{H}_{t}=\mathbb{H}_{t}^{+}\sqcup\mathbb{H}_{t}^{-0},
\]
with\hfill{}(3.1.16)
\[
\mathbb{H}_{t}^{+}=\left\{ \left(a,b\right)\in\mathbb{H}_{t}:Im\left(a\right)^{2}>t\left|b\right|^{2}\right\} ,
\]
and
\[
\mathbb{H}_{t}^{-0}=\left\{ \left(a,b\right)\in\mathbb{H}_{t}:Im\left(a\right)^{2}\leq t\left|b\right|^{2}\right\} ,
\]
set-theoretically. Thus, the partition (3.1.13) holds by (3.1.16).

By Theorem 19 and Corollary 20, if $\left(a,b\right)\in\mathbb{H}_{t}^{+}$,
then
\[
spec\left(\left[\left(a,b\right)\right]_{t}\right)\subset\left(\mathbb{C}\setminus\mathbb{R}\right),
\]
meanwhile, if $\left(a,b\right)\in\mathbb{H}_{t}^{-0}$, then
\[
spec\left(\left[\left(a,b\right)\right]_{t}\right)\subset\mathbb{R},\;\;\mathrm{in\;\;}\mathbb{C}.
\]
So, the relations (3.1.14) and (3.1.15) are proven.
\end{proof}
The above theorem specifies Theorem 19 for the cases where a fixed
scale $t$ is positive in $\mathbb{R}$, by (3.1.14) and (3.1.15),
up to the decomposition (3.1.13). 

In fact, one can realize that, for ``all'' $t\in\mathbb{R}$, the
corresponding $t$-scaled hypercomplex ring $\mathbb{H}_{t}$ is partitioned
to be
\[
\mathbb{H}_{t}=\mathbb{H}_{t}^{+}\sqcup\mathbb{H}_{t}^{-0},
\]
where $\mathbb{H}_{t}^{+}$ and $\mathbb{H}_{t}^{-0}$ are in the
sense of (3.1.13). Especially, Theorems 22, 23 and 24 characterize
the above decomposition case-by-case, based on Theorem 19 and Corollary
20. So, we obtain the following universal spectral properties on $\mathbb{H}_{t}$.
\begin{cor}
Let $t\in\mathbb{R}$ be an arbitrarily fixed scale for $\mathbb{H}_{t}$.
Then 
\[
\mathbb{H}_{t}=\mathbb{H}_{t}^{+}\sqcup\mathbb{H}_{t}^{-0},\;\textrm{set-theoretically,}
\]
where $\left\{ \mathbb{H}_{t}^{+},\:\mathbb{H}_{t}^{-0}\right\} $
is a partition in the sense of (3.1.13) for $t$. Moreover, if $\left(a,b\right)\in\mathbb{H}_{t}^{+}$,
then
\[
spec\left(\left[\left(a,b\right)\right]_{t}\right)\subset\left(\mathbb{C}\setminus\mathbb{R}\right),
\]
meanwhile, if $\left(a,b\right)\in\mathbb{H}_{t}^{-0}$, then
\[
spec\left(\left[\left(a,b\right)\right]_{t}\right)\subset\mathbb{R}\;\;\mathrm{in\;\;}\mathbb{C}.
\]
Especially, if $t<0$, then $\mathbb{H}_{t}^{-0}=\left\{ \left(0,0\right)\right\} $,
equivalently, $\mathbb{H}_{t}^{\times}=\mathbb{H}_{t}^{+}$.
\end{cor}

\begin{proof}
This corollary is nothing but a summary of Theorems 22, 23 and 24.
\end{proof}
It is not hard to check the converses of the statements of Corollary
25 hold true, too. 
\begin{thm}
Let $\mathbb{H}_{t}=\mathbb{H}_{t}^{+}\sqcup\mathbb{H}_{t}^{-0}$
be the fixed $t$-scaled hypercomplex ring for $t\in\mathbb{R}$.

\noindent (3.1.17) $\left(a,b\right)\in\mathbb{H}_{t}^{+}$, if and
only if $spec\left(\left[\left(a,b\right)\right]_{t}\right)\subset\left(\mathbb{C}\setminus\mathbb{R}\right)$.

\noindent (3.1.18) $\left(a,b\right)\in\mathbb{H}_{t}^{-0}$, if and
only if $spec\left(\left[\left(a,b\right)\right]_{t}\right)\subset\mathbb{R}$.
\end{thm}

\begin{proof}
First, assume that $\left(a,b\right)\in\mathbb{H}_{t}^{+}$ in $\mathbb{H}_{t}$.
Then, by Corollary 25,
\[
spec\left(\left[a,b\right]_{t}\right)\subset\left(\mathbb{C}\setminus\mathbb{R}\right).
\]
Now, suppose that
\[
spec\left(\left[a,b\right]_{t}\right)\subset\mathbb{R}\;\;\mathrm{in\;\;}\mathbb{C},
\]
and assume that $\left(a,b\right)\in\mathbb{H}_{t}^{+}$. Then, $\left(a,b\right)$
is contained in $\mathbb{H}_{t}^{-0}$, equivalently, it cannot be
an element of $\mathbb{H}_{t}^{+}$, by (3.1.6), (3.1.7), (3.1.10),
(3.1.12) and (3.1.15). It contradicts our assumption. Therefore,
\[
\left(a,b\right)\in\mathbb{H}_{t}^{+}\Longleftrightarrow spec\left(\left[\left(a,b\right)\right]_{t}\right)\subset\left(\mathbb{C}\setminus\mathbb{R}\right).
\]
Thus, the statement (3.1.17) holds.

By the decomposition (3.1.13), the statement (3.1.18) holds true,
by (3.1.17).
\end{proof}
By the above theorem, we obtain the following result.
\begin{cor}
Let $\mathbb{H}_{t}$ be the $t$-scaled hypercomplex ring for an
arbitrary $t\in\mathbb{R}$, and suppose it is decomposed to be
\[
\mathbb{H}_{t}=\mathbb{H}_{t}^{+}\sqcup\mathbb{H}_{t}^{-0},
\]
as in (3.1.13). Assume that a given element $\left(a,b\right)$ satisfies
the condition (3.1.5). Then

\noindent (3.1.19) $\left(a,b\right)\in\mathbb{H}_{t}^{+}$, if and
only if
\[
spec\left(\left[\left(a,b\right)\right]_{t}\right)=\left\{ x\pm i\sqrt{y^{2}-tu^{2}-tv^{2}}\right\} \subset\left(\mathbb{C}\setminus\mathbb{R}\right).
\]
(3.1.20) $\left(a,b\right)\in\mathbb{H}_{t}^{-0}$, if and only if
either
\[
spec\left(\left[\left(a,b\right)\right]_{t}\right)=\left\{ \begin{array}{ccc}
\left\{ x\right\}  &  & \mathrm{if\;}Im\left(a\right)^{2}=t\left|b\right|^{2}\\
\\
\left\{ x\pm\sqrt{tu^{2}+tv^{2}-y^{2}}\right\}  &  & \mathrm{if\;}Im\left(a\right)^{2}<t\left|b\right|^{2},
\end{array}\right.
\]
in $\mathbb{R}$.
\end{cor}

\begin{proof}
The statement (3.1.19) holds by (3.1.9) and (3.1.17). Meanwhile, the
statement (3.1.20) holds by (3.1.10) and (3.1.18).
\end{proof}
Recall that a Hilbert-space operator $T\in B\left(H\right)$ is self-adjoint,
if $T^{*}=T$ in $B\left(H\right)$, where $T^{*}$ is the adjoint
of $T$ (See Section 5 below). It is well-known that $T$ is self-adjoint,
if and only if its spectrum is contained in $\mathbb{R}$ in $\mathbb{C}$.
So, one obtains the following result.
\begin{prop}
A hypercomplex number $\left(a,b\right)\in\mathbb{H}_{t}^{-0}$ in
$\mathbb{H}_{t}$, if and only if the realization $\left[\left(a,b\right)\right]_{t}\in\mathcal{H}_{2}^{t}$
is self-adjoint ``in $M_{2}\left(\mathbb{C}\right)$.''
\end{prop}

\begin{proof}
($\Rightarrow$) Suppose $\left(a,b\right)\in\mathbb{H}_{t}^{-0}$
in $\mathbb{H}_{t}$. Then $spec\left(\left[\left(a,b\right)\right]_{t}\right)\subset\mathbb{R}$
in $\mathbb{C}$, implying that $\left[\left(a,b\right)\right]_{t}$
is self-adjoint in $M_{2}\left(\mathbb{C}\right)$.

\noindent ($\Leftarrow$) Suppose $\left[\left(a,b\right)\right]_{t}\in\mathcal{H}_{2}^{t}$
is self-adjoint in $M_{2}\left(\mathbb{C}\right)$, and assume that
$\left(a,b\right)\notin\mathbb{H}_{t}^{-0}$, equivalently, $\left(a,b\right)\in\mathbb{H}_{t}^{+}$
in $\mathbb{H}_{t}$. Then,
\[
spec\left(\left[\left(a,b\right)\right]_{t}\right)\subset\left(\mathbb{C}\setminus\mathbb{R}\right)\;\mathrm{in\;}\mathbb{C},
\]
and hence, $\left[\left(a,b\right)\right]_{t}$ is not self-adjoint
in $M_{2}\left(\mathbb{C}\right)$. It contradicts our assumption
that it is self-adjoint.
\end{proof}
Equivalent to the above proposition, one can conclude that $\left(a,b\right)\in\mathbb{H}_{t}^{+}$
in $\mathbb{H}_{t}$, if and only if $\left[\left(a,b\right)\right]_{t}$
is not be self-adjoint in $M_{2}\left(\mathbb{C}\right)$. The self-adjointness
of realizations of hypercomplex numbers would be considered more in
detail in Section 5. 

\subsection{The Scaled-Spectralizations \textmd{$\left\{ \sigma_{t}\right\} _{t\in\mathbb{R}}$}}

In this section, we fix an arbitrary scale $t\in\mathbb{R}$, and
the corresponding hypercomplex ring $\mathbb{H}_{t}$, containing
the $t$-scaled hypercomplex monoid $\mathbb{H}_{t}^{\times}=\left(\mathbb{H}_{t}\setminus\left\{ \left(0,0\right)\right\} ,\:\cdot_{t}\right)$.
Recall that $\mathbb{H}_{t}^{\times}$ is algebraically decomposed
to be
\[
\mathbb{H}_{t}^{\times}=\mathbb{H}_{t}^{inv}\sqcup\mathbb{H}_{t}^{\times sing},
\]
with\hfill{}(3.2.1)
\[
\mathbb{H}_{t}^{inv}=\left\{ \left(a,b\right):\left|a\right|^{2}\neq t\left|b\right|^{2}\right\} ,\;\textrm{the group-part,}
\]
and
\[
\mathbb{H}_{t}^{\times sing}=\left\{ \left(a,b\right):\left|a\right|^{2}=t\left|b\right|^{2}\right\} ,\;\textrm{the semigroup-part,}
\]
as in (2.5.4). Therefore, the $t$-scaled hypercomplex ring is set-theoretically
decomposed to be

\medskip{}

\hfill{}$\mathbb{H}_{t}=\mathbb{H}_{t}^{inv}\sqcup\left\{ \left(0,0\right)\right\} \sqcup\mathbb{H}_{t}^{\times sing}=\mathbb{H}_{t}^{inv}\sqcup\mathbb{H}_{t}^{sing},$\hfill{}(3.2.2)

\medskip{}

\noindent by (3.2.1), where
\[
\mathbb{H}_{t}^{sing}\overset{\textrm{denote}}{=}\left\{ \left(0,0\right)\right\} \sqcup\mathbb{H}_{t}^{\times sing}\;\;\mathrm{in\;\;}\textrm{(3.2.2).}
\]

Also, the ring $\mathbb{H}_{t}$ is spectrally decomposed to be
\[
\mathbb{H}_{t}=\mathbb{H}_{t}^{+}\sqcup\mathbb{H}_{t}^{-0},
\]
with\hfill{}(3.2.3)
\[
\mathbb{H}_{t}^{+}=\left\{ \left(a,b\right):Im\left(a\right)^{2}>t\left|b\right|^{2}\right\} ,
\]
and
\[
\mathbb{H}_{t}^{-0}=\left\{ \left(a,b\right):Im\left(a\right)^{2}\leq t\left|b\right|^{2}\right\} ,
\]
satisfying that: $\left(a,b\right)\in\mathbb{H}_{t}^{+}$ if and only
if $spec\left(\left[\left(a,b\right)\right]_{t}\right)\subset\left(\mathbb{C}\setminus\mathbb{R}\right)$;
meanwhile, $\left(a,b\right)\in\mathbb{H}_{t}^{-0}$ if and only if
$spec\left(\left[\left(a,b\right)\right]_{t}\right)\subset\mathbb{R}$,
by (3.1.19) and (3.1.20).
\begin{cor}
Let $\mathbb{H}_{t}$ be the $t$-scaled hypercomplex ring for $t\in\mathbb{R}$.
Then it is decomposed to be

\medskip{}

\hfill{}$\begin{array}{cc}
\mathbb{H}_{t}= & \left(\mathbb{H}_{t}^{inv}\cap\mathbb{H}_{t}^{+}\right)\sqcup\left(\mathbb{H}_{t}^{inv}\cap\mathbb{H}_{t}^{-0}\right)\\
\\
 & \left(\mathbb{H}_{t}^{sing}\cap\mathbb{H}_{t}^{+}\right)\sqcup\left(\mathbb{H}_{t}^{sing}\cap\mathbb{H}_{t}^{-0}\right),
\end{array}$\hfill{}(3.2.4)

\medskip{}

\noindent set-theoretically.
\end{cor}

\begin{proof}
It is proven by (3.2.2) and (3.2.3).
\end{proof}
Observe now that if $\left(a,0\right)\in\mathbb{H}_{t}$, then
\[
\left[\left(a,0\right)\right]_{t}=\left(\begin{array}{cc}
a & 0\\
0 & \overline{a}
\end{array}\right)\;\mathrm{in\;}\mathcal{H}_{2}^{t},
\]
satisfying\hfill{}(3.2.5)
\[
spec\left(\left[\left(a,0\right)\right]_{t}\right)=\left\{ a,\:\overline{a}\right\} \;\mathrm{in\;}\mathbb{C}.
\]
Indeed, by (3.1.4), if $\left(a,0\right)\in\mathbb{H}_{t}$ satisfying
$a=x+yi\in\mathbb{C}$ with $x,y\in\mathbb{R}$, then
\[
spec\left(\left[\left(a,b\right)\right]_{t}\right)=\left\{ x\pm i\sqrt{y^{2}}\right\} =\left\{ x\pm\left|y\right|i\right\} =\left\{ x\pm yi\right\} ,
\]
implying (3.2.5), where $\left|y\right|$ is the absolute value of
$y$ in $\mathbb{R}$.

Motivated by (3.2.3), (3.2.4) and (3.2.5), we define a certain $\mathbb{C}$-valued
function $\sigma_{t}$ from $\mathbb{H}_{t}$. Define a function,
\[
\sigma_{t}:\mathbb{H}_{t}\rightarrow\mathbb{C},
\]
by\hfill{}(3.2.6)
\[
\sigma_{t}\left(\left(a,b\right)\right)\overset{\textrm{def}}{=}\left\{ \begin{array}{ccc}
a=x+yi &  & \mathrm{if\;}b=0\;\mathrm{in\;}\mathbb{C}\\
\\
x+i\sqrt{y^{2}-tu^{2}-tv^{2}} &  & \mathrm{if\;}b\neq0\;\mathrm{in\;}\mathbb{C},
\end{array}\right.
\]
for all $\left(a,b\right)\in\mathbb{H}_{t}$ satisfying the condition
(3.1.5):
\[
a=x+yi\;\mathrm{and\;}b=u+vi\;\mathrm{in\;}\mathbb{C},
\]
with $x,y,u,v\in\mathbb{R}$ and $i=\sqrt{-1}$. 

Remark that such a morphism $\sigma_{t}$ is indeed a well-defined
function assigning all hypercomplex numbers of $\mathbb{H}_{t}$ to
complex numbers of $\mathbb{C}$. Moreover, by the very definition
(3.2.6), it is surjective. But it is definitely not injective. For
instance, even though
\[
\xi=\left(1+3i,\;-1+i\right)\;\mathrm{and\;\eta=\left(1-3i,\;1-i\right)}
\]
are distinct in $\mathbb{H}_{t}$, one has
\[
\sigma_{t}\left(\xi\right)=1+i\sqrt{9-2t}=\sigma_{t}\left(\eta\right),
\]
by (3.2.6).
\begin{defn}
The surjection $\sigma_{t}:\mathbb{H}_{t}\rightarrow\mathbb{C}$ of
(3.2.6) is called the $t$(-scaled)-spectralization on $\mathbb{H}_{t}$.
The images $\left\{ \sigma_{t}\left(\xi\right)\right\} _{\xi\in\mathbb{H}_{t}}$
are said to be $t$(-scaled)-spectral values. From below, we also
understand each $t$-spectral value $\sigma_{t}\left(\xi\right)\in\mathbb{C}$
of a hypercomplex number $\xi\in\mathbb{H}_{t}$ as a hypercomplex
number $\left(\sigma_{t}\left(\xi\right),0\right)$ in $\mathbb{H}_{t}$.
i.e., such an assigned hypercomplex number $\left(\sigma_{t}\left(\xi\right),0\right)$
from the $t$-spectral value $\sigma_{t}\left(\xi\right)$ of $\xi$
is also called the $t$-spectral value of $\xi$.
\end{defn}

By definition, all $t$-spectral values are not only $\mathbb{C}$-quantities
for many hypercomplex numbers of $\mathbb{H}_{t}$ whose realizations
of $\mathcal{H}_{2}^{t}$ share the same eigenvalues, but also hypercomplex
numbers of $\mathbb{H}_{t}$, whose first coordinates are the value
and the second coordinates are 0.
\begin{defn}
Let $\xi\in\mathbb{H}_{t}$ be a hypercomplex number inducing its
$t$-spectral value $w\overset{\textrm{denote}}{=}\sigma_{t}\left(\xi\right)\in\mathbb{C}$,
also understood to be $\eta=\left(w,0\right)\in\mathbb{H}_{t}$. The
corresponding realization, 
\[
\left[\eta\right]_{t}=\left(\begin{array}{cc}
w & t\cdot0\\
0 & \overline{w}
\end{array}\right)=\left(\begin{array}{ccc}
\sigma_{t}\left(\xi\right) &  & 0\\
\\
0 &  & \overline{\sigma_{t}\left(\xi\right)}
\end{array}\right)\in\mathcal{H}_{2}^{t}
\]
 is called the $t$(-scaled)-spectral form of $\xi$. By $\Sigma_{t}\left(\xi\right)$,
we denote the $t$-spectral form of $\xi\in\mathbb{H}_{t}$.
\end{defn}

Note that the conjugate-notation in Definition 30 is symbolic in the
sense that: if $t>0$, and
\[
\sigma_{t}\left(\xi\right)=1+i\sqrt{1-5t}=1-\sqrt{5t-1},
\]
(and hence, $\sigma_{t}\left(\xi\right)\in\mathbb{R}$), then the
symbol,
\[
\overline{\sigma_{t}\left(\xi\right)}\overset{\textrm{means}}{=}1-i\sqrt{1-5t}=1+\sqrt{5t-1},
\]
in $\mathbb{R}$. i.e., the conjugate-notation in Definition 30 has
a symbolic meaning containing not only the usual conjugate on $\mathbb{C}$,
but also the above computational meaning on $\mathbb{R}$.

\bigskip{}

\noindent $\mathbf{Remark\textrm{-}and\textrm{-}Assumption\;3.2.1.}$
(From below, $\mathbf{RA}\;3.2.1$) The conjugate-notation in Definition
30 is symbolic case-by-case. If the $t$-spectral value $\sigma_{t}\left(\xi\right)$
is in $\mathbb{C}$, then $\overline{\sigma_{t}\left(\xi\right)}$
means the usual conjugate. Meanwhile, if $t$-spectral value 
\[
\sigma_{t}\left(\xi\right)=x+\sqrt{tu^{2}+tv^{2}-y^{2}},
\]
with
\[
tu^{2}+tv^{2}-y^{2}\geq0,\;\mathrm{in}\;\mathbb{R},
\]
then
\[
\overline{\sigma_{t}\left(\xi\right)}=x-\sqrt{tu^{2}+tv^{2}-y^{2}}\;\mathrm{in\;}\mathbb{R},
\]
where $\xi\in\mathbb{H}_{t}$ satisfies the condition (3.1.5). \hfill{}\textifsymbol[ifgeo]{64} 

\bigskip{}

For instance, if $\xi_{1}=\left(-2-i,0\right)\in\mathbb{H}_{t}$,
then the $t$-spectral value is
\[
\sigma_{t}\left(\xi_{1}\right)=-2-i\;\mathrm{in\;}\mathbb{C},
\]
inducing the $t$-spectral form,
\[
\Sigma_{t}\left(\xi_{1}\right)=\left(\begin{array}{ccc}
-2-i &  & 0\\
\\
0 &  & -2+i
\end{array}\right)\;\mathrm{in\;}\mathcal{H}_{2}^{t};
\]
meanwhile, if $\xi_{2}=\left(-2-i,\:1+3i\right)\in\mathbb{H}_{t}$,
then the $t$-spectral value is
\[
w\overset{\textrm{denote}}{=}\sigma_{t}\left(\xi_{2}\right)=-2+i\sqrt{1-10t},
\]
inducing the $t$-spectral form,
\[
\Sigma_{t}\left(\xi_{2}\right)=\left(\begin{array}{cc}
w & 0\\
0 & \overline{w}
\end{array}\right)=\left(\begin{array}{ccc}
-2+i\sqrt{1-10t} &  & 0\\
\\
0 &  & -2-i\sqrt{1-10t}
\end{array}\right),
\]
where $\overline{w}$ is symbolic in the sense of $\mathbf{RA\;3.2.1}$;
if $t\leq0$, then
\[
\Sigma_{t}\left(\xi_{2}\right)=\left(\begin{array}{ccc}
-2+i\sqrt{1-10t} &  & 0\\
\\
0 &  & -2-i\sqrt{1-10t}
\end{array}\right),
\]
meanwhile, if $t>0$, then
\[
\Sigma_{t}\left(\xi_{2}\right)=\left(\begin{array}{ccc}
-2+\sqrt{10t-1} &  & 0\\
\\
0 &  & -2-\sqrt{10t-1}
\end{array}\right),
\]
in $\mathcal{H}_{2}^{t}$.
\begin{defn}
Two hypercomplex numbers $\xi,\eta\in\mathbb{H}_{t}$ are said to
be $t$(-scaled)-spectral-related, if
\[
\sigma_{t}\left(\xi\right)=\sigma_{t}\left(\eta\right)\;\;\mathrm{in\;\;}\mathbb{C},
\]
equivalently,
\[
\Sigma_{t}\left(\xi\right)=\Sigma_{t}\left(\eta\right)\;\;\mathrm{in\;\;}\mathcal{H}_{2}^{t}.
\]
\end{defn}

On the $t$-hypercomplex ring $\mathbb{H}_{t}$, the $t$-spectral
relation of Definition 31 is an equivalent relation. Indeed,
\[
\sigma_{t}\left(\xi\right)=\sigma_{t}\left(\xi\right),\;\forall\xi\in\mathbb{H}_{t};
\]
and if $\xi\;\mathrm{and\;}\eta$ are $t$-spectral related in $\mathbb{H}_{t}$,
then 
\[
\sigma_{t}\left(\xi\right)=\sigma_{t}\left(\eta\right)\Longleftrightarrow\sigma_{t}\left(\eta\right)=\sigma_{t}\left(\xi\right),
\]
and hence, $\eta$ and $\xi$ are $t$-spectral related in $\mathbb{H}_{t}$;
and if $\xi_{1}$ and $\xi_{2}$ are $t$-spectral related, and if
$\xi_{2}$ and $\xi_{3}$ are $t$-spectral related, then
\[
\sigma_{t}\left(\xi_{1}\right)=\sigma_{t}\left(\xi_{2}\right)=\sigma_{t}\left(\xi_{3}\right)\;\mathrm{in\;}\mathbb{C},
\]
and hence, $\xi_{1}$ and $\xi_{3}$ are $t$-spectral related.
\begin{prop}
The $t$-spectral relation on $\mathbb{H}_{t}$ is an equivalence
relation.
\end{prop}

\begin{proof}
The $t$-spectral relation is reflexive, symmetric and transitive
on $\mathbb{H}_{t}$, by the discussion of the very above paragraph.
\end{proof}
Since the $t$-spectral relation is an equivalence relation, each
element $\xi$ of $\mathbb{H}_{t}$ has its equivalence class,
\[
\widetilde{\xi}\overset{\textrm{def}}{=}\left\{ \eta\in\mathbb{H}_{t}:\eta\textrm{ is }t\textrm{-related to }\xi\right\} ,
\]
and hence, the corresponding quotient set,

\medskip{}

\hfill{}$\widetilde{\mathbb{H}_{t}}\overset{\textrm{def}}{=}\left\{ \widetilde{\xi}:\xi\in\mathbb{H}_{t}\right\} ,$\hfill{}(3.2.7)

\medskip{}

\noindent is well-defined to be the set of all equivalence classes.
\begin{thm}
Let $\widetilde{\mathbb{H}_{t}}$ be the quotient set (3.2.7) induced
by the $t$-spectral relation on $\mathbb{H}_{t}$. Then

\medskip{}

\hfill{}$\widetilde{\mathbb{H}_{t}}\;\mathrm{and\;}\mathbb{C}\;\mathrm{are\;equipotent.}$\hfill{}(3.2.8)
\end{thm}

\begin{proof}
It is not difficult to check that, for any $z\in\mathbb{C}$, there
exist $\xi\in\mathbb{H}_{t}$, such that $z=\sigma_{t}\left(\xi\right)$
by the surjectivity of the $t$-spectralization $\sigma_{t}$. It
implies that there exists $\left(z,0\right)\in\mathbb{H}_{t}$, such
that
\[
\widetilde{\left(z,0\right)}=\widetilde{\xi}\;\mathrm{in\;}\widetilde{\mathbb{H}_{t}},\;\mathrm{whenever\;}z=\sigma_{t}\left(\xi\right).
\]
Thus, set-theoretically, we have
\[
\widetilde{\mathbb{H}_{t}}=\left\{ \widetilde{\left(z,0\right)}:z\in\mathbb{C}\right\} \overset{\textrm{equip}}{=}\mathbb{C},
\]
where ``$\overset{\textrm{equip}}{=}$'' means ``being equipotent
(or, bijective) to.'' Therefore, the relation (3.2.8) holds.
\end{proof}
The above equipotence (3.2.8) of the quotient set $\widetilde{\mathbb{H}_{t}}$
of (3.2.7) with the complex numbers $\mathbb{C}$ shows that the set
$\mathbb{C}$ classifies $\mathbb{H}_{t}$, for ``every'' $t\in\mathbb{R}$,
up to the $t$-spectral relation.

\subsection{Similarity on $M_{2}\left(\mathbb{C}\right)$ and The $t$-Scaled-Spectral
Relation on $\mathbb{H}_{t}$}

In Section 3.2, we defined the $t$-spectralization $\sigma_{t}$
on the $t$-scaled hypercomplex ring $\mathbb{H}_{t}$, for a fixed
scale $t$$\in\mathbb{R}$, and it induces the $t$-spectral forms
$\left\{ \Sigma_{t}\left(\xi\right)\right\} _{\xi\in\mathbb{H}_{t}}$
in $\mathcal{H}_{2}^{t}$ as complex diagonal matrices whose main
diagonals are the eigenvalues of the realizations $\left\{ \left[\xi\right]_{t}\right\} _{\xi\in\mathbb{H}_{t}}$,
under the symbolic understanding $\mathbf{RA\;3.2.1}$. Moreover,
$\sigma_{t}$ lets the set $\mathbb{C}$ classify $\mathbb{H}_{t}$
by (3.2.8) under the $t$-spectral relation. 

Independently, we showed in {[}2{]} and {[}3{]} that: on the quaternions
$\mathbb{H}=\mathbb{H}_{-1}$, the ($-1$)-spectral relation, called
the quaternion-spectral relation in {[}2{]} and {[}3{]}, is equivalent
to the similarity ``on $\mathcal{H}_{2}^{-1}$,'' as equivalence
relations. Here, the similarity ``on $\mathcal{H}_{2}^{-1}$'' means
that: the realizations $\left[q_{1}\right]_{-1}$ and $\left[q_{2}\right]_{-1}$
of two quaternions $q_{1},q_{2}\in\mathbb{H}_{-1}$ are similar ``in
$\mathcal{H}_{2}^{-1}$,'' if there exists invertible element $U$
``in $\mathcal{H}_{2}^{-1}$,'' such that
\[
\left[q_{2}\right]_{-1}=U^{-1}\left[q_{1}\right]_{-1}U\;\mathrm{in\;}\mathcal{H}_{2}^{-1}.
\]

\noindent Here, we consider such property for an arbitrary scale $t\in\mathbb{R}$.
Recall that, we showed in {[}2{]} and {[}3{]} that: the $\left(-1\right)$-spectral
form $\Sigma_{-1}\left(\eta\right)$ and the realization $\left[\eta\right]_{-1}$
are similar ``in $\mathcal{H}_{2}^{-1}$,'' for ``all'' quaternions
which are the $\left(-1\right)$-scaled hypercomplex numbers $\eta\in\mathbb{H}_{-1}=\mathbb{H}$.
Are the $t$-spectral relation on $\mathbb{H}_{t}$ and the similarity
on $\mathcal{H}_{2}^{t}$ same as equivalence relations? In conclusion,
the answer is negative in general.

Two matrices $A$ and $B$ of $M_{n}\left(\mathbb{C}\right)$, for
any $n\in\mathbb{N}$, are said to be similar, if there exists an
invertible matrix $U\in M_{n}\left(\mathbb{C}\right)$, such that
\[
B=U^{-1}AU\;\;\mathrm{in\;\;}M_{n}\left(\mathbb{C}\right).
\]
Remember that if two matrices $A$ and $B$ are similar, then (i)
they share the same eigenvalues, (ii) they have the same traces, and
(iii) their determinants are same (e.g., {[}8{]} and {[}9{]}). We
here focus on the fact (iii): the similarity of matrices implies their
identical determinants, equivalently, if
\[
det\left(A\right)\neq det\left(B\right),
\]
then $A$ and $B$ are not similar in $M_{n}\left(\mathbb{C}\right)$.
\begin{defn}
Let $A,B\in\mathcal{H}_{2}^{t}$ be realizations of certain hypercomplex
numbers of $\mathbb{H}_{t}$, for $t\in\mathbb{R}$. They are said
to be similar ``in $\mathcal{H}_{2}^{t}$,'' if there exists an
invertible $U\in\mathcal{H}_{2}^{t}$, such that
\[
B=U^{-1}AU\;\;\mathrm{in\;\;}\mathcal{H}_{2}^{t}.
\]
By abusing notation, we say that two hypercomplex numbers $\xi\;\mathrm{and\;}\eta$
are similar in $\mathbb{H}_{t}$, if their realizations $\left[\xi\right]_{t}$
and $\left[\eta\right]_{t}$ are similar in $\mathcal{H}_{2}^{t}$.
\end{defn}

Let $\left(a,b\right)\in\mathbb{H}_{t}$ be a hypercomplex number
satisfying the condition (3.1.5) and $\left(a,b\right)\neq\left(0,0\right)$.
Then it has
\[
\left[\left(a,b\right)\right]_{t}=\left(\begin{array}{cc}
a & tb\\
\overline{b} & \overline{a}
\end{array}\right)\in\mathcal{H}_{2}^{t},
\]
\[
\sigma_{t}\left(\left(a,b\right)\right)=x+i\sqrt{y^{2}-tu^{2}-tv^{2}}\overset{\textrm{let}}{=}w\in\mathbb{C},
\]
and\hfill{}(3.3.1)
\[
\Sigma_{t}\left(\left(a,b\right)\right)=\left(\begin{array}{cc}
w & 0\\
0 & \overline{w}
\end{array}\right)\in\mathcal{H}_{2}^{t},
\]
where $\overline{w}$ is symbolic under $\mathbf{RA\;3.2.1}$. Observe
that
\[
det\left(\left[\left(a,b\right)\right]_{t}\right)=\left|a\right|^{2}-t\left|b\right|^{2}=\left(x^{2}+y^{2}\right)-t\left(u^{2}+v^{2}\right),
\]
and\hfill{}(3.3.2)
\[
det\left(\Sigma_{t}\left(\left(a,b\right)\right)\right)=\left|w\right|^{2}=x^{2}+\left|y^{2}-tu^{2}-tv^{2}\right|,
\]
by (3.3.1). These computations in (3.3.2) show that, in general, $\left[\left(a,b\right)\right]_{t}$
and $\Sigma_{t}\left(\left(a,b\right)\right)$ are ``not'' similar
``as matrices of $M_{2}\left(\mathbb{C}\right)$,'' and hence, not
similar in $\mathcal{H}_{2}^{t}$. Indeed, for instance, if 
\[
t>0,\;\mathrm{and\;}\left|a\right|^{2}<t\left|b\right|^{2},
\]
then $det\left(\left[\left(a,b\right)\right]_{t}\right)<0$, but $det\left(\Sigma_{t}\left(\left(a,b\right)\right)\right)>0$
in $\mathbb{R}$, by (3.3.2), implying that
\[
det\left(\left[\left(a,b\right)\right]_{t}\right)\neq det\left(\Sigma_{t}\left(\left(a,b\right)\right)\right)\;\textrm{in general,}
\]
showing that $\left[\left(a,b\right)\right]_{t}$ and $\Sigma_{t}\left(\left(a,b\right)\right)$
are not similar in $M_{2}\left(\mathbb{C}\right)$, and hence, they
are not similar in $\mathcal{H}_{2}^{t}$, in general.
\begin{prop}
Let $\left(a,b\right)\in\mathbb{H}_{t}$ be ``nonzero'' hypercomplex
number satisfying $\left|a\right|^{2}<t\left|b\right|^{2}$ in $\mathbb{R}$.
Then the realization $\left[\left(a,b\right)\right]_{t}$ and the
$t$-spectral form $\Sigma_{t}\left(\left(a,b\right)\right)$ are
not similar ``in $\mathcal{H}_{2}^{t}$.''
\end{prop}

\begin{proof}
Suppose $\left(a,b\right)\in\mathbb{H}_{t}$ satisfies $\left(a,b\right)\neq\left(0,0\right)$
and $\left|a\right|^{2}<t\left|b\right|^{2}$, for $t>0$. And assume
that $\left[\left(a,b\right)\right]_{t}$ and $\Sigma_{t}\left(\left(a,b\right)\right)$
are similar in $\mathcal{H}_{2}^{t}$. Since they are assumed to be
similar, their determinants are identically same. However,
\[
det\left(\left[\left(a,b\right)\right]_{t}\right)<0\;\mathrm{and\;}det\left(\Sigma_{t}\left(\left(a,b\right)\right)\right)>0,
\]
by (3.3.2). It contradicts our assumption that they are similar in
$\mathcal{H}_{2}^{t}$.
\end{proof}
The above proposition confirms that the realizations and the corresponding
$t$-spectral forms of a $t$-scaled hypercomplex number are not similar
in $\mathcal{H}_{2}^{t}$, in general.

Consider that, in the quaternions $\mathbb{H}=\mathbb{H}_{-1}$, since
the scale is $t=-1<0$ in $\mathbb{R}$,
\[
det\left(\left[\xi\right]_{-1}\right)=det\left(\Sigma_{-1}\left(\xi\right)\right)\geq0,\;\forall\xi\in\mathbb{H}_{-1},
\]
and it is proven that $\left[\xi\right]_{-1}$ and $\Sigma_{-1}\left(\xi\right)$
are indeed similar in $\mathcal{H}_{2}^{-1}$, for ``all'' $\xi\in\mathbb{H}_{-1}$
in {[}2{]} and {[}3{]}, which motivates a question: if a scale $t<0$
in $\mathbb{R}$, then
\[
det\left(\left[\eta\right]_{t}\right)=det\left(\Sigma_{t}\left(\eta\right)\right)\geq0,\;\forall\eta\in\mathbb{H}_{t},
\]
by (3.3.2); so, are the realizations $\left[\eta\right]$$_{t}$ and
the corresponding $t$-spectral forms $\Sigma_{t}\left(\eta\right)$
similar in $\mathcal{H}_{2}^{t}$ as in the case of $t=-1$?

First of all, we need to recall that if $t<0$, then the $t$-scaled
hypercomplex ring $\mathbb{H}_{t}$ forms a noncommutative field,
since the $t$-scaled hypercomplex monoid $\mathbb{H}_{t}^{\times}$
is a non-abelian group, by (2.4.8). It allows us to use similar techniques
of {[}2{]} and {[}3{]}.

\medskip{}

\noindent $\mathbf{Assumption.}$ In the rest part of this section,
a given scale $t\in\mathbb{R}$ is automatically assumed to be negative
in $\mathbb{R}$. \hfill{}\textifsymbol[ifgeo]{64}

\medskip{}

Assume that $\left(a,0\right)\in\mathbb{H}_{t}$, where $t<0$. Then
\[
\left[\left(a,0\right)\right]_{t}=\left(\begin{array}{cc}
a & 0\\
0 & \overline{a}
\end{array}\right)=\Sigma_{t}\left(\left(a,0\right)\right),
\]
in $\mathcal{H}_{2}^{t}$, since $\sigma_{t}\left(\left(a,0\right)\right)=a$
in $\mathbb{C}$. So, clearly, $\left[\left(a,0\right)\right]_{t}$
and $\Sigma_{t}\left(\left(a,0\right)\right)$ are similar in $\mathcal{H}_{2}^{t}$,
because they are equal in $\mathcal{H}_{2}^{t}$. Indeed, there exist
diagonal matrices with nonzero real entries,
\[
X=\left[\left(x,0\right)\right]_{t}\in\mathcal{H}_{2}^{t},\;\mathrm{with\;}x=x+0i\in\mathbb{C},\;x\neq0,
\]
such that
\[
\left[\left(a,0\right)\right]_{t}=X^{-1}\left(\Sigma_{t}\left(a,0\right)\right)X\;\mathrm{in\;}\mathcal{H}_{2}^{t}.
\]
Thus, we are interested in the cases where $\left(a,b\right)\in\mathbb{H}_{t}$
with $b\in\mathbb{C}^{\times}=\mathbb{C}\setminus\left\{ 0\right\} $.
\begin{lem}
Let $t<0$ in $\mathbb{R}$, and $\left(a,0\right)\in\mathbb{H}_{t}$,
a hypercomplex number. Then the realization $\left[\left(a,0\right)\right]_{t}$
and the $t$-spectral form $\Sigma_{t}\left(\left(a,0\right)\right)$
are identically same in $\mathcal{H}_{2}^{t}$, and hence, they are
similar in $\mathcal{H}_{2}^{t}$. (Remark that, in fact, the scale
$t$ is not necessarily negative in $\mathbb{R}$ here.)
\end{lem}

\begin{proof}
It is proven by the discussion of the very above paragraph. Indeed,
one has
\[
\left[\left(a,0\right)\right]_{t}=\Sigma_{t}\left(\left(a,0\right)\right)\;\mathrm{in\;}\mathcal{H}_{2}^{t},
\]
since $\sigma_{t}\left(\left(a,0\right)\right)=a$ in $\mathbb{C}$.
\end{proof}
Let $h=\left(a,b\right)\in\mathbb{H}_{t}$ with $b\in\mathbb{C}^{\times}$,
satisfying the condition (3.1.5), where $t<0$, having its realization,
\[
\left[h\right]_{t}=\left(\begin{array}{cc}
a & tb\\
\overline{b} & \overline{a}
\end{array}\right)=\left(\begin{array}{ccc}
x+yi &  & t\left(u+vi\right)\\
\\
u-vi &  & x-yi
\end{array}\right),
\]
and its $t$-spectral form,
\[
\Sigma_{t}\left(h\right)=\left(\begin{array}{ccc}
x+i\sqrt{y^{2}-tu^{2}-tv^{2}} &  & 0\\
\\
0 &  & x-i\sqrt{y^{2}-tu^{2}-tv^{2}}
\end{array}\right)\overset{\textrm{let}}{=}\left(\begin{array}{cc}
w & 0\\
0 & \overline{w}
\end{array}\right),
\]
in $\mathcal{H}_{2}^{t}$. Since $t<0$ and $b\neq0$ (by assumption),
the $t$-spectral value $w=\sigma_{t}\left(h\right)$ is a $\mathbb{C}$-quantity
with its conjugate $\overline{w}$. Define now a matrix,
\[
Q_{h}\overset{\textrm{def}}{=}\left(\begin{array}{ccc}
1 &  & t\left(\overline{\frac{w-a}{tb}}\right)\\
\\
\frac{w-a}{tb} &  & 1
\end{array}\right)\;\mathrm{in\;}M_{2}\left(\mathbb{C}\right).
\]
Remark that, by the assumption that $t<0$ and $b\neq0$, this matrix
is well-defined. Furthermore, one can immediately recognize that $Q_{h}\in\mathcal{H}_{2}^{t}$.
i.e.,

\medskip{}

\hfill{}$Q_{h}=\left[\left(1,\:\overline{\left(\frac{w-a}{tb}\right)}\right)\right]_{t}\in\mathcal{H}_{2}^{t}.$\hfill{}(3.3.3)

\medskip{}

\noindent One can find that the element $Q_{h}\in\mathcal{H}_{2}^{t}$
of (3.3.3) is indeed invertible by our negative-scale assumption,
since
\[
det\left(Q_{h}\right)=1-t\left|\frac{w-a}{tb}\right|^{2}\geq1,\;\mathrm{since\;}t<0,
\]
implying that
\[
det\left(Q_{h}\right)\neq0\Longleftrightarrow Q_{h}\textrm{ is invertible in }\mathcal{H}_{2}^{t}.
\]
Observe now that
\[
Q_{h}\Sigma_{t}\left(h\right)=\left(\begin{array}{ccc}
w &  & t\left(\overline{\frac{w^{2}-aw}{tb}}\right)\\
\\
\frac{w^{2}-aw}{tb} &  & \overline{w}
\end{array}\right)
\]

\noindent and\hfill{}(3.3.4)
\[
\left[h\right]_{t}Q_{h}=\left(\begin{array}{ccc}
w &  & t\left(a\left(\overline{\frac{w-a}{tb}}\right)+b\right)\\
\\
\overline{a\left(\frac{w-a}{tb}\right)+b} &  & \overline{w}
\end{array}\right),
\]
in $\mathcal{H}_{2}^{t}$. Now, let's compare the $\left(1,2\right)$-entries
of resulted matrices in (3.3.4). The $\left(1,2\right)$-entry of
the element $Q_{h}\Sigma_{t}\left(h\right)$ is

\medskip{}

$\;\;\;\;$$t\left(\overline{\frac{w^{2}-aw}{tb}}\right)=\overline{\frac{w(w-a)}{b}}=\overline{\frac{\left(x+i\sqrt{y^{2}-tu^{2}-tv^{2}}\right)\left(i\sqrt{y^{2}-tu^{2}-tv^{2}}-yi\right)}{u+vi}}$

\medskip{}

$\;\;\;\;\;\;\;\;\;\;\;\;\;\;\;\;\;\;\;\;$$=\overline{\frac{ix\sqrt{R}-xyi-R+y\sqrt{R}}{u+vi}}$,

\noindent where\hfill{}(3.3.5)
\[
R\overset{\textrm{denote}}{=}y^{2}-tu^{2}-tv^{2}\;\mathrm{in\;}\mathbb{R},
\]

\medskip{}

\noindent and the $\left(1,2\right)$-entry of the matrix $\left[h\right]_{t}Q_{h}$
is

\medskip{}

$\;\;\;\;$$t\left(a\left(\overline{\frac{w-a}{tb}}\right)+b\right)=t\left(\overline{\overline{a}\left(\frac{w-a}{tb}\right)+\overline{b}}\right)$

\medskip{}

$\;\;\;\;\;\;\;\;$$=t\overline{\left(\frac{\overline{a}w-\left|a\right|^{2}+t\left|b\right|^{2}}{tb}\right)}=\overline{\frac{\overline{a}w-\left|a\right|^{2}+t\left|b\right|^{2}}{b}}$

\medskip{}

$\;\;\;\;\;\;\;\;$$=\overline{\frac{\left(x-yi\right)\left(x+i\sqrt{y^{2}-tu^{2}-tv^{2}}\right)-\left(x^{2}+y^{2}\right)-t\left(u^{2}+v^{2}\right)}{u+vi}}$

\medskip{}

$\;\;\;\;\;\;\;\;$$=\overline{\frac{x^{2}+ix\sqrt{R}-xyi+y\sqrt{R}-x^{2}-y^{2}-tu^{2}-tv^{2}}{u+vi}}=\overline{\frac{x^{2}+ix\sqrt{R}-xyi+y\sqrt{R}-x^{2}-R}{u+vi}}$

\medskip{}

$\;\;\;\;\;\;\;\;$$=\overline{\frac{ix\sqrt{R}-xyi-R+y\sqrt{R}}{u+vi}},$\hfill{}(3.3.6)

\medskip{}

\noindent where the $\mathbb{R}$-quantity $R$ is in the sense of
(3.3.5). As one can see in (3.3.5) and (3.3.6), the $\left(1,2\right)$-entries
of $\left[h\right]_{t}Q_{h}$ and $Q_{h}\Sigma_{t}\left(h\right)$
are identically same. i.e.,

\medskip{}

\noindent \hfill{}$Q_{h}\Sigma_{t}\left(h\right)=\left[h\right]_{t}Q_{h}\;\;\mathrm{in\;\;}\mathcal{H}_{2}^{t},$\hfill{}(3.3.7)

\medskip{}

\noindent where the matrix $Q_{h}\in\mathcal{H}_{2}^{t}$ is in the
sense of (3.3.3).
\begin{lem}
Let $t<0$ in $\mathbb{R}$, and let $h=\left(a,b\right)\in\mathbb{H}_{t}$
with $b\in\mathbb{C}^{\times}$. Then the realization $\left[h\right]_{t}$
and the $t$-spectral form $\Sigma_{t}\left(h\right)$ are similar
in $\mathcal{H}_{2}^{t}$. In particular, there exists
\[
q_{h}=\left(1,\;t\left(\overline{\frac{w-a}{tb}}\right)\right)\in\mathbb{H}_{t},
\]
having its realization,
\[
Q_{h}=\left[q_{h}\right]_{t}=\left(\begin{array}{ccc}
1 &  & t\left(\overline{\frac{w-a}{tb}}\right)\\
\\
\frac{w-a}{tb} &  & 1
\end{array}\right)\in\mathcal{H}_{2}^{t},
\]
such that\hfill{}(3.3.8)
\[
\Sigma_{t}\left(h\right)=Q_{h}^{-1}\left[h\right]_{t}Q_{h}\;\;\mathrm{in\;\;}\mathcal{H}_{2}^{t}.
\]
\end{lem}

\begin{proof}
Under the hypothesis, one obtains that
\[
Q_{h}\Sigma_{t}\left(h\right)=\left[h\right]_{t}Q_{b}\;\;\mathrm{in\;\;}\mathcal{H}_{2}^{t},
\]
by (3.3.7). By the invertibility of $Q_{h}$, we have
\[
\Sigma_{t}\left(h\right)=Q_{h}^{-1}\left[h\right]_{t}Q_{h}\;\;\mathrm{in\;\;}\mathcal{H}_{2}^{t},
\]
implying the relation (3.3.8).
\end{proof}
The above lemma shows that if a scale $t$ is negative in $\mathbb{R}$,
then the realization $\left[h\right]_{t}$ and the $t$-spectral form
$\Sigma_{t}\left(h\right)$ are similar in $\mathcal{H}_{2}^{t}$,
whenever $h=\left(a,b\right)\in\mathbb{H}_{t}$ satisfies $b\neq0$
in $\mathbb{C}$.
\begin{thm}
If $t<0$ in $\mathbb{R}$, then every hypercomplex number $h\in\mathbb{H}_{t}$
is similar to its $t$-spectral value $\left(\sigma_{t}\left(h\right),0\right)\in\mathbb{H}_{t}$,
in the sense that:

\medskip{}

\hfill{}$\left[h\right]_{t}$ and $\Sigma_{t}\left(h\right)$ are
similar in $\mathcal{H}_{2}^{t}$.\hfill{}(3.3.9)
\end{thm}

\begin{proof}
Let $h=\left(a,b\right)\in\mathbb{H}_{t}$, for $t<0$. If $b=0$
in $\mathbb{C}$, then $\left[\left(a,0\right)\right]_{t}$ and $\Sigma_{t}\left(\left(a,0\right)\right)$
are similar in $\mathcal{H}_{2}^{t}$, by Lemma 38. Indeed, if $b=0$,
then these matrices are identically same in $\mathcal{H}_{2}^{t}$.
Meanwhile, if $b\neq0$ in $\mathbb{C}$, then $\left[h\right]_{t}$
and $\Sigma_{t}\left(h\right)$ are similar in $\mathcal{H}_{2}^{t}$
by Lemma 39. In particular, if $b\neq0$, then there exists
\[
q_{h}=\left(1,\:\overline{\frac{w-a}{tb}}\right)\in\mathbb{H}_{t},
\]
such that
\[
\Sigma_{t}\left(h\right)=\left[q_{h}\right]_{t}^{-1}\left[h\right]_{t}\left[q_{h}\right]_{t},
\]
in $\mathcal{H}_{2}^{t}$, by (3.3.8). Therefore, if $t<0$, then
$\left[h\right]_{t}$ and $\Sigma_{t}\left(h\right)$ are similar
in $\mathcal{H}_{2}^{t}$, equivalently, two hypercomplex numbers
$h$ and $\left(\sigma_{t}\left(h\right),0\right)$ are similar in
$\mathbb{H}_{t}$, for all $h\in\mathbb{H}_{t}$.
\end{proof}
The above theorem guarantees that the negative-scale condition on
hypercomplex numbers implies the similarity of the realizations and
the scaled-spectral forms of them, just like the quaternionic case
(whose scale is $-1$), shown in {[}2{]} and {[}3{]}.
\begin{thm}
If $t<0$ in $\mathbb{R}$, then the $t$-spectral relation on $\mathbb{H}_{t}$
and the similarity on $\mathbb{H}_{t}$ are same as equivalence relations
on $\mathbb{H}_{t}$. i.e.,

\medskip{}

\hfill{}$t<0\Longrightarrow t\textrm{-spectral relation}\overset{\textrm{equi}}{=}\textrm{similarity on }\mathbb{H}_{t},$\hfill{}(3.3.10)

\medskip{}

\noindent where ``$\overset{\textrm{equi}}{=}$'' means ``being
equivalent to, as equivalence relations.''
\end{thm}

\begin{proof}
Suppose a negative scale $t<0$ is fixed, and let $\mathbb{H}_{t}$
be the corresponding $t$-scaled hypercomplex ring. Assume that two
hypercomplex numbers $h_{1}$ and $h_{2}$ are $t$-spectral related.
Then their $t$-spectral values are identical in $\mathbb{C}$, i.e.,
\[
\sigma_{t}\left(h_{1}\right)=\sigma_{t}\left(h_{2}\right)\overset{\textrm{let}}{=}w\;\mathrm{in\;}\mathbb{C}.
\]
Thus the realizations $\left[h_{1}\right]_{t}$ and $\left[h_{2}\right]_{t}$
are similar to
\[
\Sigma_{t}\left(h_{1}\right)=\left(\begin{array}{cc}
w & 0\\
0 & \overline{w}
\end{array}\right)=\Sigma_{t}\left(h_{2}\right)\overset{\textrm{let}}{=}W,
\]
in $\mathcal{H}_{2}^{t}$, by (3.3.9). i.e., there exist $q_{1},q_{2}\in\mathbb{H}_{t}$
such that
\[
\left[q_{1}\right]_{t}^{-1}\left[h_{1}\right]_{t}\left[q_{1}\right]_{t}=W=\left[q_{2}\right]_{t}^{-1}\left[h_{2}\right]_{t}\left[q_{2}\right]_{t},
\]
in $\mathcal{H}_{2}^{t}$. So, one obtains that
\[
\left[h_{1}\right]_{t}=\left(\left[q_{1}\right]_{t}\left[q_{2}\right]_{t}^{-1}\right)\left[h_{2}\right]_{t}\left(\left[q_{2}\right]_{t}\left[q_{1}\right]_{t}^{-1}\right),
\]
$\Longleftrightarrow$
\[
\left[h_{1}\right]_{t}=\left(\left[q_{2}\right]_{t}\left[q_{1}\right]_{t}^{-1}\right)^{-1}\left[h_{2}\right]_{t}\left(\left[q_{2}\right]_{t}\left[q_{1}\right]_{t}^{-1}\right),
\]
in $\mathcal{H}_{2}^{t}$, implying that $\left[h_{1}\right]_{t}$
and $\left[h_{2}\right]_{t}$ are similar in $\mathcal{H}_{2}^{t}$.
Thus, if $h_{1}$ and $h_{2}$ are $t$-spectral related, then they
are similar in $\mathbb{H}_{t}$.

Conversely, suppose $T_{1}\overset{\textrm{denote}}{=}\left[h_{1}\right]_{t}$
and $T_{2}\overset{\textrm{denote}}{=}\left[h_{2}\right]_{t}$ are
similar in $\mathcal{H}_{2}^{t}$. Then there exists $U\in\mathcal{H}_{2}^{t}$,
such that
\[
T_{1}=U^{-1}T_{2}U\;\;\mathrm{in\;\;}\mathcal{H}_{2}^{t}.
\]
Since the realizations $T_{l}$ and the corresponding $t$-spectral
forms $S_{l}\overset{\textrm{denote}}{=}\Sigma_{t}\left(h_{l}\right)$
are similar by (3.3.9), say,
\[
T_{l}=V_{l}^{-1}S_{l}V_{l}\;\mathrm{in\;}\mathcal{H}_{2}^{t},\;\mathrm{for\;some\;}V_{l}\in\mathcal{H}_{2}^{t},
\]
for all $l=1,2$. Thus,
\[
T_{1}=U^{-1}T_{2}U=U^{-1}\left(V_{2}^{-1}S_{2}V_{2}\right)U,
\]
$\Longleftrightarrow$
\[
V_{1}S_{1}V_{1}^{-1}=T_{1}=\left(V_{2}U\right)^{-1}S_{2}\left(V_{2}U\right),
\]
$\Longleftrightarrow$
\[
S_{1}=V_{1}^{-1}\left(V_{2}U\right)^{-1}S_{2}\left(V_{2}U\right)V_{1},
\]
$\Longleftrightarrow$
\[
S_{1}=\left(V_{2}UV_{1}\right)^{-1}S_{2}\left(V_{2}UV_{1}\right),
\]
and hence, two matrices $S_{1}$ and $S_{2}$ are similar in $\mathcal{H}_{2}^{t}$.
It means that $S_{1}$ and $S_{2}$ share the same eigenvalues. So,
it ie either
\[
S_{1}=\left(\begin{array}{cc}
w & 0\\
0 & \overline{w}
\end{array}\right)=S_{2},
\]
for some $w\in\mathbb{C}$, or
\[
S_{1}=\left(\begin{array}{cc}
w & 0\\
0 & \overline{w}
\end{array}\right),\;\mathrm{and\;}S_{2}=\left(\begin{array}{cc}
\overline{w} & 0\\
0 & w
\end{array}\right),
\]
in $\mathcal{H}_{2}^{t}$. However, by the assumption that $t<0$,
we have
\[
S_{1}=S_{2}\;\;\mathrm{in\;\;}\mathcal{H}_{2}^{t},
\]
by (3.1.8). It shows that, if the realizations $T_{1}$ and $T_{2}$
are similar, then the $t$-spectral forms $S_{1}$ and $S_{2}$ are
identically same in $\mathcal{H}_{2}^{t}$, implying that
\[
\sigma_{t}\left(h_{1}\right)=\sigma_{t}\left(h_{2}\right)\;\;\mathrm{in\;\;}\mathbb{C},
\]
thus $h_{1}$ and $h_{2}$ are $t$-spectral related in $\mathbb{H}_{t}$.

Therefore, the equivalence (3.3.10) between the $t$-spectral relation
and the similarity on $\mathbb{H}_{t}$ holds, whenever $t<0$ in
$\mathbb{R}$.
\end{proof}
The above theorem generalizes the equivalence between the quaternion-spectral
relation, which is the $\left(-1\right)$-spectral relation, and the
similarity on the quaternions $\mathbb{H}_{-1}=\mathbb{H}$ (e.g.,
{[}2{]} and {[}3{]}).

\bigskip{}

\noindent $\mathbf{Discussion.}$ How about the cases where given
scale $t$ are nonnegative in $\mathbb{R}$, i.e., $t\geq0$? One
may need to consider the decomposition (3.2.4),
\[
\begin{array}{cc}
\mathbb{H}_{t}= & \left(\mathbb{H}_{t}^{inv}\cap\mathbb{H}_{t}^{+}\right)\sqcup\left(\mathbb{H}_{t}^{inv}\cap\mathbb{H}_{t}^{-0}\right)\\
\\
 & \left(\mathbb{H}_{t}^{sing}\cap\mathbb{H}_{t}^{+}\right)\sqcup\left(\mathbb{H}_{t}^{sing}\cap\mathbb{H}_{t}^{-0}\right),
\end{array}
\]
of $\mathbb{H}_{t}$, for $t\geq0$, where
\[
\mathbb{H}_{t}^{inv}=\left\{ \left(a,b\right):\left|a\right|^{2}\neq t\left|b\right|^{2}\right\} ,
\]
\[
\mathbb{H}_{t}^{sing}=\left\{ \left(a,b\right):\left|a\right|^{2}=t\left|b\right|^{2}\right\} ,
\]
\[
\mathbb{H}_{t}^{+}=\left\{ \left(a,b\right):Im\left(a\right)^{2}>t\left|b\right|^{2}\right\} ,
\]
and
\[
\mathbb{H}_{t}^{-0}=\left\{ \left(a,b\right):Im\left(a\right)^{2}\leq t\left|b\right|^{2}\right\} ,
\]
block-by-block. In particular, if
\[
h\in\mathbb{H}_{t}^{inv}\cap\mathbb{H}_{t}^{+},
\]
then it ``seems'' that the realization $\left[h\right]_{t}$ and
the $t$-spectral form $\Sigma_{t}\left(h\right)$ are similar in
$\mathcal{H}_{2}^{t}$. The proof ``may'' be similar to the above
proofs for negative scales. We leave this problem for a future project.
\hfill{}\textifsymbol[ifgeo]{64}

\subsection{The $t$-Spectral Mapping Theorem}

In this section, we let a scale $t$ be arbitrary in $\mathbb{R}$,
and let $\mathbb{H}_{t}$ be the $t$-scaled hypercomplex ring. Let
$h=\left(a,b\right)\in\mathbb{H}_{t}$ satisfy the condition (3.1.5),
and suppose it has its $t$-spectral value,
\[
\sigma_{t}\left(h\right)=x+i\sqrt{y^{2}-tu^{2}-tv^{2}}\overset{\textrm{let}}{=}w,
\]
and hence, its $t$-spectral form
\[
\Sigma_{t}\left(h\right)=\left(\begin{array}{cc}
w & 0\\
0 & \overline{w}
\end{array}\right)\;\;\mathrm{in\;\;}\mathcal{H}_{2}^{t},
\]
under $\mathbf{NA\;3.2.1}$.

Now recall that if $n\in\mathbb{N}$, and $A\in M_{n}\left(\mathbb{C}\right)$,
and if
\[
f\in\mathbb{C}[z]\overset{\textrm{def}}{=}\left\{ g:\begin{array}{c}
g=\overset{m}{\underset{k=0}{\sum}}z_{k}z^{k},\;\mathrm{with}\\
z_{1},...,z_{m}\in\mathbb{C},\;\mathrm{for\;}m\in\mathbb{N}
\end{array}\right\} ,
\]
then\hfill{}(3.4.1)
\[
spec\left(f\left(A\right)\right)=\left\{ f\left(w\right):w\in spec\left(A\right)\right\} ,
\]
in $\mathbb{C}$, where $\mathbb{C}[z]$ is the polynomial ring in
a variable $z$ over $\mathbb{C}$, consisting of all polynomials
in $z$ whose coefficients are from $\mathbb{C}$, and
\[
f\left(A\right)=\overset{N}{\underset{k=0}{\sum}}s_{k}A^{k},\;\mathrm{with\;}A^{0}=I_{n},
\]
whenever
\[
f\left(z\right)=\overset{N}{\underset{k=0}{\sum}}s_{k}z^{k}\in\mathbb{C}\left[z\right],\;\mathrm{with\;}s_{1},...,s_{N}\in\mathbb{C},
\]
where $I_{n}$ is the identity matrix of $M_{n}\left(\mathbb{C}\right)$,
by the spectral mapping theorem (e.g., {[}8{]} and {[}9{]}). By (3.4.1),
if $\mathbb{R}[x]$ is the polynomial ring in a variable $x$ over
the real field $\mathbb{R}$, then

\medskip{}

\hfill{}$spec\left(g\left(A\right)\right)=\left\{ g\left(w\right):w\in spec\left(A\right)\right\} \;\mathrm{in\;}\mathbb{C},$\hfill{}(3.4.2)

\medskip{}

\noindent for all $g\in\mathbb{R}[x]$, because $\mathbb{R}[z]$ is
a subring of $\mathbb{C}[z]$ if we identify $x$ to $z$.

It is shown in {[}2{]} and {[}3{]} that, for $f\in\mathbb{C}[z]$,
\[
spec\left(f\left(\left[\xi\right]_{-1}\right)\right)=\left\{ f\left(\sigma_{-1}\left(\xi\right)\right),\;f\left(\overline{\sigma_{-1}\left(\xi\right)}\right)\right\} 
\]
in $\mathbb{C}$, by (3.4.1), but
\[
f\left(\overline{\sigma_{-1}\left(\xi\right)}\right)\neq\overline{f\left(\sigma_{-1}\left(\xi\right)\right)},\;\textrm{in general,}
\]
and hence, even though the spectral mapping theorem (3.4.1) holds
``on $M_{2}\left(\mathbb{C}\right)$, for $\left[\xi\right]_{-1}\in\mathcal{H}_{2}^{-1}$,''
it does not hold ``on $\mathcal{H}_{2}^{-1}$,'' in general. It
demonstrates that, in a similar manner, the spectral mapping theorem
(3.4.1) holds ``on $M_{2}\left(\mathbb{C}\right),$'' but it does
not hold ``on the $t$-scaled realization $\mathcal{H}_{2}^{t}$
of $\mathbb{H}_{t}$,'' for $t\in\mathbb{R}$, because the spectra
of hypercomplex numbers satisfy
\[
spec\left(\left[\eta\right]_{t}\right)=\left\{ w,\overline{w}\right\} ,\;\mathrm{with\;}w=\sigma_{t}\left(\eta\right),
\]
by (3.1.4), for all $\eta\in\mathbb{H}_{t}$ under $\mathbf{RA\;3.2.1}$,
just like the quaternionic case of {[}2{]} and {[}3{]}.

\medskip{}

\noindent $\mathbf{Observation.}$ For an arbitrary scale $t\in\mathbb{R}$,
the spectral mapping theorem (3.4.1) does not hold ``on $\mathcal{H}_{2}^{t}$.''
\hfill{}\textifsymbol[ifgeo]{64}

\medskip{}

However, in {[}2{]} and {[}3{]}, it is proven that, for all $g\in\mathbb{R}[x]$,
one has
\[
spec\left(g\left(\left[\xi\right]_{-1}\right)\right)=\left\{ g\left(\sigma_{t}\left(\xi\right)\right),\:\overline{g\left(\sigma_{t}\left(\xi\right)\right)}\right\} ,
\]
in $\mathbb{C}$, by (3.4.2), since
\[
g\in\mathbb{R}[x]\Longrightarrow g\left(\overline{w}\right)=\overline{g\left(w\right)},\;\forall w\in\mathbb{C}.
\]
It means that the ``restricted'' spectral mapping theorem of (3.4.2)
holds ``on the realization $\mathcal{H}_{2}^{-1}$ of the quaternions
$\mathbb{H}_{-1}$.'' Similarly, we obtain the following result.
\begin{thm}
Let $\xi\in\mathbb{H}_{t}$, realized to be $\left[\xi\right]_{t}\in\mathcal{H}_{2}^{t}$.
Then, for any $g\in\mathbb{R}[x]$,
\[
spec\left(g\left(\left[\xi\right]_{t}\right)\right)=\left\{ g\left(\sigma_{t}\left(\xi\right)\right),\:\overline{g\left(\sigma_{t}\left(\xi\right)\right)}\right\} ,
\]
i.e.,\hfill{}(3.4.3)
\[
spec\left(g\left(\left[\xi\right]_{t}\right)\right)=\left\{ g\left(w\right):w\in spec\left(\left[\xi\right]_{t}\right)\right\} \;\mathrm{in\;}\mathbb{C},\;\forall t\in\mathbb{R}.
\]
\end{thm}

\begin{proof}
By (3.1.4) and (3.2.6), if $\xi\in\mathbb{H}_{t}$, then
\[
spec\left(\left[\xi\right]_{t}\right)=\left\{ w,\overline{w}\right\} ,\;\mathrm{with\;}w=\sigma_{t}\left(\xi\right),
\]
in $\mathbb{C}$ (under the symbolic understanding of $\mathbf{RA\;3.2.1}$).
For any $g=\overset{N}{\underset{k=1}{\sum}}s_{k}x^{k}\in\mathbb{R}[x]$,
with $s_{1},...,s_{N}\in\mathbb{R}$, and $N\in\mathbb{N}$, one has
that

\medskip{}

\hfill{}$g\left(\overline{w}\right)=\overset{N}{\underset{k-1}{\sum}}s_{k}\overline{w}^{k}=\overset{N}{\underset{k=1}{\sum}}\overline{s_{k}w^{k}}=\overline{\overset{N}{\underset{k=1}{\sum}}s_{k}w^{k}}=\overline{g\left(w\right)},$\hfill{}(3.4.4)

\medskip{}

\noindent in $\mathbb{C}$. It implies that
\[
spec\left(g\left(\left[\xi\right]_{t}\right)\right)=\left\{ g\left(w\right),g\left(\overline{w}\right)\right\} =\left\{ g\left(w\right),\overline{g\left(w\right)}\right\} ,
\]
in $\mathbb{C}$, by (3.4.2) and (3.4.4). Therefore, the relation
(3.4.3) holds true.
\end{proof}
One may call the relation (3.4.3), the hypercomplex-spectral mapping
theorem, since it holds for all scales $t\in\mathbb{R}$.

\section{The Usual Adjoint on $\mathcal{H}_{2}^{t}$ in $M_{2}\left(\mathbb{C}\right)$}

In this section, we consider how the usual adjoint on $M_{2}\left(\mathbb{C}\right)=B\left(\mathbb{C}^{2}\right)$
acts on the $t$-scaled realization $\mathcal{H}_{2}^{t}$ of the
$t$-scaled hypercomplex numbers. Throughout this section, we fix
an arbitrary scale $t\in\mathbb{R}$, and the corresponding $t$-scaled
hypercomplex ring $\mathbb{H}_{t}$ realized to be $\mathcal{H}_{2}^{t}$
in $M_{2}\left(\mathbb{C}\right)$ under the representation $\Pi_{t}=\left(\mathbb{C}^{2},\pi_{t}\right)$.
Recall that every Hilbert-space operator $T$ acting on a Hilbert
space $H$ has its unique adjoint $T^{*}$ on $H$. Especially, if
$T\in M_{n}\left(\mathbb{C}\right)=B\left(\mathbb{C}^{n}\right)$,
for $n\in\mathbb{N}$, is a matrix which is an operator on $\mathbb{C}^{n}$,
then its adjoint $T^{*}$ is determined to be the conjugate-transpose
of $T$ in $M_{n}\left(\mathbb{C}\right)$. For instance,
\[
T=\left(\begin{array}{cc}
a_{11} & a_{12}\\
a_{21} & a_{22}
\end{array}\right)\in M_{2}\left(\mathbb{C}\right)\Longleftrightarrow T^{*}=\left(\begin{array}{cc}
\overline{a_{11}} & \overline{a_{21}}\\
\overline{a_{12}} & \overline{a_{22}}
\end{array}\right)\in M_{2}\left(\mathbb{C}\right).
\]
It says that, if we understand our $t$-scaled realization $\mathcal{H}_{2}^{t}$
as a sub-structure of $M_{2}\left(\mathbb{C}\right)$, then each hypercomplex
number $\left(a,b\right)\in\mathbb{H}_{t}$ assigns a unique adjoint
$\left[\left(a,b\right)\right]_{t}^{*}$ of the realization $\left[\left(a,b\right)\right]_{t}$
``in $M_{2}\left(\mathbb{C}\right)$.''

Let $\left(a,b\right)\in\mathbb{H}_{t}$ realized to be 
\[
\left[\left(a,b\right)\right]_{t}=\left(\begin{array}{cc}
a & tb\\
\overline{b} & \overline{a}
\end{array}\right)\in\mathcal{H}_{2}^{t}.
\]
Then, as a matrix of $M_{2}\left(\mathbb{C}\right)$, this realization
has its adjoint,
\[
\left[\left(a,b\right)\right]_{t}^{*}=\left(\begin{array}{cc}
\overline{a} & b\\
t\overline{b} & a
\end{array}\right)\;\mathrm{in\;}M_{2}\left(\mathbb{C}\right).
\]
It shows that the usual adjoint (conjugate-transpose) of $\left[\left(a,b\right)\right]_{t}$
is not contained ``in $\mathcal{H}_{2}^{t}$,'' in general. In particular,
if
\[
t^{2}\neq1\Longleftrightarrow\mathrm{either\;}t\neq1\;\mathrm{or\;}t\neq-1,\;\mathrm{in\;}\mathbb{R},
\]
then
\[
\left[\left(a,b\right)\right]_{t}\notin\mathcal{H}_{2}^{t}\;\mathrm{in\;general.}
\]

\begin{thm}
The scale $t\in\mathbb{R}$ satisfies that $t^{2}=1$ in $\mathbb{R}$,
if and only if the adjoint of every realization of a hypercomplex
number $\mathbb{H}_{t}$ is contained in $\mathcal{H}_{2}^{t}$. i.e., 

\medskip{}

\hfill{}either $t=1$, or $t=-1$$\Longleftrightarrow\left[\xi\right]_{t}^{*}\in\mathcal{H}_{2}^{t},\;\forall\xi\in\mathbb{H}_{t}$.
\hfill{}(4.1)
\end{thm}

\begin{proof}
For an arbitrary scale $t\in\mathbb{R}$, if $\left(a,b\right)\in\mathbb{H}_{t}$,
then
\[
\left[\left(a,b\right)\right]_{t}^{*}=\left(\begin{array}{cc}
\overline{a} & b\\
t\overline{b} & a
\end{array}\right)\;\mathrm{in\;}M_{2}\left(\mathbb{C}\right).
\]
($\Rightarrow$) Assume that either $t=1$, or $t=-1$, equivalently,
suppose $t^{2}=1$ in $\mathbb{R}$. Then
\[
\left[\left(a,b\right)\right]_{t}^{*}=\left(\begin{array}{cc}
\overline{a} & b\\
t\overline{b} & a
\end{array}\right)=\left(\begin{array}{cc}
\overline{a} & t\left(\frac{b}{t}\right)\\
t^{2}\overline{\left(\frac{b}{t}\right)} & a
\end{array}\right)=\left(\begin{array}{cc}
\overline{a} & t\left(\frac{b}{t}\right)\\
\overline{\left(\frac{b}{t}\right)} & a
\end{array}\right),
\]
contained in $\mathcal{H}_{2}^{t}$. So, if either $t=1$, or $t=-1$,
then $\left[\left(a,b\right)\right]_{t}^{*}\in\mathcal{H}_{2}^{t}$,
for all $\left(a,b\right)\in\mathbb{H}_{t}$. Moreover, in such a
case,

\medskip{}

\hfill{}$\left[\left(a,b\right)\right]_{t}^{*}=\left[\left(\overline{a},\:\frac{b}{t}\right)\right]_{t}\;\mathrm{in\;}\mathcal{H}_{2}^{t}.$\hfill{}(4.2)

\medskip{}

\noindent ($\Leftarrow$) Assume now that $t^{2}\neq1$ in $\mathbb{R}$.
Then the adjoint $\left[\left(a,b\right)\right]_{t}^{*}$ of $\left[\left(a,b\right)\right]_{t}$
is identical to the matrix,
\[
\left[\left(a,b\right)\right]_{t}^{*}=\left(\begin{array}{cc}
\overline{a} & b\\
t\overline{b} & a
\end{array}\right)\;\mathrm{in\;}M_{2}\left(\mathbb{C}\right),
\]
which ``can'' be
\[
\left(\begin{array}{cc}
\overline{a} & t\left(\frac{b}{t}\right)\\
t^{2}\left(\overline{\frac{b}{t}}\right) & a
\end{array}\right)\;\mathrm{in\;}\mathcal{H}_{2}^{t}.
\]
However, by the assumption that $t^{2}\neq1$, the adjoint $\left[\left(a,b\right)\right]_{t}^{*}$
is not contained in $\mathcal{H}_{2}^{t}$, in general. In particular,
if $b\neq0$ in $\mathbb{C}$, then the adjoint $\left[\left(a,b\right)\right]_{t}^{*}\notin\mathcal{H}_{2}^{t}$
in $M_{2}\left(\mathbb{C}\right)$, i.e.,

\medskip{}

\hfill{}$t^{2}\neq1\;\mathrm{and\;}b\neq0\;\mathrm{in\;}\mathbb{C}\Longrightarrow\left[\left(a,b\right)\right]_{t}^{*}\in\left(M_{2}\left(\mathbb{C}\right)\setminus\mathcal{H}_{2}^{t}\right).$\hfill{}(4.3)

\medskip{}

Therefore, the characterization (4.1) holds by (4.2) and (4.3).
\end{proof}
Note that, if $t=-1$, then $\mathbb{H}_{-1}$ is the quaternions;
and if $t=1$, then $\mathbb{H}_{1}$ is the bicomplex numbers. The
above theorem shows that, only when the scaled hypercomplex ring $\mathbb{H}_{t}$
is either the quaternions $\mathbb{H}_{-1}$, or the bicomplex numbers
$\mathbb{H}_{1}$, the usual adjoint ($*$) is closed on $\mathcal{H}_{2}^{t}$,
as a well-defined unary operation, by (4.1).

\section{Free Probability on $\mathbb{H}_{t}$}

In this section, we establish a universal free-probabilistic model
on our $t$-scaled hypercomplex ring $\mathbb{H}_{t}$, for ``every''
scale $t\in\mathbb{R}$. First, recall that, on $M_{2}\left(\mathbb{C}\right)$,
we have the usual trace $\mathit{tr}$, defined by
\[
tr\left(\left(\begin{array}{cc}
a_{11} & a_{12}\\
a_{21} & a_{22}
\end{array}\right)\right)=a_{11}+a_{22},
\]
for all $\left(\begin{array}{cc}
a_{11} & a_{12}\\
a_{21} & a_{22}
\end{array}\right)\in M_{2}\left(\mathbb{C}\right)$; and the normalized trace $\tau$,
\[
\tau=\frac{1}{2}tr\;\;\mathrm{on\;\;}M_{2}\left(\mathbb{C}\right).
\]
i.e., we have two typical free-probabilistic models,
\[
\left(M_{2}\left(\mathbb{C}\right),\:tr\right)\;\mathrm{and\;}\left(M_{2}\left(\mathbb{C}\right),\:\tau\right).
\]

\subsection{Free Probability}

For more about free probability theory, see e.g., {[}19{]} and {[}22{]}.
Let $A$ be an noncommutative algebra over $\mathbb{C}$, and $\varphi:A\rightarrow\mathbb{C}$,
a linear functional on $A$. Then the pair $\left(A,\varphi\right)$
is called a (noncommutative) free probability space. By definition,
free probability spaces are the noncommutative version of classic
measure spaces $\left(X,\mu\right)$ consisting of a set $X$ and
a measure $\mu$ on the $\sigma$-algebra of $X$. As in measure theory,
the (noncommutative) free probability on $\left(A,\varphi\right)$
is dictated by the linear functional $\varphi$. Meanwhile, if $\left(A,\varphi\right)$
is unital in the sense that (i) the unity $1_{A}$ of $A$ exists,
and (ii) $\varphi\left(1_{A}\right)=1$, then it is called a unital
free probability space. These unital free probability spaces are the
noncommutative analogue of classical probability spaces $\left(Y,\rho\right)$
where the given measures $\rho$ are the probability measures satisfying
$\rho\left(Y\right)=1$.

If $A$ is a topological algebra, and if $\varphi$ is bounded (and
hence, continuous under linearity), then the corresponding free probability
space $\left(A,\varphi\right)$ is said to be a topological free probability
space. Similarly, if $A$ is a topological $*$-algebra equipped with
the adjoint ($*$), then the topological free probability space $\left(A,\varphi\right)$
is said to be a topological (free) $*$-probability space. More in
detail, if $A$ is a $C^{*}$-algebra, or a von Neumann algebra, or
a Banach $*$-algebra, we call $\left(A,\varphi\right)$, a $C^{*}$-probability
space, respectively, a $W^{*}$-probability space, respectively, a
Banach $*$-probability space, etc.. For our main purposes, we focus
on $C^{*}$-probability spaces from below.

If $\left(A,\varphi\right)$ is a $C^{*}$-probability space, and
$a\in A$, then the algebra-element $a$ is said to be a free random
variable of $\left(A,\varphi\right)$. For any arbitrarily fixed free
random variables $a_{1},...,a_{s}\in\left(A,\varphi\right)$ for $s\in\mathbb{N}$,
one can get the corresponding free distribution of $a_{1},...,a_{s}$,
characterized by the joint free moments,
\[
\varphi\left(\overset{n}{\underset{l=1}{\prod}}a_{i_{l}}^{r_{i}}\right)=\varphi\left(a_{i_{1}}^{r_{1}}a_{i_{2}}^{r_{2}}...a_{i_{n}}^{r_{n}}\right),
\]
for all $\left(i_{1},...,i_{n}\right)\in\left\{ 1,...,s\right\} ^{n}$
and $\left(r_{1},...,r_{n}\right)\in\left\{ 1,*\right\} ^{n}$, for
all $n\in\mathbb{N}$, where $a_{l}^{*}$ are the adjoints of $a_{l}$,
for all $l=1,...,s$. For instance, if $a\in\left(A,\varphi\right)$
is a free random variable, then the free distribution of $a$ is fully
characterized by the joint free moments of $\left\{ a,a^{*}\right\} $,
\[
\varphi\left(\overset{n}{\underset{l=1}{\prod}}a^{r_{l}}\right)=\varphi\left(a^{r_{1}}a^{r_{2}}...a^{r_{n}}\right),
\]
for all $\left(r_{1},...,r_{n}\right)\in\left\{ 1,*\right\} ^{n}$,
for all $n\in\mathbb{N}$ (e.g., {[}19{]} and {[}22{]}). So, if a
free random variable $a\in\left(A,\varphi\right)$ is self-adjoint
in the sense that: $a^{*}=a$ in $A$, then the free distribution
of $a$ is determined by the free-moment sequence,
\[
\left(\varphi\left(a^{n}\right)\right)_{n=1}^{\infty}=\left(\varphi\left(a\right),\varphi\left(a^{2}\right),\varphi\left(a^{3}\right),...\right)
\]
(e.g., {[}19{]} and {[}22{]}).

\subsection{Free-Probabilistic Models Induced by $\mathbb{H}_{t}$}

By identifying the $t$-scaled hypercomplex ring $\mathbb{H}_{t}$
and its realization $\mathcal{H}_{2}^{t}$ as the same ring, we identify
the $t$-scaled hypercomplex monoid $\mathbb{H}_{t}^{\times}$ and
its realization $\mathcal{H}_{2}^{t\times}$ as the same monoid. As
a subset in $M_{2}\left(\mathbb{C}\right)$, we define a subset,
\[
\mathcal{H}_{2}^{t\times}(*)\overset{\textrm{def}}{=}\left\{ \left[\xi\right]_{t}^{*}\in M_{2}\left(\mathbb{C}\right):\xi\in\mathbb{H}_{t}^{\times}\right\} ,
\]
i.e.,\hfill{}(5.2.1)
\[
\mathcal{H}_{2}^{t\times}\left(*\right)=\left\{ \left(\begin{array}{cc}
\overline{a} & b\\
t\overline{b} & a
\end{array}\right)\in M_{2}\left(\mathbb{C}\right):\left(a,b\right)\in\mathbb{H}_{t}^{\times}\right\} ,
\]
by the subset of all adjoints of realizations in $\mathcal{H}_{2}^{\times t}$.
Indeed,
\[
\left[\left(a,b\right)\right]_{t}^{*}=\left(\begin{array}{cc}
a & tb\\
\overline{b} & \overline{a}
\end{array}\right)^{*}=\left(\begin{array}{cc}
\overline{a} & b\\
t\overline{b} & a
\end{array}\right)\;\mathrm{in\;}M_{2}\left(\mathbb{C}\right).
\]
As we have seen in Section 4, the adjoint is not closed on $\mathcal{H}_{2}^{t}$
in general, and hence,
\[
\mathcal{H}_{2}^{t\times}(*)\neq\mathcal{H}_{2}^{t\times}\;\mathrm{in\;}M_{2}\left(\mathbb{C}\right),
\]
in general. In particular, the scale $t$ satisfies $t^{2}\neq1$
in $\mathbb{R}$, if and only if the above non-equality holds in $M_{2}\left(\mathbb{C}\right)$,
by (4.1). Now, let
\[
\mathcal{H}_{2}^{t\times}\left(1,*\right)\overset{\textrm{denote}}{=}\mathcal{H}_{2}^{t\times}\cup\mathcal{H}_{2}^{t\times}(*),
\]
i.e.,\hfill{}(5.2.2)
\[
\mathcal{H}_{2}^{t\times}\left(1,*\right)=\left\{ \left(\begin{array}{cc}
a & tb\\
\overline{b} & \overline{a}
\end{array}\right),\:\left(\begin{array}{cc}
\overline{a} & b\\
t\overline{b} & a
\end{array}\right):\left(a,b\right)\in\mathbb{H}_{t}^{\times}\right\} ,
\]
in $M_{2}\left(\mathbb{C}\right)$, set-theoretically. By (4.1), (5.2.1)
and (5.2.2),
\[
\mathcal{H}_{2}^{t\times}\left(1,*\right)\supsetneqq\mathcal{H}_{2}^{t\times}\;\mathrm{in\;}M_{2}\left(\mathbb{C}\right),\;\textrm{in general.}
\]

Define now the $C^{*}$-algebra $\mathfrak{H}_{2}^{t}$ by the $C^{*}$-subalgebra
of $M_{2}\left(\mathbb{C}\right)$ generated by the set $\mathcal{H}_{2}^{t\times}(1,*)$
of (5.2.2). i.e.,

\medskip{}

\hfill{}$\mathfrak{H}_{2}^{t}\overset{\textrm{denote}}{=}C^{*}\left(\mathcal{H}_{2}^{t\times}\right)\overset{\textrm{def}}{=}\overline{\mathbb{C}\left[\mathcal{H}_{2}^{t\times}\left(1,*\right)\right]},$\hfill{}(5.2.3)

\medskip{}

\noindent in $M_{2}\left(\mathbb{C}\right)$, where $C^{*}\left(Z\right)$
means the $C^{*}$-subalgebra of $B\left(\mathbb{C}^{2}\right)$ generated
by the subset $Z$ and their adjoints, and $\mathbb{C}[X]$ is the
(pure-algebraic) algebra (over $\mathbb{C}$) generated by a subset
$X$ of $M_{2}\left(\mathbb{C}\right)$, and $\overline{Y}$ means
the operator-norm-topology closure of a subset $Y$ of the operator
algebra $M_{2}\left(\mathbb{C}\right)=B\left(\mathbb{C}^{2}\right)$,
which is a $C^{*}$-algebra over $\mathbb{C}$.
\begin{defn}
The $C^{*}$-algebra $\mathfrak{H}_{2}^{t}$ of (5.2.3), generated
by the $t$-scaled hypercomplex monoid $\mathbb{H}_{t}^{\times}\overset{\textrm{monoid}}{=}\mathcal{H}_{2}^{t\times}$,
is called the $t$-scaled(-hypercomplex)-monoidal $C^{*}$-algebra
of $\mathbb{H}_{t}^{\times}$ (or, of $\mathbb{H}_{t}$).
\end{defn}

Clearly, by the definition (5.2.3), the $t$-scaled-monoidal $C^{*}$-algebra
$\mathfrak{H}_{2}^{t}$ is well-determined in $M_{2}\left(\mathbb{C}\right)$.
So, the usual trace $tr$ and the normalized trace $\tau$ on $M_{2}\left(\mathbb{C}\right)$
are well-defined on $\mathfrak{H}_{2}^{t}$. i.e., we have two trivial
free-probabilistic models of $\mathfrak{H}_{2}^{t}$,
\[
\left(\mathfrak{H}_{2}^{t},\:tr\right)\;\;\mathrm{and\;\;}\left(\mathfrak{H}_{2}^{t},\:\tau\right),
\]
as $C^{*}$-probability spaces (e.g., see Section 5.1). Note that
such free-probabilistic structures are independent from the choice
of the scales $t\in\mathbb{R}$.

Observe that, if $\left(\begin{array}{cc}
\overline{a_{l}} & b_{l}\\
t\overline{b_{l}} & a_{l}
\end{array}\right)\in\mathcal{H}_{2}^{t\times}(*)$ in $\mathfrak{H}_{2}^{t}$, for $l=1,2$, then
\[
\left(\begin{array}{cc}
\overline{a_{1}} & b_{1}\\
t\overline{b_{1}} & a_{1}
\end{array}\right)\left(\begin{array}{cc}
\overline{a_{2}} & b_{2}\\
t\overline{b_{2}} & a_{2}
\end{array}\right)=\left(\begin{array}{ccc}
\overline{a_{1}a_{2}}+tb_{1}\overline{b_{2}} &  & \overline{a_{1}}b_{2}+b_{1}a_{2}\\
\\
t\left(\overline{b_{1}a_{2}}+a_{1}\overline{b_{2}}\right) &  & t\overline{b_{1}}b_{2}+a_{1}a_{2}
\end{array}\right),
\]
identifying to be\hfill{}(5.2.4)
\[
\left(\begin{array}{ccc}
\overline{a_{1}a_{2}+t\overline{b_{1}}b_{2}} &  & b_{1}a_{2}+\overline{a_{1}}b_{2}\\
\\
t\left(\overline{b_{1}a_{2}+\overline{a_{1}}b_{2}}\right) &  & a_{1}a_{2}+t\overline{b_{1}}b_{2}
\end{array}\right)\;\mathrm{in\;}\mathfrak{H}_{2}^{t}.
\]
Therefore,
\[
\left(\begin{array}{cc}
\overline{a_{1}} & b_{1}\\
t\overline{b_{1}} & a_{1}
\end{array}\right)\left(\begin{array}{cc}
\overline{a_{2}} & b_{2}\\
t\overline{b_{2}} & a_{2}
\end{array}\right)\in\mathcal{H}_{2}^{t\times}(*),\;\textrm{too}.
\]
i.e., the matricial multiplication is closed on the set $\mathcal{H}_{2}^{t\times}(*)$
of (5.2.2), by (5.2.4). In fact, under the closed-ness (5.2.4), the
algebraic pair,
\[
\mathcal{H}_{2}^{t\times}(*)\overset{\textrm{denote}}{=}\left(\mathcal{H}_{2}^{t\times}(*),\;\cdot\right),
\]
forms a monoid with its identity $I_{2}$. So, the generating set
$\mathcal{H}_{2}^{t\times}(1,*)$ of the $t$-scaled-monoidal $C^{*}$-algebra
$\mathfrak{H}_{2}^{t}$ is the set-theoretical union of two monoids
$\mathcal{H}_{2}^{t\times}$ and $\mathcal{H}_{2}^{t\times}(*)$,
under the matricial multiplication. Note, however, that the matricial
multiplication is not closed on the generating set $\mathcal{H}_{2}^{t\times}(1,*)$
of (5.2.2). Indeed, if
\[
\left(\begin{array}{cc}
a_{1} & tb_{1}\\
\overline{b_{1}} & \overline{a_{1}}
\end{array}\right)\in\mathcal{H}_{2}^{t\times},\;\left(\begin{array}{cc}
\overline{a_{2}} & b_{2}\\
t\overline{b_{2}} & a_{2}
\end{array}\right)\in\mathcal{H}_{2}^{t\times}(*)
\]
in $\mathfrak{H}_{2}^{t}$, then
\[
\left(\begin{array}{cc}
a_{1} & tb_{1}\\
\overline{b_{1}} & \overline{a_{1}}
\end{array}\right)\left(\begin{array}{cc}
\overline{a_{2}} & b_{2}\\
t\overline{b_{2}} & a_{2}
\end{array}\right)=\left(\begin{array}{ccc}
a_{1}\overline{a_{2}}+t^{2}b_{1}\overline{b_{2}} &  & a_{1}b_{2}+ta_{2}b_{1}\\
\\
\overline{a_{2}b_{1}}+t\overline{a_{1}b_{2}} &  & \overline{b_{1}}b_{2}+\overline{a_{1}}a_{2}
\end{array}\right),
\]
and\hfill{}(5.2.5)
\[
\left(\begin{array}{cc}
\overline{a_{2}} & b_{2}\\
t\overline{b_{2}} & a_{2}
\end{array}\right)\left(\begin{array}{cc}
a_{1} & tb_{1}\\
\overline{b_{1}} & \overline{a_{1}}
\end{array}\right)=\left(\begin{array}{ccc}
a_{1}\overline{a_{2}}+\overline{b_{1}}b_{2} &  & tb_{1}\overline{a_{2}}+\overline{a_{1}}b_{2}\\
\\
ta_{1}\overline{b_{2}}+\overline{b_{1}}a_{2} &  & t^{2}b_{1}\overline{b_{2}}+\overline{a_{1}}a_{2}
\end{array}\right),
\]
in $\mathfrak{H}_{2}^{t}$. However, the resulted products of (5.2.5),
contained in $\mathfrak{H}_{2}^{t}$, are not contained in $\mathcal{H}_{2}^{t\times}(1,*)$,
in general.

\bigskip{}

\noindent $\mathbf{Observation.}$ By (5.2.4) and (5.2.5), one can
realize that: (i) if $A,B\in\mathcal{H}_{2}^{t\times}$, then $AB\in\mathcal{H}_{2}^{t\times}$,
(ii) if $C,D\in\mathcal{H}_{2}^{t\times}(*)$, then $CD\in\mathcal{H}_{2}^{t\times}(*)$,
and (iii) if $T,S\in\mathcal{H}_{2}^{t\times}(1,*)$, then $TS\notin\mathcal{H}_{2}^{t\times}(1,*)$,
in general, as elements of the $t$-scaled-monoidal $C^{*}$-algebra
$\mathfrak{H}_{2}^{t}$. Even though the non-closed rule (iii) is
satisfied ``on $\mathcal{H}_{2}^{t}\left(1,*\right)$,'' at least,
we have a multiplication rule (5.2.5) ``in the $\mathit{C^{*}}$-algebra
$\mathfrak{H}_{2}^{t}$.'' \hfill{}\textifsymbol[ifgeo]{64}

\bigskip{}

Assume that $\left[\left(a,b\right)\right]_{t}\in\mathcal{H}_{2}^{t\times}$
in $\mathfrak{H}_{2}^{t}$. Then
\[
tr\left(\left[\left(a,b\right)\right]_{t}\right)=a+\overline{a}=2Re\left(a\right),
\]
and\hfill{}(5.2.6)
\[
\tau\left(\left[\left(a,b\right)\right]_{t}\right)=\frac{1}{2}tr\left(\left[\left(a,b\right)\right]_{t}\right)=Re\left(a\right),
\]
where $Re\left(a\right)$ is the real part of $a$ in $\mathbb{C}$.
Similarly, if $\left[\left(a,b\right)\right]_{t}^{*}\in\mathcal{H}_{2}^{t\times}(*)$
in $\mathfrak{H}_{2}^{t}$, then we have
\[
tr\left(\left[\left(a,b\right)\right]_{t}^{*}\right)=tr\left(\begin{array}{cc}
\overline{a} & b\\
t\overline{b} & a
\end{array}\right)=\overline{a}+a=2Re\left(a\right),
\]
and\hfill{}(5.2.7)
\[
\tau\left(\left[\left(a,b\right)\right]_{t}^{*}\right)=\frac{1}{2}\left(2Re\left(a\right)\right)=Re\left(a\right).
\]
Remark that, since $tr$ and $\tau$ are well-defined linear functional
on the $C^{*}$-algebra $\mathfrak{H}_{2}^{t}$, they satisfy
\[
tr\left(T^{*}\right)=\overline{tr\left(T\right)},\;\mathrm{and\;}\tau\left(T^{*}\right)=\overline{\tau\left(T\right)},
\]
for all $T\in\mathfrak{H}_{2}^{t}$. So, the relation (5.2.7) is well-verified,
too.

Also, if $\left[\left(a_{1},b_{1}\right)\right]_{t},\;\left[\left(a_{2},b_{2}\right)\right]_{t}^{*}\in\mathcal{H}_{2}^{t\times}(1,*)$
in $\mathfrak{H}_{2}^{t}$, then

\medskip{}

$\;\;\;\;$$tr\left(\left[\left(a_{1},b_{1}\right)\right]_{t}\left[\left(a_{2},b_{2}\right)\right]_{t}^{*}\right)=tr\left(\left(\begin{array}{ccc}
a_{1}\overline{a_{2}}+t^{2}b_{1}\overline{b_{2}} &  & a_{1}b_{2}+ta_{2}b_{1}\\
\\
\overline{a_{2}b_{1}}+t\overline{a_{1}b_{2}} &  & \overline{b_{1}}b_{2}+\overline{a_{1}}a_{2}
\end{array}\right)\right)$

\medskip{}

\noindent by (5.2.5)

\medskip{}

$\;\;\;\;\;\;\;\;\;\;\;\;\;\;\;\;\;\;\;\;$$=a_{1}\overline{a_{2}}+t^{2}b_{1}\overline{b_{2}}+\overline{b_{1}}b_{2}+\overline{a_{1}}a_{2}$

\medskip{}

$\;\;\;\;\;\;\;\;\;\;\;\;\;\;\;\;\;\;\;\;$$=2Re\left(a_{1}\overline{a_{2}}\right)+t^{2}b_{1}\overline{b_{2}}+\overline{b_{1}}b_{2}$,

\noindent and similarly,\hfill{}(5.2.8)
\[
tr\left(\left[\left(a_{1},b_{1}\right)\right]_{t}^{*}\left[\left(a_{2},b_{2}\right)\right]_{t}\right)=2Re\left(\overline{a_{1}}a_{2}\right)+t^{2}\overline{b_{1}}b_{2}+b_{1}\overline{b_{2}},
\]
and hence,
\[
\tau\left(\left[\left(a_{1},b_{1}\right)\right]_{t}\left[\left(a_{2},b_{2}\right)\right]_{t}^{*}\right)=Re\left(a_{1}\overline{a_{2}}\right)+\frac{t^{2}b_{1}\overline{b_{2}}+\overline{b_{1}}b_{2}}{2},
\]
and\hfill{}(5.2.9)
\[
\tau\left(\left[\left(a_{1},b_{1}\right)\right]_{t}^{*}\left[\left(a_{2},b_{2}\right)\right]_{t}\right)=Re\left(\overline{a_{1}}a_{2}\right)+\frac{t^{2}\overline{b_{1}}b_{2}+b_{1}\overline{b_{2}}}{2},
\]
by (5.2.8).
\begin{prop}
Let $\left(a,b\right),\left(a_{l},b_{l}\right)\in\mathbb{H}_{t}$,
for $l=1,2$, and let $A=\left[\left(a,b\right)\right]_{t}$ and $A_{l}=\left[\left(a_{l},b_{l}\right)\right]_{t}$
be the corresponding realizations of $\mathcal{H}_{2}^{t}$, regarded
as elements of the $t$-scaled-monoidal $C^{*}$-algebra $\mathfrak{H}_{2}^{t}$.
Then
\[
\tau\left(A\right)=\frac{1}{2}tr\left(A\right)=Re\left(a\right)=\frac{1}{2}tr\left(A^{*}\right)=\tau\left(A^{*}\right),
\]
and\hfill{}(5.2.10)
\[
\tau\left(A_{1}A_{2}^{*}\right)=\frac{1}{2}tr\left(A_{1}A_{2}^{*}\right)=Re\left(a_{1}\overline{a_{2}}\right)+\frac{t^{2}b_{1}\overline{b_{2}}+\overline{b_{1}}b_{2}}{2},
\]
and
\[
\tau\left(A_{1}^{*}A_{2}\right)=\frac{1}{2}tr\left(A_{1}^{*}A_{2}\right)=Re\left(\overline{a_{1}}a_{2}\right)+\frac{t^{2}\overline{b_{1}}b_{2}+b_{1}\overline{b_{2}}}{2}.
\]
\end{prop}

\begin{proof}
The joint free moments in (5.2.10) are proven by (5.2.6), (5.2.7),
(5.2.8) and (5.2.9).
\end{proof}
The above computations in (5.2.10) provide a general way to compute
free-distributional data, in particular, the joint free moments of
matrices in the $t$-scaled-monoidal $C^{*}$-algebra $\mathfrak{H}_{2}^{t}$,
up to the trace $tr$, and up to the normalized trace $\tau$. And,
they demonstrate that computing such free-distributional data is not
easy. So, we will restrict our interests to a certain specific case.

\subsection{Free Probability on $\left(\mathfrak{H}_{2}^{t},tr\right)$}

In this section, we fix a scale $t\in\mathbb{R}$, and the corresponding
$t$-scaled-monoidal $C^{*}$-algebra $\mathfrak{H}_{2}^{t}$ generated
by the $t$-scaled hypercomplex monoid $\mathbb{H}_{t}^{\times}$.
Let $\left(\mathfrak{H}_{2}^{t},tr\right)$ be the $C^{*}$-probability
space with respect to the usual trace $tr$ on $\mathfrak{H}_{2}^{t}$. 

Recall that if a scale $t$ is negative, then the realization $\left[\xi\right]_{t}$
and the $t$-spectral form $\Sigma_{t}\left(\xi\right)$ are similar
``in $\mathcal{H}_{2}^{t}$'' by (3.3.9), for all $\xi\in\mathbb{H}_{t}$.
It implies that the similarity ``on $\mathcal{H}_{2}^{t}$'' is
equivalent to the $t$-spectral relation on $\mathbb{H}_{t}$ by (3.3.10).
Also, recall that if two matrices $A$ and $B$ are similar in $M_{n}\left(\mathbb{C}\right)$,
for any $n\in\mathbb{N}$,
\[
tr\left(A\right)=tr\left(B\right).
\]
So, if the realization $\left[\xi\right]_{t}$ and the $t$-spectral
form $\Sigma_{t}\left(\xi\right)$ are similar in $\mathcal{H}_{2}^{t}$,
then the free-moment computations would be much simpler than the computations
of (5.2.10). Note again that if $\left(a,b\right)\in\mathbb{H}_{t}$
satisfies the condition (3.1.5), then
\[
tr\left(\left[\left(a,b\right)\right]_{t}\right)=2Re\left(a\right)=2x=\left(x+i\sqrt{R}\right)+\left(x-i\sqrt{R}\right)=tr\left(\Sigma_{t}\left(a,b\right)\right),
\]
where\hfill{}(5.3.1)
\[
R=y^{2}-tu^{2}-tv^{2}\;\mathrm{in\;\mathbb{R},}
\]
under $\mathbf{RA\:3.2.1}$. Even though the identical results hold
in (5.3.1) (without similarity), if $\left[\left(a,b\right)\right]_{t}$
and $\Sigma_{t}\left(a,b\right)$ are not similar in $\mathcal{H}_{2}^{t}$,
then
\[
tr\left(\left[\left(a,b\right)\right]_{t}^{n}\right)\neq tr\left(\left(\Sigma_{t}\left(a,b\right)\right)^{n}\right),
\]
for some $n\in\mathbb{N}$, by (5.2.5). It implies that some (joint)
free-moments of $\left[\left(a,b\right)\right]_{t}$ and those of
$\Sigma_{t}\left(a,b\right)$ are not identical, and hence, the free
distributions of them are distinct.
\begin{lem}
Suppose the realization $\left[\left(a,b\right)\right]_{t}$ and the
$t$-spectral form $\Sigma_{t}\left(a,b\right)$ are similar in $\mathcal{H}_{2}^{t}$
for $\left(a,b\right)\in\mathbb{H}_{t}$. Then

\medskip{}

\hfill{}$tr\left(\left[\left(a,b\right)\right]_{t}^{n}\right)=2Re\left(\sigma_{t}\left(a,b\right)^{n}\right)=tr\left(\left(\left[\left(a,b\right)\right]_{t}^{*}\right)^{n}\right)$\hfill{}(5.3.2)

\medskip{}

\noindent for all $n\in\mathbb{N}$, where $\sigma_{t}\left(a,b\right)$
is the $t$-spectral value of $\left(a,b\right)$.
\end{lem}

\begin{proof}
Suppose $\left(a,b\right)\in\mathbb{H}_{t}$ satisfies the condition
(3.1.5). Then
\[
\left[\left(a,b\right)\right]_{t}=\left(\begin{array}{cc}
a & tb\\
\overline{b} & \overline{a}
\end{array}\right)\;\mathrm{and\;}\Sigma_{t}\left(\left(a,b\right)\right)=\left(\begin{array}{cc}
\sigma_{t}\left(a,b\right) & 0\\
0 & \overline{\sigma_{t}\left(a,b\right)}
\end{array}\right),
\]
in $\mathcal{H}_{2}^{t}$, where
\[
\sigma_{t}\left(a,b\right)=x+i\sqrt{y^{2}-tu^{2}-tv^{2}},
\]
under $\mathbf{RA\:3.2.1}$. Assume that $\left[\left(a,b\right)\right]_{t}$
and $\Sigma_{t}\left(\left(a,b\right)\right)$ are similar in $\mathcal{H}_{2}^{t}$.
Then the matrices $\left[\left(a,b\right)\right]_{t}^{n}$ and $\Sigma_{t}\left(\left(a,b\right)\right)^{n}$
are similar in $\mathcal{H}_{2}^{t}$, for all $n\in\mathbb{N}$.
Indeed, if two elements $A$ and $B$ are similar in $\mathcal{H}_{2}^{t}$,
satisfying $B=U^{-1}AU$ in $\mathcal{H}_{2}^{t}$, for an invertible
element $U\in\mathcal{H}_{2}^{t}$, then
\[
B^{n}=\left(U^{-1}AU\right)^{n}=U^{-1}A^{n}U\;\mathrm{in\;}\mathcal{H}_{2}^{t},
\]
implying the similarity of $A^{n}$ and $B^{n}$, for $n\in\mathbb{N}$.
Thus,
\[
tr\left(\left[\left(a,b\right)\right]_{t}^{n}\right)=tr\left(\Sigma_{t}\left(\left(a,b\right)\right)^{n}\right),
\]
and
\[
tr\left(\Sigma_{t}\left(\left(a,b\right)\right)^{n}\right)=tr\left(\left(\begin{array}{cc}
\sigma_{t}\left(a,b\right)^{n} & 0\\
0 & \overline{\sigma_{t}\left(a,b\right)^{n}}
\end{array}\right)\right),
\]
implying that
\[
tr\left(\left[\left(a,b\right)\right]_{t}^{n}\right)=tr\left(\Sigma_{t}\left(\left(a,b\right)\right)^{n}\right)=2Re\left(\sigma_{t}\left(a,b\right)^{n}\right),
\]
for all $n\in\mathbb{N}$. Therefore, the first equality in (5.3.2)
holds.

Since $tr$ is a well-defined linear functional on the $C^{*}$-algebra
$\mathfrak{H}_{2}^{t}$, one has
\[
tr\left(A^{*}\right)=\overline{tr\left(A\right)},\;\mathrm{for\;all\;}A\in\mathfrak{H}_{2}^{t}.
\]
Since
\[
tr\left(\left(\left[\left(a,b\right)\right]_{t}^{*}\right)^{n}\right)=tr\left(\left(\left[\left(a,b\right)\right]_{t}^{n}\right)^{*}\right)=\overline{tr\left(\left[\left(a,b\right)\right]_{t}^{n}\right)},
\]
one has
\[
tr\left(\left(\left[\left(a,b\right)\right]_{t}^{*}\right)^{n}\right)=\overline{2Re\left(\sigma_{t}\left(a,b\right)^{n}\right)}=2Re\left(\sigma_{t}\left(a,b\right)^{n}\right),
\]
for all $n\in\mathbb{N}$. So, the second equality in (5.3.2) holds,
too.
\end{proof}
Note that the formula (5.3.2) holds true under the similarity assumption
of the realization and the $t$-spectral form. 

Remark that every complex number $w\in\mathbb{C}$ is polar-decomposed
to be
\[
w=\left|w\right|w_{o}\;\mathrm{with}\;w_{o}\in\mathbb{T},
\]
uniquely, where $\mathbb{T}=\left\{ z\in\mathbb{C}:\left|z\right|=1\right\} $
is the unit circle in $\mathbb{C}$. So, all our $t$-spectral values
$\sigma_{t}\left(\xi\right)$ are polar-decomposed to be
\[
\sigma_{t}\left(\xi\right)=\left|\sigma_{t}\left(\xi\right)\right|\sigma_{t}\left(\xi\right)_{o}\;\mathrm{with\;}\sigma_{t}\left(\xi\right)_{o}\in\mathbb{T},
\]
for all $\xi\in\mathbb{H}_{t}$. In such a sense, we have that
\[
tr\left(\left[\xi\right]_{t}^{n}\right)=2\left|\sigma_{t}\left(\xi\right)\right|^{n}Re\left(\sigma_{t}\left(\xi\right)_{o}^{n}\right),
\]
for all $n\in\mathbb{N}$, by (5.3.2).
\begin{cor}
Suppose the realization $\left[\xi\right]_{t}$ and the $t$-spectral
form $\Sigma_{t}\left(\xi\right)$ are similar in $\mathcal{H}_{2}^{t}$
for $\xi\in\mathbb{H}_{t}$. Then

\medskip{}

\hfill{}$tr\left(\left[\xi\right]_{t}^{n}\right)=2\left|\sigma_{t}\left(\xi\right)\right|^{n}Re\left(\sigma_{t}\left(\xi\right)_{o}^{n}\right)=tr\left(\left(\left[\xi\right]_{t}^{*}\right)^{n}\right),$\hfill{}(5.3.3)

\medskip{}

\noindent for all $n\in\mathbb{N}$, where $\sigma_{t}\left(\xi\right)=\left|\sigma_{t}\left(\xi\right)\right|\sigma_{t}\left(\xi\right)_{o}$
is the polar decomposition of $\sigma_{t}\left(\xi\right)$, with
$\sigma_{t}\left(\xi\right)_{o}\in\mathbb{T}$.
\end{cor}

\begin{proof}
The free-distributional data (5.3.3) is immediately obtained by (5.3.2)
under the polar decomposition of the $t$-spectral value $\sigma_{t}\left(\xi\right)$
in $\mathbb{C}$.
\end{proof}
Assume again that a hypercomplex number $\left(a,b\right)\in\mathbb{H}_{t}$
satisfies our similarity assumption, i.e., $T\overset{\textrm{denote}}{=}\left[\left(a,b\right)\right]_{t}$
and $S\overset{\textrm{denote}}{=}\Sigma_{t}\left(\left(a,b\right)\right)$
are similar in $\mathcal{H}_{2}^{t}$. Then, for any 
\[
\left(r_{1},...,r_{n}\right)\in\left\{ 1,*\right\} ^{n},\;\mathrm{for\;}n\in\mathbb{N},
\]
the matrix $\overset{n}{\underset{l=1}{\prod}}T^{r_{l}}$ is similar
to $\overset{n}{\underset{l=1}{\prod}}S^{r_{l}}$ in $\mathcal{H}_{2}^{t}$
(and hence, in $\mathfrak{H}_{2}^{t}$).
\begin{thm}
Let $\left(a,b\right)\in\mathbb{H}_{t}$ satisfy the similarity assumption
that: $T\overset{\textrm{denote}}{=}\left[\left(a,b\right)\right]_{t}$
and $S\overset{\textrm{denote}}{=}\Sigma_{t}\left(\left(a,b\right)\right)$
are similar in $\mathcal{H}_{2}^{t}$. If
\[
\sigma_{t}\left(a,b\right)=rw_{o},\;\textrm{polar decomposition,}
\]
with\hfill{}(5.3.4)
\[
r=\left|\sigma_{t}\left(a,b\right)\right|\;\mathrm{and\;}w_{o}\in\mathbb{T},
\]
then

\medskip{}

\hfill{}$tr\left(\overset{n}{\underset{l=1}{\prod}}T^{r_{l}}\right)=2r^{n}Re\left(w_{o}^{\overset{n}{\underset{l=1}{\sum}}e_{l}}\right),$\hfill{}(5.3.5)

\medskip{}

\noindent for all $\left(r_{1},...,r_{n}\right)\in\left\{ 1,*\right\} ^{n}$,
for all $n\in\mathbb{N}$, where
\[
e_{l}=\left\{ \begin{array}{ccc}
1 &  & \mathrm{if\;}r_{l}=1\\
-1 &  & \mathrm{if\;}r_{l}=*,
\end{array}\right.
\]
for all $l=1,...,n$.
\end{thm}

\begin{proof}
Since the realization $T$ and the $t$-spectral form $S$ are assumed
to be similar in $\mathcal{H}_{2}^{t}$, their adjoints $T^{*}$ and
$S^{*}$ are similar in $\mathcal{H}_{2}^{t\times}(*)\cup\left\{ \left[\left(0,0\right)\right]_{t}\right\} $;
and hence, the matrix $\overset{n}{\underset{l=1}{\prod}}T^{r_{l}}$
and $\overset{n}{\underset{l=1}{\prod}}S^{r_{l}}$ are similar ``in
$\mathfrak{H}_{2}^{t}$.'' Consider that
\[
S=\left(\begin{array}{ccc}
\sigma_{t}\left(a,b\right) &  & 0\\
\\
0 &  & \overline{\sigma_{t}\left(a,b\right)}
\end{array}\right)=\left(\begin{array}{cc}
rw_{o} & 0\\
0 & r\overline{w_{o}}
\end{array}\right)=r\left(\begin{array}{cc}
w_{o} & 0\\
0 & w_{o}^{-1}
\end{array}\right),
\]
under hypotheses, because $\overline{z}=\frac{1}{z}=z^{-1}$ in $\mathbb{T}$,
whenever $z\in\mathbb{T}$ in $\mathbb{C}$. It shows that
\[
S^{j}=r^{j}\left(\begin{array}{cc}
w_{o}^{j} & 0\\
0 & w_{o}^{-j}
\end{array}\right),\;\mathrm{for\;all\;}j\in\mathbb{N}\cup\left\{ 0\right\} ,
\]
and
\[
S^{*}=\overline{r}\left(\begin{array}{cc}
\overline{w_{o}} & 0\\
0 & w_{o}
\end{array}\right)=r\left(\begin{array}{cc}
w_{o}^{-1} & 0\\
0 & w_{o}
\end{array}\right),
\]
satisfying that
\[
\left(S^{*}\right)^{j}=\left(S^{j}\right)^{*},\;\mathrm{for\;all\;}j\in\mathbb{N}.
\]
It implies that, for any $\left(r_{1},...,r_{n}\right)\in\left\{ 1,*\right\} ^{n}$,
for $n\in\mathbb{N}$, there exists $\left(e_{1},...,e_{n}\right)\in\left\{ \pm1\right\} ^{n}$,
such that
\[
e_{l}=\left\{ \begin{array}{ccc}
1 &  & \mathrm{if\;}r_{l}=1\\
-1 &  & \mathrm{if\;}r_{l}=*,
\end{array}\right.
\]
for all $l=1,...,n$, and\hfill{}(5.3.6)
\[
\overset{n}{\underset{l=1}{\prod}}S^{r_{l}}=r^{n}\left(\begin{array}{ccc}
w_{o}^{\overset{n}{\underset{l=1}{\sum}}e_{l}} &  & 0\\
\\
0 &  & w_{o}^{-\left(\overset{n}{\underset{l=1}{\sum}}e_{l}\right)}
\end{array}\right),
\]
in $\mathfrak{H}_{2}^{t}$. Thus, under our similarity assumption,
\[
tr\left(\overset{n}{\underset{l=1}{\prod}}T^{r_{l}}\right)=tr\left(\overset{n}{\underset{l=1}{\prod}}S^{r_{l}}\right)=r^{n}\left(w_{o}^{\overset{n}{\underset{l=1}{\sum}}e_{l}}+w_{o}^{-\left(\overset{n}{\underset{l=1}{\sum}}e_{l}\right)}\right),
\]
implying that
\[
tr\left(\overset{n}{\underset{l=1}{\prod}}T^{r_{l}}\right)=r^{n}\left(2Re\left(w_{o}^{\overset{n}{\underset{l=1}{\sum}}e_{l}}\right)\right),
\]
for all $\left(r_{1},...,r_{n}\right)\in\left\{ 1,*\right\} ^{n}$,
for all $n\in\mathbb{N}$, where $\left(e_{1},...,e_{n}\right)\in\left\{ \pm1\right\} ^{n}$
satisfies (5.3.6).

Therefore, under our similarity assumption and the polar decomposition
(5.3.4), the free-distributional data (5.3.5) holds.
\end{proof}
By the above theorem, one immediately obtain the following result.
\begin{cor}
Let $\left(a,b\right)\in\mathbb{H}_{t}$ satisfy the similarity assumption
that: $T\overset{\textrm{denote}}{=}\left[\left(a,b\right)\right]_{t}$
and $S\overset{\textrm{denote}}{=}\Sigma_{t}\left(\left(a,b\right)\right)$
are similar in $\mathcal{H}_{2}^{t}$. If
\[
\sigma_{t}\left(a,b\right)=rw_{o},\;\textrm{polar decomposition,}
\]
with\hfill{}(5.3.7)
\[
r=\left|\sigma_{t}\left(a,b\right)\right|\;\mathrm{and\;}w_{o}\in\mathbb{T},
\]
then

\medskip{}

\hfill{}$\tau\left(\overset{n}{\underset{l=1}{\prod}}T^{r_{l}}\right)=r^{n}Re\left(w_{o}^{\overset{n}{\underset{l=1}{\sum}}e_{l}}\right),$\hfill{}(5.3.8)

\medskip{}

\noindent for all $\left(r_{1},...,r_{n}\right)\in\left\{ 1,*\right\} ^{n}$,
for all $n\in\mathbb{N}$, where
\[
e_{l}=\left\{ \begin{array}{ccc}
1 &  & \mathrm{if\;}r_{l}=1\\
-1 &  & \mathrm{if\;}r_{l}=*,
\end{array}\right.
\]
for all $l=1,...,n$.
\end{cor}

\begin{proof}
By (5.3.5), the free-distributional data (5.3.8) holds up to the normalized
trace $\tau=\frac{1}{2}tr$ on $\mathfrak{H}_{2}^{t}$, under (5.3.7).
\end{proof}
Under our similarity assumption and the condition (5.3.7), the free-distributional
data (5.3.8) fully characterizes the free distribution of $\left[\left(a,b\right)\right]_{t}\in\mathcal{H}_{2}^{t}$
in the $C^{*}$-probability space $\left(\mathfrak{H}_{2}^{t},\:\tau\right)$.
\begin{cor}
Suppose a given scale $t$ is negative in $\mathbb{R}$. Let $\left(a,b\right)\in\mathbb{H}_{t}$,
and let $T\overset{\textrm{denote}}{=}\left[\left(a,b\right)\right]_{t}$
and $S\overset{\textrm{denote}}{=}\Sigma_{t}\left(\left(a,b\right)\right)$
in $\mathcal{H}_{2}^{t}$. If
\[
\sigma_{t}\left(a,b\right)=rw_{o},\;\textrm{polar decomposition,}
\]
with\hfill{}(5.3.9)
\[
r=\left|\sigma_{t}\left(a,b\right)\right|\;\mathrm{and\;}w_{o}\in\mathbb{T},
\]
then

\medskip{}

\hfill{}$tr\left(\overset{n}{\underset{l=1}{\prod}}T^{r_{l}}\right)=2r^{n}Re\left(w_{o}^{\overset{n}{\underset{l=1}{\sum}}e_{l}}\right)=2\tau\left(\overset{n}{\underset{l=1}{\prod}}T^{r_{l}}\right),$\hfill{}(5.3.10)

\medskip{}

\noindent for all $\left(r_{1},...,r_{n}\right)\in\left\{ 1,*\right\} ^{n}$,
for all $n\in\mathbb{N}$, where
\[
e_{l}=\left\{ \begin{array}{ccc}
1 &  & \mathrm{if\;}r_{l}=1\\
-1 &  & \mathrm{if\;}r_{l}=*,
\end{array}\right.
\]
for all $l=1,...,n$.
\end{cor}

\begin{proof}
In Theorem 48 and Corollary 49, we showed that if $T$ and $S$ are
similar in $\mathcal{H}_{2}^{t}$, then the free-distributional data
(5.3.10) holds under the condition (5.3.9), by (5.3.5) and (5.3.8),
respectively. So, it suffices to show that the realization $T$ and
the $t$-spectral form $S$ are similar in $\mathcal{H}_{2}^{t}$.
However, since $t<0$ in $\mathbb{R}$, the matrices $T$ and $S$
are similar in $\mathcal{H}_{2}^{t}$ by (3.3.9).
\end{proof}
The above corollary shows that, if a given scale $t$ is negative
in $\mathbb{R}$, then the free-distributional data (5.3.10) fully
characterizes the free distributions of the realizations $\left[\xi\right]_{t}$
in the $t$-scaled-monoidal $C^{*}$-algebra $\mathfrak{H}_{2}^{t}$
up to the usual trace $tr$, and the normalized trace $\tau$, for
``all'' $\xi\in\mathbb{H}_{t}$. In other words, it illustrates
that, if $t<0$ in $\mathbb{R}$, then the free-distributional data
on the $C^{*}$-probability spaces,
\[
\left(\mathfrak{H}_{2}^{t},\:tr\right)\;\mathrm{and\;}\left(\mathfrak{H}_{2}^{t},\:\tau\right),
\]
are fully characterized by the spectra of hypercomplex numbers of
$\mathbb{H}_{t}$, by (5.3.9) and (5.3.10).

But, if $t\geq0$, and hence, there are some hypercomplex numbers
$\eta$ of $\mathbb{H}_{t}$ whose realization and spectral form are
not similar in $\mathcal{H}_{2}^{t}$, then computing joint free moments
of $\left[\eta\right]_{t}$ in $\mathfrak{H}_{2}^{t}$ would not be
easy e.g., see (5.2.10).

\subsection{More Free-Distributional Data on $\left(\mathfrak{H}_{2}^{t},\:\tau\right)$
for $t<0$}

In this section, a fixed scale $t$ is automatically assumed to be
negative, i.e., $t<0$ in $\mathbb{R}$. At this moment, we emphasize
that most main results of this section would hold even though $t$
is not negative in $\mathbb{R}$. However, we assume a given scale
$t$ is negative for convenience (e.g., see (5.3.10)). Let $\mathfrak{H}_{2}^{t}$
be the $t$-scaled-monoidal $C^{*}$-algebra inducing a $C^{*}$-probability
space $\left(\mathfrak{H}_{2}^{t},\tau\right)$, where $\tau$ is
the normalized trace on $\mathfrak{H}_{2}^{t}$. Since $t$ is assumed
to be negative in $\mathbb{R}$, the realizations $T=\left[\eta\right]_{t}$
and the $t$-spectral forms $S=\Sigma_{t}\left(\eta\right)$ are similar
in $\mathcal{H}_{2}^{t}$ by (3.3.9), and hence,
\[
\tau\left(\overset{n}{\underset{l=1}{\prod}}T^{r_{l}}\right)=r^{n}Re\left(w_{o}^{\overset{n}{\underset{l=1}{\sum}}e_{l}}\right)=\tau\left(\overset{n}{\underset{l=1}{\prod}}S^{r_{l}}\right),
\]
by (5.3.5), where\hfill{}(5.4.1)
\[
\sigma_{t}\left(\eta\right)=rw_{o}\in\mathbb{C},\;\textrm{polar decomposition,}
\]
with $r=\left|\sigma_{t}\left(\eta\right)\right|$ and $w_{o}\in\mathbb{T}$,
for all $\left(r_{1},...,r_{n}\right)\in\left\{ 1,*\right\} ^{n}$,
where $\left(e_{1},...,e_{n}\right)\in\left\{ \pm1\right\} ^{n}$
satisfies (5.3.6), for all $n\in\mathbb{N}$, for ``all'' $\eta\in\mathbb{H}_{t}$.
And the free-distributional data (5.4.1) fully characterizes the free
distribution of $\left[\eta\right]_{t}\in\left(\mathfrak{H}_{2}^{t},\tau\right)$,
for all $\eta\in\mathbb{H}_{t}$.

In this section, we refine (5.4.1) case-by-case, up to operator-theoretic
properties of elements of $\left(\mathfrak{H}_{2}^{t},\tau\right)$. 
\begin{defn}
Let $\mathcal{A}$ be a unital $C^{*}$-algebra with its unity $1_{\mathcal{A}}$,
and let $T\in\mathcal{A}$, and $T^{*}\in\mathcal{A}$, the adjoint
of $T$. 

\noindent (1) $T$ is said to be self-adjoint, if $T^{*}=T$ in $\mathcal{A}$.

\noindent (2) $T$ is a projection, if $T^{*}=T=T^{2}$ in $\mathcal{A}$.

\noindent (3) $T$ is normal, if $T^{*}T=TT^{*}$ in $\mathcal{A}$.

\noindent (4) $T$ is a unitary, if $T^{*}T=1_{\mathcal{A}}=TT^{*}$
in $\mathcal{A}$.
\end{defn}

Let $\left(a,b\right)\in\mathbb{H}_{t}$, satisfying the condition
(3.1.5), and $T\overset{\textrm{denote}}{=}\left[\left(a,b\right)\right]_{t}\in\mathcal{H}_{2}^{t}$,
as an element of $\left(\mathfrak{H}_{2}^{t},\tau\right)$. Then its
adjoint,
\[
T^{*}=\left(\begin{array}{cc}
\overline{a} & b\\
t\overline{b} & a
\end{array}\right)\in\mathcal{H}_{2}^{t}(*),
\]
is well-defined in $\left(\mathfrak{H}_{2}^{t},\tau\right)$, and
the corresponding $t$-spectral form,
\[
S\overset{\textrm{denote}}{=}\Sigma_{t}\left(\left(a,b\right)\right)=\left(\begin{array}{cc}
w & 0\\
0 & \overline{w}
\end{array}\right)\in\mathcal{H}_{2}^{t},
\]
is contained in $\left(\mathfrak{H}_{2}^{t},\tau\right)$, where $\overline{w}$
is determined by $\mathbf{RA\;3.2.1}$, and 
\[
w=\sigma_{t}\left(a,b\right)=x+i\sqrt{y^{2}-tu^{2}-tv^{2}}
\]
is the $t$-spectral value, uniqely polar-decomposed to be 
\[
w=rw_{o}\;\mathrm{with\;}r=\left|\sigma_{t}\left(a,b\right)\right|\;\mathrm{and\;}w_{o}\in\mathbb{T}.
\]

\medskip{}

\noindent $\mathbf{Assumption\;and\;Notation\;5.4.1.}$ (from below
$\mathbf{AN\;5.4.1}$) From now on, if we say that ``a given hypercomplex
number $\left(a,b\right)\in\mathbb{H}_{t}$ satisfies $\mathbf{AN\;5.4.1}$,''
then it means it has its realization denoted by $T$, its $t$-spectral
form denoted by $S$, determined by the $t$-spectral value denoted
by $w$, which is polar-decomposed to be $w=rw_{o}$, as indicated
in the very above paragraph. \hfill{}\textifsymbol[ifgeo]{64}

\medskip{}

Let $\left(a,b\right)\in\mathbb{H}_{t}$ satisfy $\mathbf{AN\;5.4.1}$.
Then the self-adjointness of the realization $T\in\mathcal{H}_{2}^{t}$
in $\mathfrak{H}_{2}^{t}$ says that 
\[
T^{*}=T\Longleftrightarrow\left(\begin{array}{cc}
\overline{a} & b\\
t\overline{b} & a
\end{array}\right)=\left(\begin{array}{cc}
a & tb\\
\overline{b} & \overline{a}
\end{array}\right),
\]
if and only if
\[
\overline{a}=a\;\mathrm{and\;}tb=b\;\mathrm{in\;}\mathbb{C},
\]
if and only if\hfill{}(5.4.2)
\[
a\in\mathbb{R}\;\mathrm{and\;}b=0.
\]
Especially, the equality $b=0$ in (5.4.2) is obtained by our negative-scale
assumption: $t<0$ in $\mathbb{R}$.
\begin{prop}
Let $\left(a,b\right)\in\mathbb{H}_{t}$ satisfy $\mathbf{AN\;5.4.1}$.
Then the realization $T\in\mathcal{H}_{2}^{t}$ is self-adjoint in
$\mathfrak{H}_{2}^{t}$, if and only if

\medskip{}

\hfill{}$a\in\mathbb{R}\;\mathrm{and\;}b=0\Longleftrightarrow\left(a,b\right)=\left(Re\left(a\right),0\right)\;\mathrm{in\;}\mathbb{H}_{t}$.\hfill{}(5.4.3)
\end{prop}

\begin{proof}
The self-adjointness (5.4.3) is shown by (5.4.2).
\end{proof}
The self-adjointness (5.4.3) illustrates that the self-adjoint generating
elements $T\in\mathcal{H}_{2}^{t}$ of $\left(\mathfrak{H}_{2}^{t},\tau\right)$
have their forms,
\[
T=\left(\begin{array}{cc}
x & 0\\
0 & x
\end{array}\right)\in\mathcal{H}_{2}^{t}\left(1,*\right)\;\mathrm{with\;}x\in\mathbb{R}.
\]

\bigskip{}

\noindent $\mathbf{Remark\;and\;Observation.}$ The above self-adjointness
characterization (5.4.3) is obtained for the case where $t<0$ in
$\mathbb{R}$. How about the other cases? Generally, one has $T$
is self-adjont in $\mathfrak{H}_{2}^{t}$, if and only if 
\[
\overline{a}=a\;\mathrm{and\;}tb=b,
\]
like (5.4.2). Thus one can verify that: (i) if $t=0$, then $T$ is
self-adjoint, if and only if $a\in\mathbb{R}$ and $b=0$, just like
(5.4.3); (ii) if $t>0$ and $t\neq1$, then $T$ is self-adjoint,
if and only if $a\in\mathbb{R}$ and $b=0$, just like (5.4.3); meanwhile,
(iii) if $t=1$ (equivalently, if $\left(a,b\right)$ is a bicomplex
number of $\mathbb{H}_{1}$), then $T$ is self-adjoint in $\mathfrak{H}_{2}^{1}$,
if and only if $a\in\mathbb{R}$, if and only if $\left(a,b\right)=\left(Re\left(a\right),b\right)$
in $\mathbb{H}_{1}$. In summary,
\[
T\textrm{ is self-adjoint in }\mathfrak{H}_{2}^{t}\Longleftrightarrow\left(a,b\right)=\left(Re\left(a\right),0\right)\textrm{ in }\mathbb{H}_{t},
\]
like (5.4.3), whenever $t\in\mathbb{R}\setminus\left\{ 1\right\} $,
meanwhile,
\[
T\textrm{ is self-adjoint in }\mathfrak{H}_{2}^{1}\Longleftrightarrow\left(a,b\right)=\left(Re\left(a\right),b\right)\in\mathbb{H}_{1}.
\]
\textifsymbol[ifgeo]{64}

\bigskip{}

Now, let $\left(a,b\right)\in\mathbb{H}_{t}$, under $\mathbf{AN\;5.4.1}$
and our negative-scale assumption, satisfy the self-adjointness (5.4.3),
i.e., it is actually $\left(a,0\right)$ with $a\in\mathbb{R}$. Then
\[
T=\left(\begin{array}{cc}
a & 0\\
0 & a
\end{array}\right)=S\;\mathrm{in\;}\mathcal{H}_{2}^{t}\left(1,*\right),
\]
as an element of $\mathfrak{H}_{2}^{t}$.
\begin{thm}
Let $\left(a,b\right)\in\mathbb{H}_{t}$ satisfy $\mathbf{AN\;5.4.1}$,
and assume that the realization $T$ is self-adjoint in $\left(\mathfrak{H}_{2}^{t},\tau\right)$.
Then

\medskip{}

\hfill{}$\tau\left(\overset{n}{\underset{l=1}{\prod}}T^{r_{l}}\right)=\tau\left(T^{n}\right)=a^{n}$
$\;$in$\;$ $\mathbb{R}$,\hfill{}(5.4.4)

\medskip{}

\noindent for all $\left(r_{1},...,r_{n}\right)\in\left\{ 1,*\right\} ^{n}$,
for all $n\in\mathbb{N}$.
\end{thm}

\begin{proof}
By the self-adjointness (5.4.3) of the realization $T$ of $\left(a,b\right)\in\mathbb{H}_{t}$,
one has $\left(a,b\right)=\left(a,0\right)$ in $\mathbb{H}_{t}$,
with $a\in\mathbb{R}$, and
\[
T=S=\left(\begin{array}{cc}
a & 0\\
0 & a
\end{array}\right)=S^{*}=T^{*}\;\mathrm{in\;}\mathfrak{H}_{2}^{t}.
\]
So,
\[
\tau\left(\overset{n}{\underset{l=1}{\prod}}T^{r_{l}}\right)=\tau\left(T^{n}\right)=\tau\left(S^{n}\right)=\tau\left(\left(\begin{array}{cc}
a^{n} & 0\\
0 & a^{n}
\end{array}\right)\right),
\]
for all $\left(r_{1},...,r_{n}\right)\in\left\{ 1,*\right\} ^{n}$,
for all $n\in\mathbb{N}$. Therefore, the free-distributional data
(5.4.4) holds true.
\end{proof}
\bigskip{}

\noindent $\mathbf{Observation.}$ Similar to the above theorem, one
can verify that: if $t\in\mathbb{R}\setminus\left\{ 1\right\} $,
then the free-distributional data (5.4.4) holds for self-adjoint realizations
$T\in\left(\mathfrak{H}_{2}^{t},\tau\right)$ of $\left(a,0\right)\in\mathbb{H}_{t}$
with $a\in\mathbb{R}$.\hfill{}\textifsymbol[ifgeo]{64}

\bigskip{}

By (5.4.3), the realization $T$ of a hypercomplex number $\left(a,b\right)\in\mathbb{H}_{t}$,
satisfying $\mathbf{AN\;5.4.1}$, is self-adjoint, if and only if
$\left(a,b\right)=\left(a,0\right)$ with $a\in\mathbb{R}$. And,
by definition, such a self-adjoint matrix $T$ can be a projection,
if and only if it is idempotent in the sense that
\[
T^{2}=T\;\;\mathrm{in\;\;}\mathfrak{H}_{2}^{t}.
\]
Observe that a self-adjoint realization $T$ satisfies the above idempotence,
if and only if
\[
T^{2}=\left(\begin{array}{cc}
a^{2} & 0\\
0 & a^{2}
\end{array}\right)=\left(\begin{array}{cc}
a & 0\\
0 & a
\end{array}\right)=T,
\]
if and only if\hfill{}(5.4.5)
\[
a^{2}=a\Longleftrightarrow a=0,\;\mathrm{or\;}a=1,\;\mathrm{in\;}\mathbb{R}.
\]

\begin{prop}
Let $\left(a,b\right)\in\mathbb{H}_{t}$ satisfy $\mathbf{AN\;5.4.1}$.
Then the realization $T$ is a projection, if and only if

\medskip{}

\hfill{}$\mathrm{either\;}T=I_{2},\;\mathrm{or\;}T=O_{2}\;\mathrm{in\;}\mathcal{H}_{2}^{t},$\hfill{}(5.4.6)

\medskip{}

\noindent where $I_{2}=\left[\left(1,0\right)\right]_{t}$ is the
identity matrix, and $O_{2}=\left[\left(0,0\right)\right]_{t}$ is
the zero matrix of $\mathfrak{H}_{2}^{t}$.
\end{prop}

\begin{proof}
The operator-equality (5.4.6) holds in $\mathcal{H}_{2}^{t}$ (and
hence, in $\mathfrak{H}_{2}^{t}$) by (5.4.5).
\end{proof}
\bigskip{}

\noindent $\mathbf{Observation.}$ Like the above proposition, one
can conclude that: if $t\in\mathbb{R}\setminus\left\{ 1\right\} $,
then the realization $T$ is a projection in $\mathfrak{H}_{2}^{t}$,
if and only if it is either the identity matrix $I_{2}$, or the zero
matrix $O_{2}$ of $\mathfrak{H}_{2}^{t}$. How about the case where
$t=1$? As we discussed above, $T\in\mathfrak{H}_{2}^{1}$ is self-adjoint,
if and only if $\left(a,b\right)=\left(Re\left(a\right),b\right)$
in $\mathbb{H}_{1}$, if and only if
\[
T=\left(\begin{array}{cc}
x & b\\
\overline{b} & x
\end{array}\right)\in\mathcal{H}_{2}^{1},\;\mathrm{and\;}S=\left(\begin{array}{ccc}
x+i\sqrt{-u^{2}-v^{2}} &  & 0\\
\\
0 &  & x-i\sqrt{-u^{2}-v^{2}}
\end{array}\right),
\]
implying that
\[
S=\left(\begin{array}{ccc}
x-\left|b\right| &  & 0\\
\\
0 &  & x+\left|b\right|
\end{array}\right)\;\;\mathrm{in\;\;}\mathfrak{H}_{2}^{1},
\]
under $\mathbf{AN\;5.4.1}$. Such a self-adjoint $T$ is a projection,
if and only if $T^{2}=T$ in $\mathfrak{H}_{2}^{1}$, if and only
if
\[
x^{2}+\left|b\right|^{2}=x\;\mathrm{\;and\;\;}2xb=b.
\]
Thus if $b=0$, then $x\in\left\{ 0,1\right\} $, meanwhile, if $b\neq0$,
then
\[
x^{2}+\left|b\right|^{2}=x\;\;\mathrm{and\;\;}x=\frac{1}{2},
\]
$\Longleftrightarrow$
\[
x=\frac{1}{2}\;\mathrm{\;and\;\;}\frac{1}{4}+\left|b\right|^{2}=\frac{1}{2},
\]
$\Longleftrightarrow$
\[
x=\frac{1}{2}\;\;\;\mathrm{and\;\;\;}\left|b\right|^{2}=\frac{1}{4},
\]
if and only if
\[
\left(a,b\right)=\left(\frac{1}{2},b\right)\;\;\mathrm{with\;\;}\left|b\right|^{2}=\frac{1}{4}.
\]
It implies that $T$ is a projection in $\mathfrak{H}_{2}^{1}$, if
and only if
\[
\left(a,b\right)=\left(0,0\right),\;\mathrm{or\;}\left(a,b\right)=\left(1,0\right),
\]
or
\[
\left(a,b\right)=\left(\frac{1}{2},\:b\right)\;\mathrm{with\;}\left|b\right|^{2}=\frac{1}{4},
\]
in $\mathbb{H}_{1}$.\hfill{}\textifsymbol[ifgeo]{64}

\bigskip{}

The above proposition says that, under our negative-scale assumption,
the only projections of $\mathfrak{H}_{2}^{t}$ induced by hypercomplex
numbers of $\mathbb{H}_{t}$ are the identity element $I_{2}=\left[\left(1,0\right)\right]_{t}$,
and the zero element $O_{2}=\left[\left(0,0\right)\right]_{t}$ in
$\mathfrak{H}_{2}^{t}.$ For any unital $C^{*}$-probability spaces
$\left(\mathcal{A},\varphi\right)$, the unity $1_{\mathcal{A}}$
has its free distributions characterized by its free-moment sequence,
\[
\left(\varphi\left(1_{\mathcal{A}}^{n}\right)=\varphi\left(1_{\mathcal{A}}\right)\right)_{n=1}^{\infty}=\left(1,1,1,1,1,...\right);
\]
and the free distribution of the zero element $0_{\mathcal{A}}$ is
nothing but the zero-free distribution, characterized by the free-moment
sequence,
\[
\left(\varphi\left(0_{\mathcal{A}}^{n}\right)=\varphi\left(0_{\mathcal{A}}\right)\right)_{n=1}^{\infty}=\left(0,0,0,0,...\right).
\]

\begin{thm}
Let $\left(a,b\right)\in\mathbb{H}_{t}$, satisfying $\mathbf{AN\;5.4.1}$,
have its realization $T\in\mathcal{H}_{2}^{t}$, which is a ``non-zero''
projection in $\mathfrak{H}_{2}^{t}$. Then
\[
\tau\left(T^{n}\right)=1,\;\;\forall n\in\mathbb{N}.
\]
(In fact, this result holds true for all $t\in\mathbb{R}\setminus\left\{ 1\right\} $.)
\end{thm}

\begin{proof}
Under hypothesis, the realization $T\in\mathcal{H}_{2}^{t}$ is a
projection in $\mathfrak{H}_{2}^{t}$, if and only if $\left(a,b\right)=\left(1,0\right)$,
or $\left(0,0\right)$ in $\mathbb{H}_{t}$, by (5.4.6). Since $T\in\mathcal{H}_{2}^{t}$
is assumed to a non-zero projection in $\mathfrak{H}_{2}^{t}$, we
have
\[
\left(a,b\right)=\left(1,0\right)\;\mathrm{in\;\mathbb{H}_{t},}\Longleftrightarrow T=I_{2}=S\;\mathrm{in\;}\mathfrak{H}_{2}^{t}.
\]
Therefore,
\[
\tau\left(T^{n}\right)=\tau\left(I_{2}^{n}\right)=1,\;\forall n\in\mathbb{N}.
\]
(Note that it holds true for all $t\in\mathbb{R}\setminus\left\{ 1\right\} $.)
\end{proof}
Let $\left(a,b\right)\in\mathbb{H}_{t}$ satisfy $\mathbf{AN\;5.4.1}$,
and let $T\in\mathcal{H}_{2}^{t}$ be the realization in $\mathfrak{H}_{2}^{t}$.
Observe that
\[
T^{*}T=\left(\begin{array}{cc}
\overline{a} & b\\
t\overline{b} & a
\end{array}\right)\left(\begin{array}{cc}
a & tb\\
\overline{b} & \overline{a}
\end{array}\right)=\left(\begin{array}{ccc}
\left|a\right|^{2}+\left|b\right|^{2} &  & \left(t+1\right)\overline{a}b\\
\\
\left(t+1\right)a\overline{b} &  & t^{2}\left|b\right|^{2}+\left|a\right|^{2}
\end{array}\right),
\]
and\hfill{}(5.4.7)
\[
TT^{*}=\left(\begin{array}{cc}
a & tb\\
\overline{b} & \overline{a}
\end{array}\right)\left(\begin{array}{cc}
\overline{a} & b\\
t\overline{b} & a
\end{array}\right)=\left(\begin{array}{ccc}
\left|a\right|^{2}+t^{2}\left|b\right|^{2} &  & \left(t+1\right)ab\\
\\
\left(t+1\right)\overline{ab} &  & \left|b\right|^{2}+\left|a\right|^{2}
\end{array}\right),
\]
in $\mathfrak{H}_{2}^{t}$. So, the realization $T$ of $\left(a,b\right)$
is normal in $\mathfrak{H}_{2}^{t}$, if and only if

\medskip{}

\hfill{}$\left|a\right|^{2}+t^{2}\left|b\right|^{2}=\left|a\right|^{2}+\left|b\right|^{2}\;\mathrm{and\;}\left(t+1\right)\overline{a}b=\left(t+1\right)ab,$\hfill{}(5.4.8)

\medskip{}

\noindent in $\mathbb{C}$, by (5.4.7).
\begin{prop}
Let $\left(a,b\right)\in\mathbb{H}_{t}$ satisfy $\mathbf{AN\;5.4.1}$.
Then the realization $T\in\mathcal{H}_{2}^{t}$ is normal in $\mathfrak{H}_{2}^{t}$,
if and only if

\medskip{}

\hfill{}$t^{2}\left|b\right|^{2}=\left|b\right|^{2}\;\mathrm{and\;}\left(t+1\right)\overline{a}b=\left(t+1\right)ab,$\hfill{}(5.4.9)

\medskip{}

\noindent in $\mathbb{C}$. In particular, if $t=-1$ (equivalently,
if $\left(a,b\right)\in\mathbb{H}_{-1}$ is a quaternion), then $T$
is normal in $\mathfrak{H}_{2}^{-1}$; if $t=1$, (equivalently, if
$\left(a,b\right)\in\mathbb{H}_{1}$ is a bicomplex number), then
$T$ is normal in $\mathfrak{H}_{2}^{1}$, if and only if

\medskip{}

\noindent \hfill{}either $\left(a,b\right)=\left(Re\left(a\right),b\right)$
or $\left(a,b\right)=\left(a,0\right)$ in $\mathbb{H}_{1}$;\hfill{}(5.4.10)

\medskip{}

\noindent meanwhile, if $t\in\mathbb{R}\setminus\left\{ \pm1\right\} $,
then $T$ is normal in $\mathfrak{H}_{2}^{t}$, if and only if

\medskip{}

\hfill{}$b=0\;\mathrm{in\;}\mathbb{C}\Longleftrightarrow\left(a,b\right)=\left(a,0\right)\in\mathbb{H}_{t}.$\hfill{}(5.4.11)
\end{prop}

\begin{proof}
By (5.4.8), the normality characterization (5.4.9) holds.

By (5.4.9), if $t=-1$ in $\mathbb{R}$, and hence, if $\left(a,b\right)\in\mathbb{H}_{-1}$
is a quaternion, then the condition (5.4.9) is identified with
\[
\left|b\right|^{2}=\left|b\right|^{2},\;\mathrm{and\;}0=0,
\]
which are the identities on $\mathbb{C}$. These identities demonstrate
that the realization of every quaternion is automatically normal in
$\mathfrak{H}_{2}^{-1}$.

Suppose $t=1$ in $\mathbb{R}$. Then the condition (5.4.9) is equivalent
to
\[
\left|b\right|^{2}=\left|b\right|^{2}\;\mathrm{and\;}2\overline{a}b=2ab,
\]
if and only if either
\[
\overline{a}=a\;\mathrm{in\;}\mathbb{C}\Longleftrightarrow\left(a,b\right)=\left(Re\left(a\right),b\right)\in\mathbb{H}_{1}\;\mathrm{(if}\;b\neq0\mathrm{),}
\]
or
\[
\left(a,b\right)=\left(a,0\right)\in\mathbb{H}_{1}\;\;\;\mathrm{(if\;}b=0\mathrm{)}.
\]
Thus, if $t=1$, then $T$ is normal, if and only if the condition
(5.4.10) holds.

Assume now that both $t\neq1$ and $t\neq-1$, i.e., suppose $t^{2}\neq1$
in $\mathbb{R}$. So, the first condition of (5.4.9) is identified
with
\[
t^{2}\left|b\right|^{2}=\left|b\right|^{2}\Longleftrightarrow b=0\;\mathrm{in\;}\mathbb{C}.
\]
So, the second condition of (5.4.9) automatically holds, since
\[
\left(t+1\right)\overline{a}\cdot0=\left(t+1\right)a\cdot0\Longleftrightarrow0=0.
\]
Therefore, the realization $T\in\mathcal{H}_{2}^{t}$ of $\left(a,b\right)\in\mathbb{H}_{t}$
is normal in $\mathfrak{H}_{2}^{t}$, if and only if $\left(a,b\right)=\left(a,0\right)$
in $\mathbb{H}_{t}$, whenever $t\in\mathbb{R}\setminus\left\{ \pm1\right\} $.
i.e., the normality (5.4.11) holds.
\end{proof}
The above proposition illustrates that: (i) the realizations of ``all''
quaternions are normal in $\mathfrak{H}_{2}^{-1}$, (ii) the realizations
of bicomplex numbers are normal in $\mathfrak{H}_{2}^{1}$, if and
only if either $\left(a,b\right)=\left(Re\left(a\right),b\right)$,
or $\left(a,b\right)=\left(a,0\right)$ in $\mathbb{H}_{1}$, by (5.4.10),
and (iii) the only realizations $\left[\left(a,0\right)\right]_{t}$
are normal in $\mathfrak{H}_{2}^{t}$, whenever $t\in\mathbb{R}\setminus\left\{ \pm1\right\} $,
by (5.4.11).
\begin{thm}
Let $\left(a,b\right)\in\mathbb{H}_{t}$ satisfy $\mathbf{AN\;5.4.1}$.

\noindent (5.4.12) Suppose $t=-1$. Then $T$ is normal in $\mathfrak{H}_{2}^{-1}$,
and its free distribution is characterized by the formula (5.3.10).

\noindent (5.4.13) Let $t\in\mathbb{R}\setminus\left\{ \pm1\right\} $.
If $T$ is ``non-zero'' normal in $\mathfrak{H}_{2}^{t}$, then
\[
\tau\left(\overset{n}{\underset{l=1}{\prod}}T^{r_{l}}\right)=R^{n}Re\left(W_{o}^{\overset{n}{\underset{l=1}{\sum}}e_{l}}\right),
\]
with\hfill{}(5.4.14)
\[
R=\left|a\right|\;\mathrm{and\;}W_{o}=\frac{a}{\left|a\right|}\in\mathbb{T},
\]
where
\[
e_{l}=\left\{ \begin{array}{ccc}
1 &  & \mathrm{if\;}r_{l}=1\\
-1 &  & \mathrm{if\;}r_{l}=*,
\end{array}\right.
\]
for $l=1,...,n$, for all $\left(r_{1},...,r_{n}\right)\in\left\{ 1,*\right\} ^{n}$,
for all $n\in\mathbb{N}$.
\end{thm}

\begin{proof}
The statement (5.4.12) holds by (5.3.10).

Of course, if $t<0$, and if $T\in\mathcal{H}_{2}^{t}$, then the
free-distributional data (5.4.14) holds by (5.3.10), because $T$
and the $t$-spectral form $S$ are similar in $\mathcal{H}_{2}^{t}$
as elements of $\left(\mathfrak{H}_{2}^{t},\tau\right)$. However,
in the statement (5.4.13), the normality works for all the scales
$t\in\mathbb{R}\setminus\left\{ \pm1\right\} $. Assume that the realization
$T$ is a ``non-zero,'' ``normal'' element of $\mathfrak{H}_{2}^{t}$.
Then
\[
\left(a,b\right)=\left(a,0\right)\in\mathbb{H}_{t},\;\mathrm{with\;}a\neq0,
\]
by (5.4.11). Therefore,
\[
T=\left(\begin{array}{cc}
a & 0\\
0 & \overline{a}
\end{array}\right)=S,
\]
because $\sigma_{t}\left(a,0\right)=a$ in $\mathbb{C}$. i.e., the
realization $T$ and the $t$-spectral form $S$ are identical in
$\mathfrak{H}_{2}^{t}$, implying the similarity of them. So, under
$\mathbf{AN\;5.4.1}$,
\[
a=w\overset{\textrm{denote}}{=}\sigma_{t}\left(a,0\right),
\]
polar-decomposed to be
\[
w=a=\left|a\right|\left(\frac{a}{\left|a\right|}\right)\in\mathbb{C},
\]
i.e., $r=\left|a\right|$ and $w_{o}=\frac{a}{\left|a\right|}$ under
$\mathbf{AN\;5.4.1}$. Therefore, similar to (5.3.10), the free-distributional
data (5.4.14) holds.
\end{proof}
Note that, in the proof of the statement (5.4.13), we did not use
our negative-scale assumption for the cases where $t<0$, but $t\neq-1$.
Indeed, even though $t\geq0$, but $t\neq1$, the normality (5.4.11)
shows that the realization $T$ is a diagonal matrix not affected
by the scale $t$. So, whatever scales $t$ are given in $\mathbb{R}\setminus\left\{ \pm1\right\} $,
the free-distributional data (5.4.14) holds in $\left(\mathfrak{H}_{2}^{t},\tau\right)$,
under normality. Then, how about the case where $t=1$? Recall that
if $t=1$, then the realization $T$ of $\left(a,b\right)\in\mathbb{H}_{1}$
is normal in $\mathfrak{H}_{2}^{1}$, if and only if either
\[
\left(a,b\right)=\left(Re\left(a\right),b\right),\;\mathrm{if\;}b\neq0,
\]
or
\[
\left(a,b\right)=\left(a,0\right),\;\mathrm{if\;}b=0,
\]
in $\mathbb{H}_{1}$, by (5.4.10). So, if $\left(a,b\right)=\left(a,0\right)$
in $\mathbb{H}_{1}$, the joint free moments of $T$ are determined
similarly by the formula (5.4.14), by the identity (and hence, the
similarity) of $T$ and $S$ (under $\mathbf{AN\;5.4.1}$). However,
if $\left(a,b\right)=\left(Re\left(a\right),b\right)$ with $b\neq0$,
then we need a better tool than (5.2.10) to compute the corresponding
free-distributional data, because we cannot use our similarity technique
(of Theorem 48) here.

By the definition of the unitarity, if an element $U$ of a $C^{*}$-algebra
$\mathcal{A}$ is a unitary, then it is automatically normal. i.e.,
the unitarity implies the normality. Let $\left(a,b\right)\in\mathbb{H}_{t}$
satisfy $\mathbf{AN\;5.4.1}$ with its realization $T\in\mathcal{H}_{2}^{t}$
in $\left(\mathfrak{H}_{2}^{t},\tau\right)$, and suppose it is a
unitary in $\mathfrak{H}_{2}^{t}$. By the assumption that $T$ is
a unitary in $\mathfrak{H}_{2}^{t}$, it is normal. 

Assume first that $t=-1$ in $\mathbb{R}$, and hence, $\left(a,b\right)\in\mathbb{H}_{-1}$
is a quaternion. Then the realization $T$ is automatically normal
in $\mathfrak{H}_{2}^{t}$ by (5.4.12). Indeed, in this case,
\[
T=\left(\begin{array}{cc}
a & -b\\
\overline{b} & \overline{a}
\end{array}\right)\;\mathrm{with\;}T^{*}=\left(\begin{array}{cc}
\overline{a} & b\\
-\overline{b} & a
\end{array}\right)=\left[\left(\overline{a},-b\right)\right]_{-1},
\]
in $\mathcal{H}_{2}^{-1}$, as elements of $\mathfrak{H}_{2}^{-1}.$
So, the normality is guaranteed;
\[
T^{*}T=\left(\begin{array}{ccc}
\left|a\right|^{2}+\left|b\right|^{2} &  & 0\\
\\
0 &  & \left|a\right|^{2}+\left|b\right|^{2}
\end{array}\right)=TT^{*},
\]
in $\mathcal{H}_{2}^{-1}$, as elements of $\mathfrak{H}_{2}^{-1}$.
It shows that $T$ is a unitary in $\mathfrak{H}_{2}^{-1}$, if and
only if

\medskip{}

\hfill{}$\left|a\right|^{2}+\left|b\right|^{2}=1.$\hfill{}(5.4.15)

\medskip{}

Meanwhile, if $t\in\mathbb{R}\setminus\left\{ \pm1\right\} $ in $\mathbb{R}$,
then $T$ is normal, if and only if $\left(a,b\right)=\left(a,0\right)$
in $\mathbb{H}_{t}$ by (5.4.11), if and only if
\[
T=\left(\begin{array}{cc}
a & 0\\
0 & \overline{a}
\end{array}\right)\in\mathcal{H}_{2}^{t},
\]
which is identical (and hence, similar) to the $t$-spectral form
$S$ of $\left(a,0\right)$ in $\mathfrak{H}_{2}^{t}$. This normal
element $T$ is a unitary in $\mathfrak{H}_{2}^{t}$, if and only
if
\[
T^{*}T=I_{2}=TT^{*}\Longleftrightarrow\left(\begin{array}{cc}
\left|a\right|^{2} & 0\\
0 & \left|a\right|^{2}
\end{array}\right)=\left(\begin{array}{cc}
1 & 0\\
0 & 1
\end{array}\right),
\]
if and only if\hfill{}(5.4.16)
\[
\left|a\right|^{2}=1\;\;\;\mathrm{in\;\;\;}\mathbb{C}.
\]

\begin{prop}
Let $\left(a,b\right)\in\mathbb{H}_{t}$ satisfy $\mathbf{AN\;5.4.1}$.

\noindent (5.4.17) Let $t=-1$. Then $T$ is a unitary in $\mathfrak{H}_{2}^{t}$,
if and only if $\left|a\right|^{2}+\left|b\right|^{2}=1$.

\noindent (5.4.18) Let $t\in\mathbb{R}\setminus\left\{ \pm1\right\} $.
Then $T$ is a unitary in $\mathfrak{H}_{2}^{t}$, if and only if
$\left|a\right|^{2}=1$ and $b=0$.
\end{prop}

\begin{proof}
The statements (5.4.17) and (5.4.18) hold by (5.4.15) and (5.4.16),
respectively, because a unitary realization $T$ of $\left(a,b\right)$
automatically satisfies the normality (5.4.9).
\end{proof}
\bigskip{}

\noindent $\mathbf{Observation.}$ Now, assume that $t=1$, and let
$\left(a,b\right)\in\mathbb{H}_{1}$ be a bicomplex number satisfying
$\mathbf{AN\;5.4.1}$. By (5.4.10), the realization $T\in\mathcal{H}_{2}^{1}$
is normal in $\mathfrak{H}_{2}^{1}$, if and only if either
\[
\left(a,b\right)=\left(Re\left(a\right),b\right),\;\mathrm{or\;}\left(a,b\right)=\left(a,0\right),
\]
in $\mathbb{H}_{1}$. So, if $\left(a,b\right)=\left(a,0\right)$
in $\mathbb{H}_{1}$, then one obtains the unitarity that: $T$ is
a unitary in $\mathfrak{H}_{2}^{1}$, if and only if $\left|a\right|^{2}=1$,
just like (5.4.18). However, if 
\[
\left(a,b\right)=\left(Re\left(a\right),b\right)=\left(x,b\right)\;\mathrm{in\;}\mathbb{H}_{1},
\]
with $b\neq0$ in $\mathbb{C}$, then $T$ is a unitary in $\mathfrak{H}_{2}^{1}$,
if and only if
\[
\left(\begin{array}{cc}
x & \overline{b}\\
b & x
\end{array}\right)\left(\begin{array}{cc}
x & b\\
\overline{b} & x
\end{array}\right)=\left(\begin{array}{ccc}
x^{2}+\overline{b^{2}} &  & 2xRe\left(b\right)\\
\\
2xRe\left(b\right) &  & x^{2}+b^{2}
\end{array}\right)=I_{2},
\]
and
\[
\left(\begin{array}{cc}
x & b\\
\overline{b} & x
\end{array}\right)\left(\begin{array}{cc}
x & \overline{b}\\
b & x
\end{array}\right)=\left(\begin{array}{ccc}
x^{2}+b^{2} &  & 2xRe\left(b\right)\\
\\
2xRe\left(b\right) &  & x^{2}+\overline{b^{2}}
\end{array}\right)=I_{2},
\]
in $\mathfrak{H}_{2}^{1}$, if and only if
\[
x^{2}+\overline{b^{2}}=x^{2}+b^{2}=1\;\mathrm{and\;}2xRe\left(b\right)=0,
\]
if and only if
\[
b^{2}=\overline{b^{2}}=1-x^{2}\;\mathrm{and\;}2xRe\left(b\right)=0,
\]
if and only if
\[
b^{2}=1-x^{2}\in\mathbb{R}\;\mathrm{and\;}x=0,
\]
because $b$ is assumed not to be zero in $\mathbb{C}$, if and only
if
\[
x=0\;\mathrm{and\;}b=\pm1\;\;\;\mathrm{in\;\;\;}\mathbb{R},
\]
if and only if
\[
T=\left(\begin{array}{cc}
0 & 1\\
1 & 0
\end{array}\right),\;\mathrm{or\;}T=\left(\begin{array}{cc}
0 & -1\\
-1 & 0
\end{array}\right)\;\mathrm{in\;}\mathcal{H}_{2}^{1},
\]
if and only if
\[
\left(a,b\right)=\left(0,1\right),\;\mathrm{or\;}\left(a,b\right)=\left(0,-1\right)\;\mathrm{in\;}\mathbb{H}_{1}.
\]
i.e., if $\left(a,b\right)=\left(Re\left(a\right),b\right)$ in $\mathbb{H}_{1}$,
then $T$ is a unitary in $\mathfrak{H}_{2}^{1}$, if and only if
\[
\left(a,b\right)=\left(0,1\right),\;\mathrm{or\;}\left(a,b\right)=\left(0,-1\right),
\]
in $\mathbb{H}_{1}$. In summary, the realization $T\in\mathcal{H}_{2}^{1}$
of a bicomplex number $\left(a,b\right)\in\mathbb{H}_{1}$ is a unitary
in $\mathfrak{H}_{2}^{t}$, if and only if either
\[
\left(a,b\right)=\left(a,0\right)\;\;\mathrm{with\;\;}\left|a\right|^{2}=1,
\]
or
\[
\left(a,b\right)=\left(0,1\right),\;\mathrm{or\;}\left(a,b\right)=\left(0,-1\right),
\]
in $\mathbb{H}_{1}$.\hfill{}\textifsymbol[ifgeo]{64}

\bigskip{}

By the unitarity (5.4.17) and (5.4.18), one has the following result. 
\begin{thm}
Let $\left(a,b\right)\in\mathbb{H}_{t}$ satisfy $\mathbf{AN\;5.4.1}$.

\noindent (5.4.19) Suppose $t=-1$. If $T$ is a unitary in $\mathfrak{H}_{2}^{t}$,
then its free distribution is characterized by the formula (5.3.10)
with $r=1$.

\noindent (5.4.20) Let $t\in\mathbb{R}\setminus\left\{ \pm1\right\} $.
If $T$ is a unitary in $\mathfrak{H}_{2}^{t}$, then
\[
\tau\left(\overset{n}{\underset{l=1}{\prod}}T^{r_{l}}\right)=Re\left(a^{\overset{n}{\underset{l=1}{\sum}}e_{l}}\right),\;\mathrm{with\;}a\in\mathbb{T}\;\mathrm{in\;}\mathbb{C},
\]
where\hfill{}(5.4.21)
\[
e_{l}=\left\{ \begin{array}{ccc}
1 &  & \mathrm{if\;}r_{l}=1\\
-1 &  & \mathrm{if\;}r_{l}=*,
\end{array}\right.
\]
for $l=1,...,n$, for all $\left(r_{1},...,r_{n}\right)\in\left\{ 1,*\right\} ^{n}$,
for all $n\in\mathbb{N}$.
\end{thm}

\begin{proof}
The statement (5.4.19) holds by (5.3.11). In particular, by the unitarity
characterization (5.4.17), the free-distributional data in (5.3.11)
must have $r=1$, since
\[
\left|\sigma_{t}\left(a,b\right)\right|=\left|w\right|\overset{\textrm{denote}}{=}r=1,
\]
under the similarity of $T$ and $S$, by (5.4.17). 

By (5.4.13), if $t\neq\pm1$, then the free-distributional data (5.4.21)
holds by (5.4.14). Indeed, under the unitarity of $\mathit{T}$, the
formula (5.4.14) satisfies
\[
R=\left|a\right|=1\;\mathrm{and\;}W_{o}=a\in\mathbb{T}.
\]
Therefore, the joint free moments (5.4.21) holds.
\end{proof}
The above theorem characterizes the free distributions of unitary
elements of $\left(\mathfrak{H}_{2}^{t},\tau\right)$ induced by $\mathbb{H}_{t}$,
where $t\in\mathbb{R}\setminus\left\{ 1\right\} $.

Suppose $t=1$, and $\left(a,b\right)\in\mathbb{H}_{1}$ satisfies
$\mathbf{AN\;5.4.1}$. In the above $\mathbf{Observation}$, we showed
that the realization $T\in\mathcal{H}_{2}^{1}$ of $\left(a,b\right)$
is a unitary, if and only if either
\[
\left(a,b\right)=\left(a,0\right)\;\mathrm{with\;}a\in\mathbb{T},
\]
or
\[
\left(a,b\right)=\left(0,1\right),\;\mathrm{or\;}\left(a,b\right)=\left(0,-1\right),
\]
in $\mathbb{H}_{1}$, equivalently, either
\[
T=\left(\begin{array}{cc}
a & 0\\
0 & \overline{a}
\end{array}\right)\;\;\mathrm{with\;\;}a\in\mathbb{T},
\]
or
\[
T=\left(\begin{array}{cc}
0 & 1\\
1 & 0
\end{array}\right),\;\mathrm{or\;}T=\left(\begin{array}{cc}
0 & -1\\
-1 & 0
\end{array}\right),
\]
in $\mathcal{H}_{2}^{1}$ (as an element of $\mathfrak{H}_{2}^{1}$).
Thus, if $\left(a,b\right)=\left(a,0\right)\in\mathbb{H}_{1}$ with
$\left|a\right|^{2}=1$, then the free distribution of $T$ is similarly
characterized by the formula (5.4.21). Meanwhile, if $T=\left[\left(0,1\right)\right]_{1}$,
then
\[
T^{*}=T\in\mathcal{H}_{2}^{1}\subset\mathcal{H}_{2}^{1}\left(1,*\right)\;\mathrm{in\;}\mathfrak{H}_{2}^{1},
\]
and
\[
T^{2}=\left(\begin{array}{cc}
0 & 1\\
1 & 0
\end{array}\right)\left(\begin{array}{cc}
0 & 1\\
1 & 0
\end{array}\right)=\left(\begin{array}{cc}
1 & 0\\
0 & 1
\end{array}\right)=I_{2},
\]
in $\mathfrak{H}_{2}^{1}$, satisfying that\hfill{}(5.4.22)
\[
\left(T^{n}\right)_{n=1}^{\infty}=\left(T,I_{2},T,I_{2},T,I_{2},...\right);
\]
and, if $T=\left[\left(0,-1\right)\right]_{1}$, then
\[
T^{*}=T\in\mathcal{H}_{2}^{1}\subset\mathcal{H}_{2}^{1}\left(1,*\right)\;\mathrm{in\;}\mathfrak{H}_{2}^{1},
\]
and
\[
T^{2}=\left(\begin{array}{cc}
0 & -1\\
-1 & 0
\end{array}\right)\left(\begin{array}{cc}
0 & -1\\
-1 & 0
\end{array}\right)=\left(\begin{array}{cc}
1 & 0\\
0 & 1
\end{array}\right)=I_{2},
\]
in $\mathfrak{H}_{2}^{1}$, satisfying that\hfill{}(5.4.23)
\[
\left(T^{n}\right)_{n=1}^{\infty}=\left(T,I_{2},T,I_{2}T,I_{2},...\right).
\]
Therefore, one obtains the following result in addition with Theorem
59.
\begin{thm}
Let $\left(a,b\right)\in\mathbb{H}_{1}$ be a bicomplex number satisfying
$\mathbf{AN\;5.4.1}$. Then the realization $T$ is a unitary in $\left(\mathfrak{H}_{2}^{1},\tau\right)$,
if and only if either
\[
\left(a,b\right)=\left(a,0\right),\;\mathrm{with\;}\left|a\right|^{2}=1,
\]
or\hfill{}(5.4.24)
\[
\left(a,b\right)=\left(0,1\right),\;\mathrm{or\;}\left(a,b\right)=\left(0,-1\right)\;\mathrm{in\;}\mathbb{H}_{1}.
\]

\noindent (5.4.25) If $\left(a,b\right)=\left(a,0\right)$, with $\left|a\right|^{2}=1$,
in $\mathbb{H}_{1}$, then the free distribution of $T$ is characterized
by the formula (5.4.21).

\noindent (5.4.26) If either $\left(a,b\right)=\left(0,1\right)$,
or $\left(a,b\right)=\left(0,-1\right)$ in $\mathbb{H}_{1}$, then
the free distribution of the unitary realization $T$ is fully characterized
by the free-moment sequence, 

\medskip{}

\hfill{}$\left(\tau\left(T^{n}\right)\right)_{n=1}^{\infty}=\left(0,1,0,1,0,1,0,1,...\right).$\hfill{}(5.4.27)
\end{thm}

\begin{proof}
By the very above $\mathbf{Observation}$ after Proposition 58, it
is shown that the realization $T\in\mathcal{H}_{2}^{1}$ of a bicomplex
number $\left(a,b\right)\in\mathbb{H}_{1}$ is a unitary in $\mathfrak{H}_{2}^{1}$,
if and only if the condition (5.4.24) holds true.

The statement (5.4.25) is shown similarly by the proof of the statement
(5.4.20). So, the free-distributional data (5.4.21) holds.

Now, if either $T=\left[\left(0,1\right)\right]_{1}$, or $T=\left[\left(0,-1\right)\right]_{1}$
in $\mathcal{H}_{2}^{1}$, it is not only a unitary, but also a self-adjoint
element of $\left(\mathfrak{H}_{2}^{1},\tau\right)$, and hence, the
free distribution of $T$ is fully characterized by the free-moment
sequence $\left(\tau\left(T^{n}\right)\right)_{n=1}^{\infty}$. However,
by (5.4.22) and (5.4.23), one immediately obtain the free-moment sequence
(5.4.27). Therefore, the statement (5.4.26) holds.
\end{proof}
The above theorem fully characterizes the free distributions of the
unitaries of $\left(\mathfrak{H}_{2}^{1},\tau\right)$ induced by
bicomplex numbers of $\mathbb{H}_{1}$.

\end{document}